\documentclass{elsarticle}

\usepackage[a4paper,width=155mm,top=25mm,bottom=30mm]{geometry}
\usepackage{subcaption}
\usepackage{mdframed}   
\usepackage{booktabs}
\usepackage{afterpage}
\usepackage{contour}
\usepackage{bm}
\usepackage{stmaryrd}
\usepackage{amsthm}
\usepackage{calrsfs}
\DeclareMathAlphabet{\pazocal}{OMS}{zplm}{m}{n}
\usepackage{mathtools}
\usepackage{framed} 
\usepackage{capt-of}
\usepackage[utf8x]{inputenc} 

\newtheorem{proposition}{Proposition}
\usepackage{graphicx,import}

\usepackage{float}
\floatstyle{plaintop}
\restylefloat{table}

\contourlength{0.4pt}
\contournumber{64}

\theoremstyle{remark}
\newtheorem{remark}{Remark}

\usepackage{esvect} 

\usepackage{tikz}
\usepackage{pgfplots}
\usepgfplotslibrary{external}
\usetikzlibrary{spy}
\usetikzlibrary{positioning}
\usepackage{graphicx}
\usepackage{amssymb}
\usepackage{lipsum}
\usepackage{mathtools}
\usepackage{amsmath}
\usepackage{bbm}
\usepackage{commath} 
\mathtoolsset{showonlyrefs} 
\usepackage{amssymb}
\usepackage{siunitx} 
\mdfdefinestyle{box}{%
    outerlinewidth=0.5pt,
    innertopmargin=5pt,
    innerbottommargin=0pt,
    innerrightmargin=0pt,
    innerleftmargin=0pt,
    topline=false,
    rightline=false,
    leftline=false,
    bottomline=false,
    backgroundcolor=white}

\newcounter{weakforms}

\newenvironment{weakform}[1]
{
\refstepcounter{weakforms}
\begin{mdframed}[style=box]
\textbf{(W\theweakforms)}
\itshape
}
{
    \end{mdframed}
    }

\usepackage{helvet}
\usepackage[EULERGREEK]{sansmath}
\definecolor{mypink}{RGB}{255,107,252}
\definecolor{myorange}{RGB}{255,116,17}
\definecolor{myorange2}{RGB}{255,174,61}
\definecolor{myyellow}{RGB}{255,221,2}
\definecolor{mylightblue}{RGB}{147,230,255}
\definecolor{mylightred}{RGB}{255,147,147}
\definecolor{mygreen}{RGB}{176,255,160}
\pgfplotscreateplotcyclelist{customlist}{%
	{blue,every mark/.append style={fill=blue!80!black},mark=*},
	{red,every mark/.append style={fill=red!80!black},mark=square*},
    {green!40!black,every mark/.append style={fill=mygreen},mark=otimes*},
	{blue!70!red,every mark/.append style={fill=mypink},mark=triangle*,mark size=2.5pt},
	{myorange,every mark/.append style={fill=myyellow},mark=diamond*,mark size=3.0pt},
	{blue,every mark/.append style={fill=mylightblue},mark=pentagon*,mark size=2.8pt},
	{red,every mark/.append style={fill=white},mark=oplus*},
	{myorange,every mark/.append style={fill=myorange2,rotate=180},mark=triangle*,mark size=2.5pt}}

\pgfplotsset{
	compat=newest,
	tick label style={font=\sansmath\sffamily\small},
    label style={font=\sansmath\sffamily\small},
	legend style={nodes={scale=0.7, transform shape},draw=none, font=\sansmath\sffamily, align=right},
	cycle list name=customlist
	}

\usepackage[colorlinks = true,
            linkcolor = black,
            urlcolor  = black,
            citecolor = black,
            anchorcolor = black]{hyperref}

\newcommand{\mat}[1]{#1}
\newcommand{\uvec}{\mathbf u} 
\newcommand{\p}{p} 
\newcommand{\fvec}{\mathbf f} 
\newcommand{\hvec}{\mathbf h}
\newcommand{\vvec}{\mathbf v}

\newcommand{\uvech}{\mathbf u^h}
\newcommand{\xvec}{\mathbf x}
\newcommand{\ph}{\p^h}
\newcommand{\gvec}{\mathbf g}
\newcommand{\nnorm}{\mathbf n} 
\newcommand{\domain}{\Omega}
\newcommand{\dirboundary}{\Gamma_D}
\newcommand{\neuboundary}{\Gamma_N}
\newcommand{\reals}{\mathbb{R}}
\newcommand{\realsd}{\mathbb{R}^{d}}
\newcommand{\algvec}[1]{\underline{\bm{\mathsf{#1}}}}
\newcommand{\nuh}{N_u^h}
\newcommand{\nph}{N_p^h}

\newcommand{\bbref}{\widehat{\domain}}
\newcommand{\paramspace}{\mathcal M}
\newcommand{\param}{\algvec{\mu}}
\newcommand{\is}[1]{{#1}} 
\newcommand{\irs}[1]{{#1}} 
\newcommand{\iin}[1]{{[#1]}} 

\newcommand{\weakref}[1]
{\hyperref[#1]{\textrm{W\ref*{#1}}}}

\journal{arXiv}

\begin{document}

\newlength
\figureheight
\newlength
\figurewidth

\title{Model order reduction of flow based on a modular geometrical approximation of blood vessels}

\author[epfl]{Luca Pegolotti\corref{cor1}}
\ead{luca.pegolotti@epfl.ch}
\author[stanford]{Martin Pfaller}
\author[stanford]{Alison Marsden}
\author[epfl]{Simone Deparis}

\cortext[cor1]{Corresponding author}

\address[epfl]{SCI-SB-SD, Institute of Mathematics, \'{E}cole Polytechnique
F\'{e}d\'{e}rale de Lausanne, Station 8, EPFL, CH–1015 Lausanne, Switzerland}
\address[stanford]{Department of Pediatrics (Cardiology), Bioengineering, Stanford University, Clark Center E1.3, 318 Campus Drive, Stanford, CA 94305, USA}

\begin{abstract}
We are interested in a reduced order method for the efficient simulation of blood flow in arteries. The blood dynamics is modeled by means of the incompressible Navier--Stokes equations. Our algorithm is based on an approximated domain-decomposition of the target geometry into a number of subdomains obtained from the parametrized deformation of geometrical building blocks (e.g. straight tubes and model bifurcations). On each of these building blocks, we build a set of spectral functions by proper orthogonal decomposition of a large number of snapshots of finite element solutions (offline phase). The global solution of the Navier--Stokes equations on a target geometry is then found by coupling linear combinations of these local basis functions by means of spectral Lagrange multipliers (online phase). Being that the number of reduced degrees of freedom is considerably smaller than their finite element counterpart, this approach allows us to significantly decrease the size of the linear system to be solved in each iteration of the Newton--Raphson algorithm. We achieve large speedups with respect to the full order simulation (in our numerical experiments, the gain is at least of one order of magnitude and grows inversely with respect to the reduced basis size), whilst still retaining satisfactory accuracy for most cardiovascular simulations.
\end{abstract}

\maketitle

\section{Introduction}
Cardiovascular disease is the leading cause of death worldwide. This term broadly encompasses a variety of pathological cases ranging from heart disease to many other peripheral vascular diseases. Consequently, the numerical simulation of blood flow in the cardiovascular system has gained considerable attention during the last twenty years as a valuable quantitative tool for the study and diagnosis of such conditions \cite{bao2014usnctam,figueroa2017}.

Blood dynamics is typically modeled by means of the incompressible Navier--Stokes equations; their discretization by numerical methods such as the finite element (FE) method leads to the \textit{Full Order Model} (FOM). Despite the rapid and constant growth in computational power of the architectures that are currently employed to run such simulations, the time and resources required are often incompatible with clinical practice.  Moreover, in many cases the numerical results of the FOM are affected by a level of uncertainty not inherently associated with the employed numerical method but rather with the inexact geometry and boundary conditions considered in the simulation \cite{sankaran2016uncertainty,sankaran2011stochastic}. \textit{Reduced Order Models} (ROMs) aim at lowering the computational burden of FOMs at the cost of settling for larger approximation errors. This is particularly desirable in multi-query scenarios, i.e. whenever the same simulation needs to be performed for multiple values of the input parameters (for instance, in order to quantify the uncertainty due to the problem data \cite{chen2013simulation}).

Among the ROMs employed in the context of cardiovascular modeling are the popular geometrical multiscale 0D/1D models \cite{malossi2013numerical,malossi2012partitioned,marsden2015multiscale}. These models consider a coarse approximation of systems of arteries as electric circuits (0D) or as segments in which the quantities of interest are obtained by an averaging process across the section of the vessels (1D). Geometrical multiscale models often prove to be remarkably accurate in approximating flow rates and pressure drops \cite{blanco2018comparison}, albeit the strong geometrical approximations inevitably entail a significant loss of local details. For this reason, algorithms to couple geometrical multiscale models with full 3D simulations---to be employed in the regions in which higher quality solutions are required---have been devised \cite{malossi2013implicit,sankaran2012patient,moghadam2013modular}.

In this paper, we aim at formalizing a ROM allowing to approximate the local features of the blood flow. The strategy is based on the combination of a domain-decomposition approach with the Reduced Basis (RB) method and can be interpreted as a specific implementation of the Reduced Basis Element (RBE) method \cite{maday2004reduced,Maday2002195}. We refer the reader to \cite{iapichino2012reduced,iapichino2010reduced,iapichino2014reduced,lovgren2006reduced3,lovgren2006reduced2} for uses of the RBE method for the approximation of the 2D steady Stokes equations in the context of the cardiovascular system. To our knowledge, our work represents the first application of the RBE method to the unsteady 3D Navier--Stokes equations.

Similarly to 1D models, the proposed method is based on geometrical approximation of the vessel geometry as a composition of simple subdomains. However, in our approach, these subdomains are three dimensional and obtained from the parametrized geometrical deformation of a handful of elementary building blocks (e.g. straight tubes and bifurcations). Each building block in its reference configuration is equipped with a set of spectral basis functions. Since in this paper we focus on the FE method to generate the FOMs, these are FE functions defined over triangulations of the building block. Specifically, the basis functions are found by means of Proper Orthogonal Decomposition (POD) of a large number of flow solutions which are computed during a computationally expensive offline data generation phase. The global flow approximations by our ROM is computed as composition of local (to the subdomains) solutions, namely linear combinations of the basis functions defined in every subdomain scaled by the divergence-free Piola transformation. The local solutions are coupled by a nonconforming domain-decomposition method based on the use of spectral Lagrange multipliers on the 2D interfaces \cite{deparis2019coupling}.

The gain in performance with respect to the FOM is given by the decreased size of the linear system to be solved at each iteration of the Newton--Raphson algorithm. Indeed, while the number of degrees of freedom in the FE model is typically large (tens or hundreds of thousands per each subdomain in our numerical simulations), only a few hundred basis functions per subdomain are sufficient to retain acceptable levels of accuracy.

The rest of this paper is structured as follows. In Section~\ref{sec:rb}, we provide a self-contained and concise introduction to the RB method in order to set the notation and terminology for the remaining sections. Section~\ref{sec:ns} is dedicated to the Navier--Stokes equations and their numerical discretization by the FE method. In Section~\ref{sec:modular}, we define the concept of modular approximation of arteries by the above-mentioned domain-decomposition approach; we also address the numerical solution of the Navier--Stokes equations on the decomposed geometries by the FE method and the nonconforming domain-decomposition method. In particular, we devise an ad-hoc preconditioner that takes advantage of the peculiar block structure of the global system matrix. It is worth noticing that, although the paper focuses on a ROM, addressing the solution of the partitioned problem with the FE method is necessary, as the RB functions in the subdomains are generated by POD of local solutions obtained from global problems in decomposed domains. This strategy of data collection (offline phase) is discussed in Section~\ref{sec:nsmodularrb}. In the same section, we also delineate the algorithm for the approximation of the global solution on a decomposed target geometry using the ROM (online phase). Our numerical results are reported in Section~\ref{sec:numericalresults}. Specifically, Section~\ref{subsec:numres_artificial} and Section~\ref{subsec:numres_aorta}, respectively, focus on two critical points: (i) the accuracy with respect to the corresponding FE solution in the decomposed geometry---which, in our case, plays the role of FOM---and the achieved speedup, and (ii) the comparison of the reduced solution to the one obtained on a physiological and non-decomposed geometry. In the latter case, we mostly aim at evaluating the effects of the geometrical approximation on the local features of the flow, for instance, in terms of the wall-shear stress (WSS).  Finally, in Section~\ref{sec:conclusions} we draw our conclusions and discuss future perspectives of the current study.

\subsection{Notation}
\paragraph{Vectors and matrices} We use the notation $\mathbf a:\domain \mapsto \reals^d$ to indicate $d$-dimensional vector fields defined over some domain $\domain \subset \reals^d$, whereas $\algvec{a} \in \reals^{m}$ refers to an algebraic column vector. When considering multiple subdomains, local matrices and global block matrices are written as $\mat{A}$ and $\pazocal A$, respectively. Global block algebraic vectors are written in capital letters and we use the notation $\algvec{A} = [\algvec{a}_1,\ldots,\algvec{a}_n]\in \reals^M$, $M = \sum_{i = 1}^{n} m_i$, to indicate the concatenation of $\algvec{a}_1 \in \reals^{m_1},\ldots,\algvec{a}_n\in \reals^{m_n}$. On the other hand, the notation $\mat{A} = [\algvec{a}_1|\ldots|\algvec{a}_n] \in \reals^{m \times n}$ indicates a matrix whose colums are represented by $\algvec{a}_1 \in \reals^{m},\ldots,\algvec{a}_n \in \reals^{m}$.

\paragraph{Superscripts, subscripts and hats} We use the superscripts $h$ or $N$ to denote quantities related to the FE or the RB approximations. In Section~\ref{sec:modular}, we introduce a modular domain-decomposition method based on reference building blocks and subdomains of the target geometry. Indices of the subdomains are written as simple superscripts (e.g. $\domain^\is{j}$); superscripts in square brackets refer to interfaces (e.g. $\Gamma^\iin{ij} = \overline{\domain}^\is{i} \cap \overline{\domain}^\is{j}$) or portions of boundaries and related quantities. In the context of the Newton--Raphson algorithm, superscripts in pharentheses are used to refer to the iteration index of the iterative method.
Subscripts are used to indicate: (i) indices of vectors or components (or blocks) of vectors or matrices, and (ii) quantities at specific timesteps (e.g. $\algvec{u}^h_i$ is the vector of FE degrees of freedom of velocity at the $i^\text{th}$ timestep).
Whenever a symbol refers to a reference building block or the reference interface (unit disk), it is indicated with a hat notation (e.g. $\bbref^{\is{i}}$ is the $i^\text{th}$ reference building block).

\section{The Reduced Basis method in a nutshell}
\label{sec:rb}
In this section, we provide a non-comprehensive introduction to the RB method which is intended to set the theoretical basis for the remainder of the paper. For a more complete overview, we refer the reader to \cite{quarteroni2015reduced,hesthaven2015certified}.

Let us consider an open and bounded domain $\domain$ and a steady differential problem of the form
\begin{equation}
    \mathcal{L}(u;\param) = \mathcal{G}(\param),
    \label{eq:differential_elliptic}
\end{equation}
where $u \in \pazocal V$ ($\pazocal V$ being a suitable functional space) is the solution, $\param \in \pazocal D \subset \reals^{N_\mu}$ is a vector of geometrical and/or physical parameters, $\mathcal L$ is a generic differential operator, and $\mathcal G$ is a functional encoding the data of the problem, such as forcing term and boundary conditions. In this section we assume that $\mathcal L$ is an elliptic operator. The extension of this setting to the Navier--Stokes equations is considered in Section~\ref{sec:nsmodularrb}.

The standard approach to solve Eq.~\eqref{eq:differential_elliptic} by a Galerkin method corresponds to transforming the continuous problem into a finite dimensional one, often referred to as \textit{Full Order Model} (FOM). For instance, in Section~\ref{subsec:numdiscns} we discuss how the numerical solution of the Navier--Stokes equations by the FE method is found as a linear combination of the FE basis functions. With respect to the model problem in Eq.~\eqref{eq:differential_elliptic}, this translates to $u^h = \sum_{i = 1}^{N^h} u^h_i \varphi^h_i$, where $\varphi^h_i \in \pazocal V^h \subset \pazocal V$ are FE basis functions;  $\algvec{u}^h = [u^h_1,\ldots,u^h_{N^h}]$ is typically called vector of degrees of freedom. Assuming that the differential operator $\mathcal L$ can be mathematically described in the weak sense by a bilinear form as $\ell(\varphi, \psi)$, for $\varphi \in \pazocal V$ and $\psi \in \pazocal V$, we identify the matrix $\mat{L}^h(\param)_{ij} = \ell(\varphi^h_j,\varphi^h_i;\param) \in \reals^{N^h \times N^h}$, and similarly $\algvec{G}^h(\param)_i = \int_{\Omega} \mathcal G(\param) \varphi^h_i \in \reals^{N^h}$. For example, for the linear differential operator $\mathcal L(u;\eta) = \eta \Delta u$, describing a Poisson equation with parameter $\eta$, we have $\ell(\varphi^h_j,\varphi^h_i;\eta) = \int_{\domain} \eta \nabla \varphi^h_j \cdot \nabla \varphi^h_i$. The resulting linear system of dimension $N^h \times N^h$
\begin{equation}
    \mat{L}^h(\param) \algvec u^h = \algvec{G}^h(\param)
    \label{eq:elliptic_ls}
\end{equation}
is possibly very large and expensive to solve. In multi-query scenarios---i.e. whenever it is required to solve Eq.~\eqref{eq:elliptic_ls} for multiple values of the parameter $\param$---it is often crucial to reduce the dimensionality of the system in order to save computational time. One way to achieve this is by employing ROMs, such as the RB method.

The main idea of the RB method is to construct a low dimensional basis for the solution $u$ out of a number of solutions $N_s$ (snapshots) of the FOM, which are computed during the so-called \textit{offline phase}. In the \textit{online phase}, the reduced solution is obtained as a linear combination of the RB functions; system \eqref{eq:elliptic_ls} is casted in the form of a small linear system where the unknowns represent the coefficients of such a linear combination. In the remainder of this section, we address the offline and online phases more in depth.

\subsection{The offline phase: basis construction}
\label{subsec:rboffline}
There exist two main strategies for the construction of the reduced basis: greedy algorithms \cite{binev2011convergence,hesthaven2014efficient} and the proper orthogonal decomposition (POD) method. The former lead to a more efficient offline phase, as they allow to minimize the number of snapshots $N_s$ to be computed. A major drawback of greedy algorithms is that they are based on an a posteriori estimate of the projection error, which is often difficult to compute in practical applications. For this reason, in this paper we opt for the POD method, which often requires a larger number of snapshots $N_s$ but is in turn more general. We refer e.g. \cite{kunisch2001galerkin,rathinam2003new} and \cite{Kunisch2002GalerkinPO} for applications of POD to parabolic and fluid problems, respectively, and \cite{rathinam2003new} for  \cite{kunisch2001galerkin} for a comprehensive study of the properties of POD when applied to the solution of Ordinary Differential Equations (ODEs). In the context of cardiac simulations, this technique has been successfully employed both in fluid (e.g. in \cite{ballarin2016fast} to simulate blood flow in patient-specific coronary artery bypass grafts) and structural simulations (e.g. in \cite{pfaller2020using}, where POD is used to reduce the space of admissible displacements of the heart muscle). In the POD approach, the reduced basis is usually constructed by singular value decomposition (SVD) \cite{golub2012matrix,trefethen1997numerical} out of the set of snapshots, which are obtained by sampling $N_s$ parameters $\param_1,\ldots,\param_{N_s}$ in $\pazocal D$ and by solving the corresponding FOM. Formally, we arrange the snapshots in matrix form as $\mat{S} = [\algvec{u}^h_1|\ldots|\algvec{u}^h_{N_s}] \in \reals^{N^h \times N_s}$ and we seek matrices $\mat{U} = [\algvec{\zeta}^h_1|\ldots|\algvec{\zeta}^h_{N_s}] \in \reals^{N^h \times N_s}$, $\Sigma \in \reals^{N_s \times N_s}$ and $\mat{Z} \in \reals^{N_s \times N_s}$ such that $\mat{S} = \mat{U} \Sigma \mat{Z}^\text{T}$; the columns of $\mat{U}$ and $\mat{Z}$ are orthonormal. In the context of POD, $\algvec{\zeta}^h_1,\ldots,\algvec{\zeta}^h_{N_s}$ are often called modes. We remark that in classic SVD the matrices $\mat{U}$ and $\Sigma$ are of size $N^h \times N^h$ and $N^h \times N_s$, respectively; here we consider the ``economic'' version of the algorithm. Matrix $\Sigma$ takes the form $\Sigma = \text{diag}(\sigma_1,\ldots,\sigma_{N_s})$ and its diagonal is composed of the singular values of matrix $\mat{S}$ ordered from largest to smallest, i.e. $\sigma_1 \geq \ldots \geq \sigma_{N_s} \leq 0$. Let us define $\mat{V}:= [\algvec{\zeta}^h_1|\ldots|\algvec{\zeta}^h_{N}] \in \reals^{N^h \times N} $ as the matrix composed of the first $N$ modes and let us recall that, given a N-dimensional orthonormal basis $\mat{W} = [\algvec{w}_1|\ldots|\algvec{w}_N] \in \reals^{N^h \times N}$, the projection of a generic vector $\algvec{x} \in \reals^{N^h}$ onto span$\{\algvec{w}_1|\ldots|\algvec{w}_N\}$ is given by $\Pi_W \algvec{x} = \mat{W}\mat{W}^\text{T} \algvec{x}$. Then, the following proposition holds.
\begin{proposition}
Let $\mathcal V_N = \{ \mat{W} \in \reals^{N^h \times N}: \mat{W}^\text{T} \mat{W} = \mat{I} \}$ be the set of all N-dimensional orthonormal bases. Then,
\begin{equation}
    \sum_{i = 1}^{N_s}\Vert \algvec{u}_i - \mat{V}\mat{V}^\text{T} \algvec{u}_i \Vert_2^2 = \min_{\mat{W} \in \mathcal V_N} \sum_{i = 1}^{N_s} \Vert \algvec{u}_i - \mat{W} \mat{W}^\text{T} \algvec{u}_i \Vert_2^2 = \sum_{i = N+1}^{N_s} \sigma_i^2.
\end{equation}
\label{prop:optimality}
\end{proposition}
We refer to \cite{quarteroni2015reduced} for a proof of Proposition~\ref{prop:optimality}. In other words, $\mat{V}$ is the N-dimensional basis minimizing the projection error of the snapshots over its column space; moreover, such error is strictly related to the magnitude of the singular values $\sigma_{N+1},\ldots \sigma_{N_s}$. Thus, a common heuristic to choose $N$ is to set it equal to the smallest integer $\widetilde{N}$ such that
\begin{equation}
\dfrac{\sum_{i = 1}^{\widetilde{N}} \sigma_i^2}{\sum_{i = 1}^{N_s} \sigma_i^2} \geq 1 - \varepsilon^2,
\label{eq:energy}
\end{equation}
where $\varepsilon$ is a user-provided tolerance. The left hand side of Eq.~\eqref{eq:energy} is the relative information content of the POD basis, namely the percentage of energy of the snapshots retained by the first $\widetilde{N}$ modes. The size of the reduced basis $N$ selected by following criterion \eqref{eq:energy} is typically much smaller than the size of the FOM $N^h$, i.e. $N \ll N^h$.

\begin{remark}
\label{remark:ortho}
Given a symmetric positive definite matrix $\mat{X}^h$ which is a norm matrix for $\Vert \cdot \Vert_{\pazocal V}$ in the FE space, i.e. $\Vert u \Vert_\pazocal{V} = (\algvec{u}^h)^\text{T} \mat{X}^h \algvec{u}^h$, it is possible to perform the POD such that the basis $\mat{U}$ is orthonormal with respect to $\mat{X}^h$ (i.e. $\mat{U}^\text{T} \mat{X}^h \mat{U} = \mat{I}$). In order to achieve this, we observe that, since $\mat{X}^h$ is symmetric positive definite, it admits a Cholesky decomposition $\mat{X}^h = \mat{H}^\text{T} \mat{H}$, $\mat{H}$ being upper triangular. Matrix $\mat{U}$ is then found as $\mat{U} = \mat{H}^{-1} \widetilde{\mat{U}}$, where $\widetilde{\mat{U}}$ is computed by SVD of $\mat{H} \mat{S} = \widetilde{\mat{U}} \widetilde{\Sigma} \widetilde{\mat{V}}^\text{T}$. When constructing the reduced basis for our particular application in Section~\ref{sec:reducedrb}, following this approach allows us to achieve the optimality expressed in Proposition~\ref{prop:optimality} with respect to norms more suited to the specific variables of interest (namely, H1 norm for the velocity and L2 norm for the pressure).
\end{remark}

\subsection{The online phase: solution of the reduced problem}
\label{subsec:rbonline}
Let us observe that it is legitimate to associate with each POD mode a corresponding functional representation $\zeta^h_j = \sum_{i = 1}^{N^h} (\algvec{\zeta}^h_j)_i \varphi^h_i$. The RB approximation then reads $u^N = \sum_{i = 1}^{N} u^N_i \zeta^h_i$, $\algvec{u}^N = [u^N_1,\ldots,u^N_N]$ being the vector of reduced degrees of freedom. Then, evaluating the weak formulation of problem~\eqref{eq:differential_elliptic} at test and trial functions in $\text{span}\{\zeta^h_i\}_{i = 1}^{N}$, we find the reduced linear system
\begin{equation}
    \mat{L}^N(\param) \algvec u^N = \algvec{G}^N(\param),
    \label{eq:rbsystem}
\end{equation}
where $\mat{L}^N(\param)_{ij} = \ell(\zeta_j^h,\zeta_i^h) \in \reals^{N \times N}$ and $\algvec{G}^N(\param)_i = \int_{\Omega} \mathcal{G}(\param) \zeta_i^h \in \reals^{N}$. The assembly and solution of system \eqref{eq:rbsystem} correspond to the online phase. The transformation of the reduced vector of degrees of freedom into its FE counterpart is simply performed by $\algvec{u}^h \approx \mat{V} \algvec{u}^N \in \reals^{N^h}$.

By exploiting the expansion $\zeta^h_j = \sum_{i = 1}^{N^h} (\algvec{\zeta}^h_j)_i \varphi^h_i$ it is easy to find that $\mat{L}^N(\param) = \mat{V}^\text{T} \mat{L}^h(\param) \mat{V}$ and $\algvec{G}^N(\param) = \mat{V}^\text{T} \algvec{G}^h(\param)$. In the most general case, therefore, the assembly of the reduced system is done by constructing the full order matrix and right hand side and by computing their projection onto the RB space. If the problem features an affine decomposition, namely there exist parameter-dependent coefficents $\theta^L_{q}$ for $q = 1,\ldots,Q_L$ and $\theta^G_{q}$ for $q = 1,\ldots,Q_G$ such that
\begin{equation}
    \mat{L}^h(\param) = \sum_{q = 1}^{Q_L} \theta^L_{q}(\param) \mat{L}^h_q, \qquad \algvec{G}^h(\param) = \sum_{q = 1}^{Q_G} \theta^G_{q}(\param) \algvec{G}^h_q,
    \label{eq:affinedec}
\end{equation}
a considerable speedup is achieved by precomputing the matrices $\mat{L}^N_q = \mat{V}^\text{T} \mat{L}^h_q \mat{V}$ and the vectors $\algvec{G}^N_q = \mat{V}^\text{T} \algvec{G}^h_q$ in the offline phase, and by assembling the reduced elements of system \eqref{eq:rbsystem} as
\begin{equation}
    \mat{L}^N(\param) = \sum_{q = 1}^{Q_L} \theta^L_{q}(\param) \mat{L}^N_q, \qquad \algvec{G}^N(\param) = \sum_{q = 1}^{Q_G} \theta^G_{q}(\param) \algvec{G}^N_q.
\end{equation}
Unfortunately, in most practical scenarios an affine decomposition of the form \eqref{eq:affinedec} is not readily available. In such cases, a common strategy to efficiently perform the assembly of system \eqref{eq:rbsystem} consists in employing the (discrete) empirical interpolation method (DEIM) \cite{barrault2004empirical, Chaturantabut2010deim} and its matrix variant MDEIM \cite{negri2015efficient}. In this paper, we do not address the optimization of the assembly of the reduced system, which will be investigated in future works.
\begin{remark}
\label{remark:quality}
The quality of the RB approximation depends on three factors: the POD tolerance $\varepsilon$, the number of considered snapshots $N_s$, and the choice of sampling space. Assuming the latter to be appropriate, however, the number of snapshots required to achieve errors of the order of the POD tolerance  in the online phase may become too large as the dimension of the parameter space $\pazocal D$ increases. In other words, in applications were the space of parameters is too rich, it is unfeasible to sample a sufficient number of snapshots, and the online error of the RB method might be considerably greater than the one obtained on the snapshots.
\end{remark}

\section{The Navier--Stokes equations}
\label{sec:ns}
In this section, we first introduce the Navier--Stokes equations in strong and weak formulations (Section~\ref{subsec:nsweak}). The latter poses the mathematical foundation for the numerical discretization in space by the FE method as described in Section~\ref{subsec:numdiscns}, where we also derive the fully-discrete model by considering a generic Backward Differentiation Formulas (BDF) scheme.
\subsection{Strong and weak formulations}
\label{subsec:nsweak}
Let us consider the problem of approximating the blood flow in a vessel, mathematically represented by an open and bounded domain $\domain \in \realsd$. In this paper, we take $d = 3$, but the discussion is also valid for the case $d = 2$. We model
the blood as an incompressible Newtonian fluid and, therefore, its dynamics is described by the Navier--Stokes equations
\begin{equation}
    \begin{alignedat}{3}
    \rho_\text{f} \dfrac{\partial \uvec}{\partial t} + \rho_\text{f} (\uvec \cdot \nabla) \uvec - 2 \mu_\text{f} \nabla \cdot \varepsilon(\uvec) + \nabla \p &= \fvec  &&\qquad \text{in }\domain \times (0,T),\\
    \nabla \cdot \uvec &= 0 &&\qquad \text{in }\domain \times (0,T),\\
    \uvec &= \gvec &&\qquad \text{on }\dirboundary \times (0,T),\\
    \sigma(\uvec, \p) \nnorm &= \hvec &&\qquad \text{on }\neuboundary \times (0,T), \\
    \uvec &= \uvec_0 &&\qquad \text{for } t = 0,
    \end{alignedat}
    \label{eq:navierstokes}
\end{equation}
where $\uvec : \domain \times (0,T) \mapsto \realsd$ and
$\p : \domain \times (0,T) \mapsto \reals$ are velocity and pressure of the fluid,
$\rho_\text{f}$ is the density, $\mu_\text{f}$ is the viscosity, $\varepsilon(\uvec) = (\nabla \uvec + \nabla \uvec^\text{T})/2$ is
the strain rate tensor, $\sigma(\uvec,p) = 2 \mu_\text{f} \varepsilon(\uvec) - p \mat{I}$ is the Cauchy stress tensor,
$\fvec : \domain \times (0,T) \mapsto \realsd$ is a forcing term,
$\gvec : \dirboundary \times (0,T) \mapsto \realsd$ and
$\hvec : \neuboundary \times (0,T) \mapsto \realsd$ are Dirichlet and Neumann data,
$\nnorm$ is the normal unit vector to the boundary $\partial \domain$, and
$\uvec_0: \domain \mapsto \realsd$ is the prescribed initial condition. Since we deal with cardiovascular applications, we take $\dirboundary =  \Gamma_{\text{w}} \cup \left (\bigcup_{i=1}^{N_\text{in}}\Gamma^\iin{i}_{\text{in}}\right )$ and $\neuboundary = \bigcup_{i=1}^{N_\text{out}}\Gamma^\iin{i}_{\text{out}}$; $\Gamma^\iin{1}_{\text{in}},\ldots, \Gamma^\iin{{N_{\text{in}}}}_{\text{in}}$, $\Gamma_\text{w}$ and $\Gamma^\iin{1}_{N}, \ldots, \Gamma^\iin{{N_\text{out}}}_{\text{out}}$ are the inlets, wall and outlets of the vessel, respectively. The inlet velocity profiles and outlet Neumann data are denoted $\gvec_1,\ldots,\gvec_{N_\text{in}}$ and $\hvec_1,\ldots,\hvec_{N_\text{out}}$; on the wall $\Gamma_\text{w}$ we consider $\uvec = \boldsymbol{0}$.
The first equation in Eq.~\eqref{eq:navierstokes} (\textit{momentum equation})
represents the generalization of Newton's second law of motion to continuums, and
the second equation (\textit{continuity equation}) is the incompressibility constraint.

The numerical solution of the Navier--Stokes equations by classical Galerkin methods
such as the FE method entails transforming Eq.~\eqref{eq:navierstokes}
into its weak formulation. Let us denote $\pazocal V_{g} := [H^1_{g,\dirboundary}(\domain)]^d$---that is, the space of functions belonging to $[H^1(\domain)]^d$ such that their trace is equal to $\gvec$ on $\dirboundary$, $\pazocal Q := L^2(\domain)$,
$\pazocal V_{0} := [H^1_{0, \dirboundary}(\domain)]^d$, and let us consider generic test functions
$\vvec \in \pazocal V_{0}$ and $q \in \pazocal Q$. The weak formulation of Eq.~\eqref{eq:navierstokes} is obtained by multiplying the momentum and continuity equations by $\vvec$ and $q$ respectively and by integrating over the domain $\domain$. Hence, we find:
\begin{weakform}{}
given $\fvec$, $\gvec$, and $\hvec$ regular enough, find $(\uvec,\p) \in \pazocal V_{g} \times \pazocal Q$, such that, for every $t \in (0,T)$,
\begin{equation}
\begin{alignedat}{3}
        \int_{\domain} \rho_\text{f} \dfrac{\partial \uvec}{\partial t} \cdot \vvec + \int_{\domain} \rho_\text{f} [(\uvec \cdot \nabla) \uvec] \cdot \vvec + \int_{\domain} \sigma(\uvec,\p) : \nabla \vvec &= \int_{\domain} \fvec \cdot \vvec + \int_{\neuboundary} \hvec \cdot \vvec &&\qquad \forall \vvec \in \pazocal V_{0},\\
        \int_{\domain} \nabla \cdot \uvec q &= 0 &&\qquad \forall q \in \pazocal Q,
\end{alignedat}
\end{equation}
and such that $\uvec = \uvec_0$ for $t = 0$.
\label{weak:navierstokes}
\end{weakform}

\subsection{Numerical discretization}
\label{subsec:numdiscns}
In order to transform the infinite dimensional problem \weakref{weak:navierstokes} into
a finite dimensional one we consider the two subspaces $\pazocal V_{g}^h := \text{span}\{\boldsymbol \varphi_i^h \}_{i=1}^{\nuh} \subset \pazocal V_{g}$ and $\pazocal Q^h := \text{span}\{\psi_i^h \}_{i=1}^{\nph} \subset \pazocal Q$, and
the finite dimensional approximations of velocity $\uvech (\xvec,t) = \sum_{i = 1}^{\nuh} u_i^h (t)\boldsymbol \varphi_i^h (\xvec)$ and pressure $\ph (\xvec,t) = \sum_{i = 1}^{\nph} p_i^h(t) \psi_i^h (\xvec)$ (in these expressions we
explicitly highlight for the sake of clarity the dependance of each term on space and time, but this
is omitted hereon). The choice of basis functions for $\pazocal V_{g}^h$ and
$\pazocal Q^h$ clearly plays a crucial role in the accuracy of the approximation and it is a characteristic
of the discretization method of choice. Moreover, in the case of saddle-point problems such as the Navier--Stokes equations, the quality of the discretization is critical to ensure the well-posedness of the
discrete problem; we refer the reader to Section~\ref{subsec:supremizers} and \cite{boffi2013mixed,brezzi1974existence} for details. In this section and in Section~\ref{sec:modular}, the basis functions $\{\boldsymbol \varphi_i^h \}_{i=1}^{\nuh}$ and  $\{\psi_i^h \}_{i=1}^{\nph}$ are standard Lagrangian P2-P1 Taylor--Hood FE basis functions \cite{hood1974navier} (i.e. quadratic and linear piece-wise polynomials for the velocity and the pressure, respectively) obtained from a triangulation of the domain $\mathcal T^h$ composed of tetrahedra. We consider other possibilities---i.e. RB functions---in Section~\ref{sec:nsmodularrb}.

By introducing the vectors of degrees of freedom $\algvec{u}^h = [u_1^h,\ldots,u^h_{\nuh}] \in \reals^{\nuh}$, $\algvec{p}^h = [p^h_1,\ldots,p^h_{\nph}] \in \reals^{\nph}$ and $\algvec f^h_i = \int_{\domain} \fvec \cdot \boldsymbol \varphi_i^h + \int_{\neuboundary} \hvec \cdot \boldsymbol \varphi_i^h \in \reals^{\nuh}$, and the matrices
$\mat{M}^h_{ij} = \int_{\domain} \mu_\text{f} \boldsymbol \varphi_j^h \cdot \boldsymbol \varphi_i^h \in \reals^{\nuh \times \nuh}$ (mass), $\mat{K}^h_{ij} = \int_{\domain} 2 \mu_\text{f} \varepsilon(\boldsymbol \varphi_j^h) : \varepsilon(\boldsymbol \varphi_i^h) \in \reals^{\nuh \times \nuh}$ (stiffness), $\mat{C}^h(\uvech)_{ij} = \int_{\domain} \rho_\text{f} [(\uvech \cdot \nabla)\boldsymbol \varphi_j^h] \cdot \boldsymbol \varphi_i^h \in \reals^{\nuh \times \nuh}$ (convective matrix), and $\mat{D}^h_{ij} = -\int_{\domain} \nabla \cdot \boldsymbol \varphi_j^h \psi_i^h \in \reals^{\nph \times \nuh}$ (divergence), the discrete version of \weakref{weak:navierstokes} is conveniently expressed in the form of linear system as
\begin{equation}
    \begin{bmatrix}
        \mat{M}^h & \\
          &
    \end{bmatrix}
    \begin{bmatrix}
        \dot{\algvec{u}}^h \\
        \dot{\algvec{p}}^h
    \end{bmatrix}
    +
    \begin{bmatrix}
        \mat{K}^h + \mat{C}^h\left (\uvech \right ) & (\mat{D}^h)^\text{T} \\
        \mat{D}^h & \\
    \end{bmatrix}
    \begin{bmatrix}
        \algvec{u}^h \\
        \algvec{p}^h
    \end{bmatrix}
    =
    \begin{bmatrix}
        \algvec{f}^h \\
        \algvec{0}
    \end{bmatrix}.
    \label{eq:semidiscretens}
\end{equation}
We exploit that $\uvech = \sum_{i = 1}^{\nuh} u^h_i \boldsymbol \varphi_i^h$
and indicate the convective term matrix as $\mat{C}^h(\algvec{u}^h)$---i.e. as a function
of the degrees of freedom $\algvec{u}^h$ instead of the approximated function $\uvech$---in the remainder of the paper.

Let us now introduce a sequence of timesteps $t_0, t_1, \ldots, t_{N_t}$ such that
$t_0 = 0$, $t_{N_t} = T$, and $t_{k+1} = t_{k} + \Delta t$ for every $k = 0,\ldots,N_t$;
$\Delta t$ is called timestep size.
We denote the value of $\algvec{u}^h$ and $\algvec{p}^h$ at timestep $t_k$ by $\algvec{u}^h(t_k) = \algvec{u}^h_k$ and $\algvec{p}^h(t_k) = \algvec{p}^h_k$, respectively. The numerical discretization in time of Eq.~\eqref{eq:semidiscretens} is performed
by means of BDF schemes. Specifically, given
$\algvec{u}^h_{k-j+1}$ and $\algvec{p}^h_{k-j+1}$ for $j = 1,\ldots,\sigma$, the numerical solution of the Navier--Stokes equations at timestep $t_{k+1}$ by a BDF scheme of order $\sigma$ satisfies
\begin{equation}
    \algvec{r}\left (\algvec{w}^h_{k+1} \right ) := \mat{H}^h \algvec{w}^h_{k+1} - \sum_{j = 1}^{\sigma} \alpha_j \mat{H}^h \algvec{w}^h_{k-j+1} - \Delta t \beta \mathring{\algvec{f}}^h\left (t_{k+1},\algvec{w}^h_{k+1} \right ) = \algvec{0},
    \label{eq:residualbdf}
\end{equation}
where $\algvec{w}^h = [\algvec{u}^h, \algvec{p}^h] \in \reals^{\nuh + \nph}$ and
\begin{equation}
    \mat{H}^h := \begin{bmatrix}
        \mat{M}^h & \\
          &
    \end{bmatrix},\quad
    \mathring{\algvec{f}}^h\left (t,\algvec w^h\right) :=
    \begin{bmatrix}
        \algvec{f}^h(t) \\
        \algvec{0}
    \end{bmatrix}
    -
    \begin{bmatrix}
        \mat{K}^h + \mat{C}^h\left (\algvec u^h\right ) & (\mat{D}^h)^\text{T} \\
        \mat{D}^h & \\
    \end{bmatrix}
    \begin{bmatrix}
        \algvec{u}^h \\
        \algvec{p}^h
    \end{bmatrix}.
\end{equation}
The coefficients $\{\alpha_j\}_{j = 1}^{\sigma}$ and $\beta$ depend on the specific BDF scheme. For example, for $\sigma = 1$ we have $\alpha_1 = \beta = 1$ (Backward Euler scheme), and, for $\sigma = 2$, $\alpha_1 = 4/3$, $\alpha_2 = -1/3$ and $\beta = 2/3$; these two choices lead to numerical methods of first and second order, respectively.

Eq.~\eqref{eq:residualbdf} is in general nonlinear and its solution requires the application of ad-hoc numerical methods; in this paper, we adopt the Newton--Raphson algorithm (see Section~\ref{subsec:preconditioner}).

\section{Modular domain-decomposition of arteries}
\label{sec:modular}
With the aim of the model order reduction presented in Section~\ref{sec:nsmodularrb}, it is beneficial to perform a geometrical approximation of the vessel based on a domain-decomposition approach. We first set the theoretical basis for the approximation of the Navier--Stokes equations on modular geometries in Section~\ref{subsec:modular}. The discretization of the continuous formulation is then performed in Section~\ref{subsec:discretization_dd}, where we also present our approach for the treatment of the coupling variables (i.e. Lagrange multipliers) at the interfaces.
\subsection{The continuous Navier--Stokes equations on modular geometries}
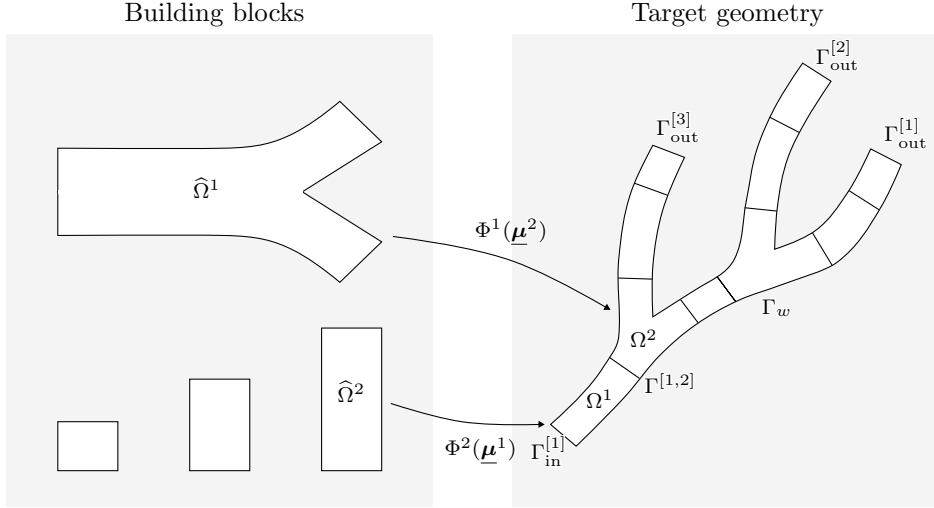
\begin{figure}
    \centering
    \tikzsetnextfilename{sketch_dd}
    \begin{tikzpicture}
        \node[inner sep=0pt] (bifs) at (0,0)
        {\includegraphics[scale = 0.4,trim = {2.6cm 1cm 2.8cm 1cm},clip]{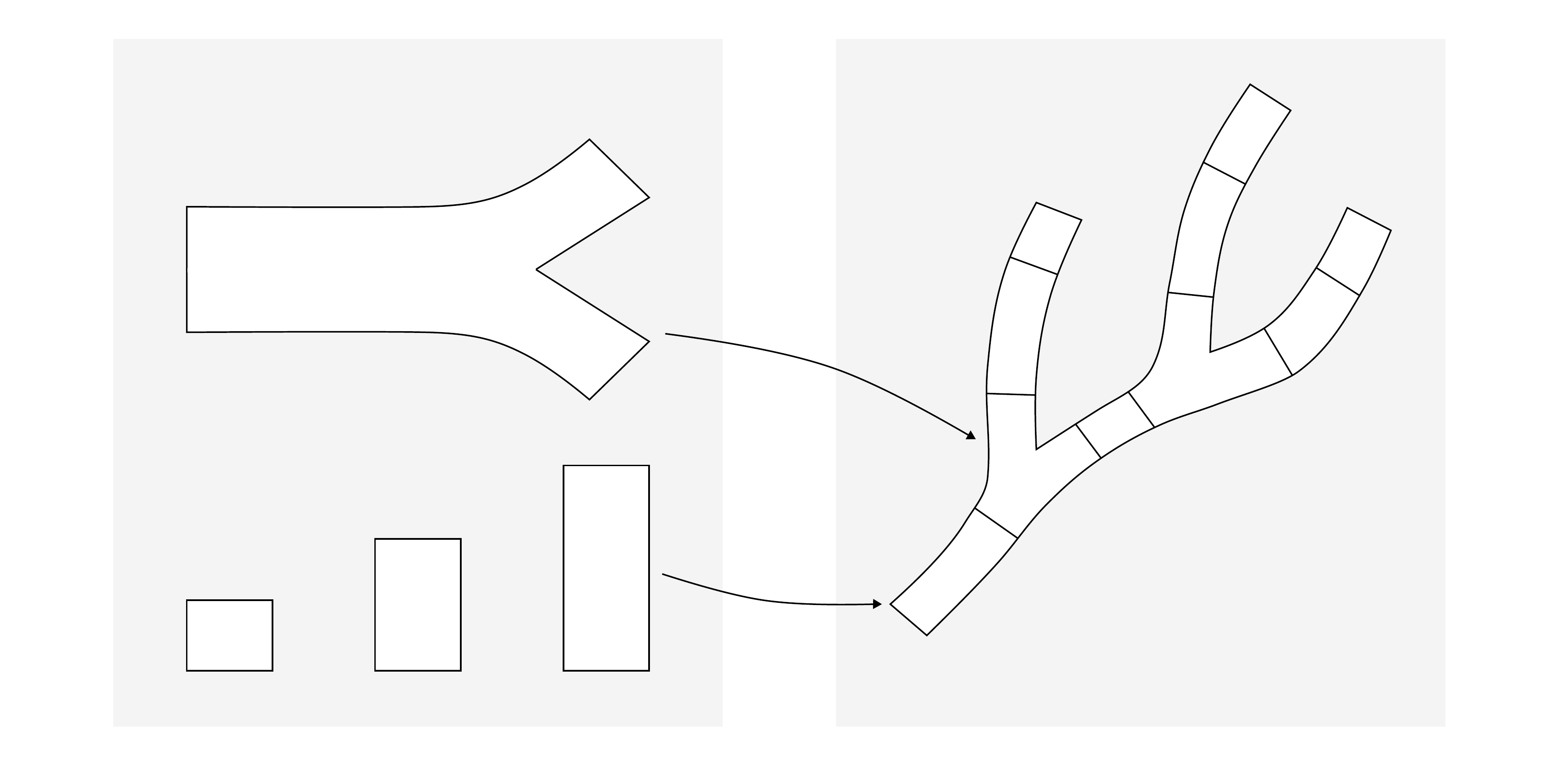}};
        \node  at (0.5,0.5) {\contour{white}{\footnotesize $\Phi^1(\param^2)$}};
        \node  at (0.1,-2.4) {\contour{white}{\footnotesize $\Phi^2(\param^1)$}};
        \node  at (-1.6,-1.6) {\footnotesize $\bbref^{\is{2}}$};
        \node  at (-3.5,1.1) {\footnotesize $\bbref^{\is{1}}$};
        \node  at (2.25,-0.9) {\footnotesize $\domain^{\is{2}}$};
        \node  at (1.67,-1.7) {\footnotesize $\domain^{\is{1}}$};
        \node  at (2.6,-1.5) {\contour{white}{\footnotesize $\Gamma^{\iin{1,2}}$}};
        \node  at (1,-2.4) {\contour{white}{\footnotesize $\Gamma_{\text{in}}^{\iin{1}}$}};
        \node  at (2.7,1.9) {\contour{white}{\footnotesize $\Gamma_{\text{out}}^{\iin{3}}$}};
        \node  at (5.7,1.82) {\contour{white}{\footnotesize $\Gamma_{\text{out}}^{\iin{1}}$}};
        \node  at (4.8,2.85) {\contour{white}{\footnotesize $\Gamma_{\text{out}}^{\iin{2}}$}};
        \node  at (4,-0.5) {\contour{white}{\footnotesize $\Gamma_{w}$}};
        \node  at (-3.4,3.4) {Building blocks};
        \node  at (3.35,3.4) {Target geometry};
    \end{tikzpicture}
    \caption{Sketch of the domain-decomposition of a target geometry. Each block in the target geometry is found from the parametrized geometrical deformation of a small number of reference building blocks.}
    \label{fig:ddsketch}
\end{figure}
\begin{figure}
    \centering
    \tikzsetnextfilename{deformations}
    \begin{tikzpicture}
        \node[inner sep=0pt] (bifs) at (0,0)
        {\includegraphics[width = \textwidth,trim = {18cm 1cm 11cm 13cm},clip]{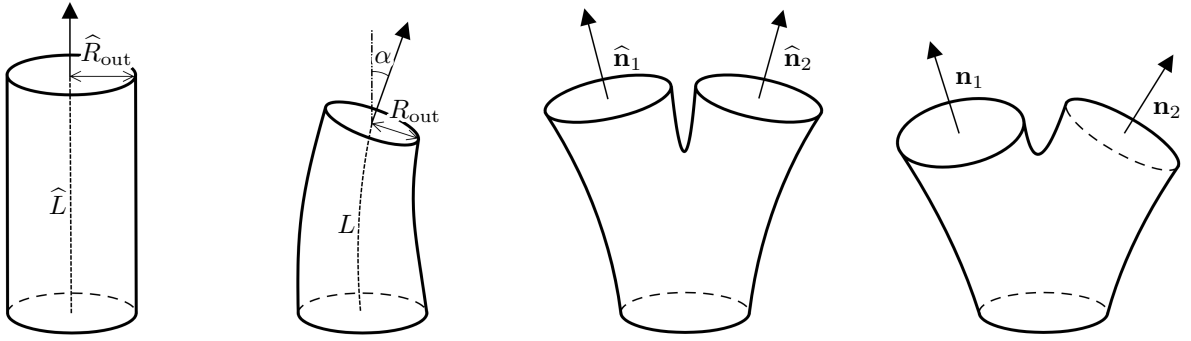}};
        \node  at (-7.1,-0.2) {\contourlength{1.3pt} \contour{white}{$\widehat{L}$}};
        \node  at (-6.4,1.8) {\contourlength{1.5pt}\contour{white}{$\widehat{R}_\text{out}$}};
        \node  at (-3.3,-0.5) {\contourlength{1.3pt} \contour{white}{$L$}};
        \node  at (-2.4,1) {\contourlength{1.3pt} \contour{white}{$R_\text{out}$}};
        \node  at (-2.82,1.7) {\contourlength{1.3pt} \contour{white}{$\alpha$}};
        \node  at (0.4,1.7) {\contourlength{1.3pt} \contour{white}{$\widehat{\mathbf{n}}_1$}};
        \node  at (2.7,1.7) {\contourlength{1.3pt}\contour{white}{$\widehat{\mathbf{n}}_2$}};
        \node  at (4.9,1.4) {\contourlength{1.3pt} \contour{white}{$\mathbf{n}_1$}};
        \node  at (7.5,1) {\contourlength{1.3pt} \contour{white}{$\mathbf{n}_2$}};
    \end{tikzpicture}
    \caption{Types of reference building blocks and affine transformations. On the left, tubes: the geometrical parameters are the angle of the outlet normal $\alpha$---due to the axial symmetry and to the rotation matrix $\mat Q$ in Eq.~\eqref{eq:transformation}, a single angle is sufficient to represent a bending in any direction---in the deformed configuration $\alpha$ and ratios between the reference and deformed lengths ($\widehat{L}/L$) and reference radiae ($\widehat{R}/R$). On the right, bifurcation: the geometrical parameters are the angles describing the rotation of the reference outlet normals $\widehat{\mathbf{n}}_1$ and $\widehat{\mathbf{n}}_2$ onto the outlet normals $\mathbf{n}_1$ and $\mathbf{n}_2$ (i.e. three Euler angles per outlet, that is six geometrical parameters in total).}
    \label{fig:deformations}
\end{figure}

\label{subsec:modular}
We introduce a library of \textit{building blocks} $\bbref^{\irs{i}}$, $i = 1,\ldots,N_{\text{bb}}$. In the context of cardiovascular simulations, these reference building blocks are model cylinders and bifurcations, as shown in Fig.~\ref{fig:ddsketch}. The target geometry is then approximated as a modular composition of subdomains $\domain \approx \domain_\text{m}(\paramspace) = \bigcup_{j = 1}^{N_\domain} \domain^{\is{j}}(\param^\is{j})$. Here, $\domain^{\is{j}} := \Phi^\irs{z(j)}(\bbref^{\irs{z(j)}};\param^\is{j})$ is an open and bounded subdomain obtained by applying a prescribed parametrized geometrical deformation $\Phi^\irs{z(j)}$ to the $z(j)^{\text{th}}$ building block, $z: [1,\ldots,N_\domain] \mapsto [1,\ldots,N_\text{bb}]$ is an injective map from the indices of the subdomains in the target geometry to the indices of the building blocks, and $\mathcal{M} := \{ \param^\is{j} \}_{j = 1}^{N_\domain}$ is the set of geometrical parameters. In the following, we indicate $z(j) := z_j$ for brevity. Each vector of parameters $\param^j$ belongs to a space $\pazocal D^\irs{z_j} \subset \reals^{N_\mu^\irs{z_j}}$ whose dimensionality depends on the corresponding reference building block. For each $i = 1,\ldots,N_\text{bb}$ and given a parameter vector $\param$, we focus on geometrical deformations of the form
\begin{equation}
\Phi^i(\widehat{\mathbf{x}};\param) = \mat Q(\param) \boldsymbol {\varphi}^i(\widehat{\mathbf x};\param) + \mathbf{t}(\param), \quad \forall \widehat{\mathbf{x}} \in \widehat{\Omega}^i,
\label{eq:transformation}
\end{equation}
where $Q(\param)$ is a rotation matrix, $\mathbf{t}$ is a translation vector, and $\boldsymbol \varphi^i(\cdot;\param)$ is a nonaffine geometrical deformation. The types of building blocks we consider and the corresponding adimissible nonaffine deformations are depicted in Fig.~\ref{fig:deformations}.

The subdomains in the target geometry satisfy $\domain^\is{l}(\param^\is{l}) \cap \domain^\is{m}(\param^\is{m}) = \emptyset$ if $l \neq m$, and we define the interface $\Gamma^{\iin{jm}}(\param^\is{j},\param^\is{m}) := \overline{\domain}^\is{j}(\param^\is{j}) \cap \overline{\domain}^\is{m}(\param^\is{m})$. The building blocks in the reference configuration are designed with circular inlet and outlet faces; the geometrical deformations are chosen such that the interfaces are circles for every possible choice of the geometrical parameters.

\begin{remark}
    Although $\Gamma^\iin{lm}$ and $\Gamma^\iin{ml}$ represent the same physical surface, it is still beneficial to differentiate between the two as we associate with each interface the vectors $\nnorm_{lm}$ and $\nnorm_{ml}$, i.e. the outward normal unit vectors with respect to $\domain^\is{l}$ and $\domain^\is{m}$ (clearly, $\nnorm_{lm} = -\nnorm_{ml}$). This distinction allows to simplify the notation in \weakref{weak:navierstokes_dd}.
\end{remark}

For every subdomain $\domain^\is{j}(\param^\is{j})$, we introduce the set of indices of the neighboring subdomains $N(j)$ and the sets $I_\text{in}(j)$ and $I_\text{out}(j)$ such that $\partial \domain^\is{j}(\param^\is{j}) \cap \Gamma_\text{in}^\iin{i} \neq \emptyset$ for all $i \in I_\text{in}(j)$, and $\partial \domain^\is{j}(\param^\is{j}) \cap \Gamma_\text{out}^\iin{i} \neq \emptyset$ for all $i \in I_\text{out}(j)$. The dependance of the subdomains $\domain^\is{j}$ and interfaces $\Gamma^\iin{jm}$, $\Gamma^\iin{i}_\text{in}$ and $\Gamma^\iin{i}_\text{out}$ on the geometrical parameters $\paramspace$ is omitted unless ambiguity arises in the remainder of the paper.

Let us rewrite the weak formulation \weakref{weak:navierstokes} to account for the modular decomposition of the original domain, i.e. $\domain_\text{m}$. For each subdomain $\domain^\is{j}$, we introduce the spaces $\pazocal V_g^{\is{j}} := [H^1_{g, \Gamma_D}(\domain^\is{j})]^d$, $\pazocal Q^{\is{j}} := L^2(\domain^\is{j})$ and $\pazocal V^{\is{j}}_{0} := [H^1_{0, \Gamma_D}(\domain^\is{j})]^d$. Moreover, for every interface $\Gamma^\iin{jm}$ we define the spaces $\pazocal L^{\iin{jm}} = [H^{-1/2}_{00}(\Gamma^\iin{jm})]^d$. For the sake of conciseness, let us use the following notation
\begin{equation}
    \begin{aligned}
        \mathcal M^{\is{j}}(\boldsymbol \varepsilon, \boldsymbol \varphi, \psi; \boldsymbol \omega) &:= \int_{\domain^\is{j}} \rho_\text{f} \boldsymbol \varepsilon \cdot \boldsymbol \omega + \int_{\domain^\is{j}} \rho_\text{f} [(\boldsymbol \varphi \cdot \nabla) \boldsymbol \varphi] \cdot \boldsymbol \omega + \int_{\domain^\is{j}} \sigma(\boldsymbol \varphi,\psi) : \nabla \boldsymbol \omega \\
        \mathcal C^{\is{j}}(\boldsymbol \varphi; \eta) &:= \int_{\domain^\is{j}} \nabla \cdot \boldsymbol \varphi \eta,
    \end{aligned}
\end{equation}
for every $\boldsymbol \varepsilon, \boldsymbol \varphi, \boldsymbol \omega \in [H^1(\domain^\is{j})]^d$ and for every $\psi,\eta \in L^2(\domain^\is{j})$. Assuming for simplicity that $\domain = \domain_\text{m}$, it can be shown that \weakref{weak:navierstokes} is equivalent---in a sense that will specified in Remark~\ref{remark:equivalence}---to the following weak formulation:
\begin{weakform}{}
given $\fvec$, $\gvec$, and $\hvec$ regular enough and for every $j = 1,\ldots,N_\domain$, find $(\uvec^{\is{j}},\p^{\is{j}},\{\boldsymbol{\lambda}^{\iin{jm}}\}_{m \in N(j)}) \in \pazocal V^{\is{j}}_{g} \times \pazocal Q^{\is{j}} \times \prod_{m \in N(j)} \pazocal L^{\iin{jm}}$, such that, for every $t \in (0,T)$,
\begin{equation}
\begin{alignedat}{3}
        \mathcal M^\is{j}\left (\dfrac{\partial \uvec^{\is{j}}}{\partial t},\uvec^{\is{j}},p^{\is{j}}; \vvec \right ) + \sum_{m \in N(j)} \int_{\Gamma^\iin{jm}} \boldsymbol \lambda^{\iin{jm}} \cdot \vvec &=  \int_{\domain^\is{j}} \fvec \cdot \vvec + \sum_{i \in I_\text{out}(j)}\int_{\Gamma_\text{out}^\iin{i}} \hvec \cdot \vvec&&\qquad \forall \vvec \in \pazocal V^{\is{j}}_{0},\\
        \mathcal C^{\is{j}}\left (\uvec^{\is{j}};q \right) &= 0 &&\qquad \forall q \in \pazocal Q^{\is{j}},
\end{alignedat}
\label{eq:ns_dd}
\end{equation}
and such that $\uvec^{\is{j}} = \uvec_0|_{\domain^\is{j}}$ for $t = 0$. Moreover, for every $m \in N(j)$, $\boldsymbol \lambda^\iin{jm} = -\boldsymbol \lambda^\iin{mj}$ and
\begin{equation}
\int_{\Gamma^\iin{jm}} \boldsymbol \eta \cdot \left (\uvec^{\is{j}} - \uvec^{\is{m}} \right ) = 0 \qquad \forall \boldsymbol \eta \in \pazocal L^{\iin{jm}}.
\label{eq:continuity_dd}
\end{equation}
\label{weak:navierstokes_dd}
\end{weakform}
\begin{remark}
    \label{remark:equivalence}
    The two weak formulations \weakref{weak:navierstokes} and \weakref{weak:navierstokes_dd} are equivalent in the following sense. If $(\uvec,p)$ satisfies \weakref{weak:navierstokes}, then $(\uvec|_{\domain^\is{j}},p|_{\domain^\is{j}},\{\sigma(\uvec,p)\nnorm_{jm}\}_{m \in N(j)})$ is also solution of the local problem \weakref{weak:navierstokes_dd}, for every $j = 1,\ldots,N_\domain$. The Lagrange multipliers therefore play the role of the stress at the interfaces; for details, see e.g. \cite{wohlmuth2000mortar,deparis2019coupling}. Conversely, if $(\uvec^{\is{j}},\p^{\is{j}},\{\boldsymbol{\lambda}^{\iin{jm}}\}_{m \in N(j)})$ are the local solutions of \weakref{weak:navierstokes_dd}, then $(\uvec,p) = (\Pi_{j = 1}^{N_\domain} \uvec^{\is{j}}, \Pi_{j = 1}^{N_\domain} p^{\is{j}})$ is solution of \weakref{weak:navierstokes}.
\end{remark}

\subsection{Discretization of the primal hybrid formulation of the flow problem}
\label{subsec:discretization_dd}

The discretization of differential problems by the FE method requires the definition of a computational mesh, as already mentioned in Section~\ref{subsec:numdiscns}. In the case of the approximated modular geometry $\domain_{\text{m}}$, each building block $\bbref^\is{i}$ is equipped with its own triangulation $\widehat{\mathcal T}^{\is{i},h}$. Therefore, the global mesh $\mathcal{T}^{h}_\text{m}(\mathcal M) = \bigcup_{j = 1}^{\text{N}_\domain} \mathcal{T}^{\is{j},h}(\param^\is{j}) = \bigcup_{j = 1}^{\text{N}_\domain} \Phi^{z_j}(\widehat{\mathcal T}^{\is{z_j},h};\param^\is{j})$ is a composition of
distinct meshes which do not necessarily satisfy conformity constraints at the interfaces. We recall that a mesh is conforming if, for every two elements, their intersection is either empty or a whole face. In other words, a conforming mesh does not feature any \textit{hanging nodes}.

The global mesh is in general nonconforming and we are compelled to consider nonconforming domain-decomposition methods for the solution of the Navier--Stokes equations. These are formally defined as domain-decomposition methods in which the search space for the discrete solution is not a subset of the continuous search space (in our case $\pazocal V_g \times \pazocal Q$). The most popular approaches rely on the introduction of suitable Lagrange multipliers at the interfaces to enforce transmission conditions, as in \weakref{weak:navierstokes_dd}; see e.g. the well-known mortar method  \cite{bernardi1989new,bernardi2005basics,braess1999multigrid} and INTERNODES \cite{deparis2016internodes,gervasio2016analysis}. In this paper, we adopt the algorithm presented in \cite{deparis2019coupling}, which is based on the discretization of the Lagrange multipliers space via a small number of spectral basis functions defined on the interfaces. For our application, this choice is convenient because (i) the method allows to recover the $h$-convergence order of the primal discretization method---i.e. the FE method---even when a small number of spectral basis functions is considered, (ii) defining a spectral basis on each interface does not require to project nor to interpolate the traces of FE basis functions from one side to the other, and (iii) as already mentioned, the interfaces are circular in the target configuration, which allows us to employ a set of standard orthonormal basis functions on the two-dimensional disk. In the remainder of this section and in Section~\ref{subsec:discretization_dd} we recall the basics of this nonconforming method applied to the Navier--Stokes equations; the interested reader is referred to \cite{deparis2019coupling} for the details.

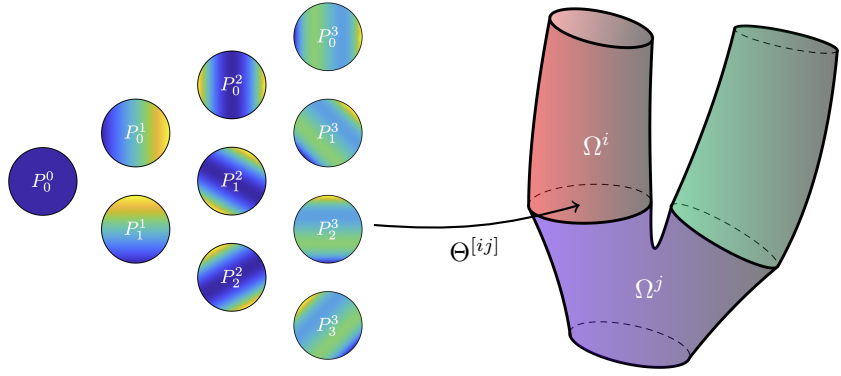
\begin{figure}
    \tikzsetnextfilename{chebbfs_mapped}
    \centering
    \begin{tikzpicture}
        \node[inner sep=0pt] (BFS) {\includegraphics[height = 0.34\textwidth,trim= {0cm 0.2cm 0 0.2cm},clip]{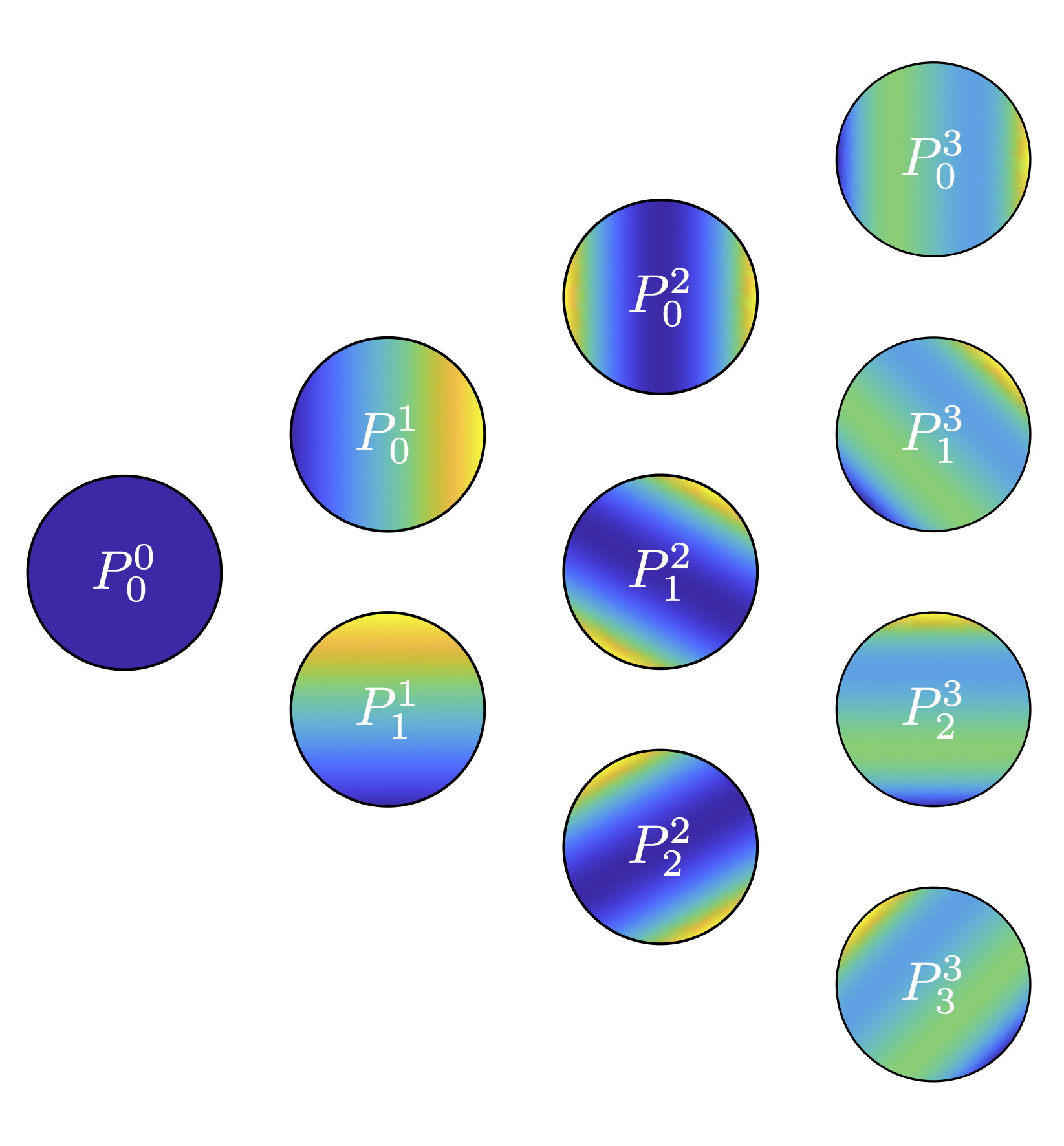}};
        \node[right=2cm of BFS, inner sep=0.5pt] (interface) {\includegraphics[height = 0.32\textwidth,trim = {6.5cm 3cm 5cm 3cm},clip]{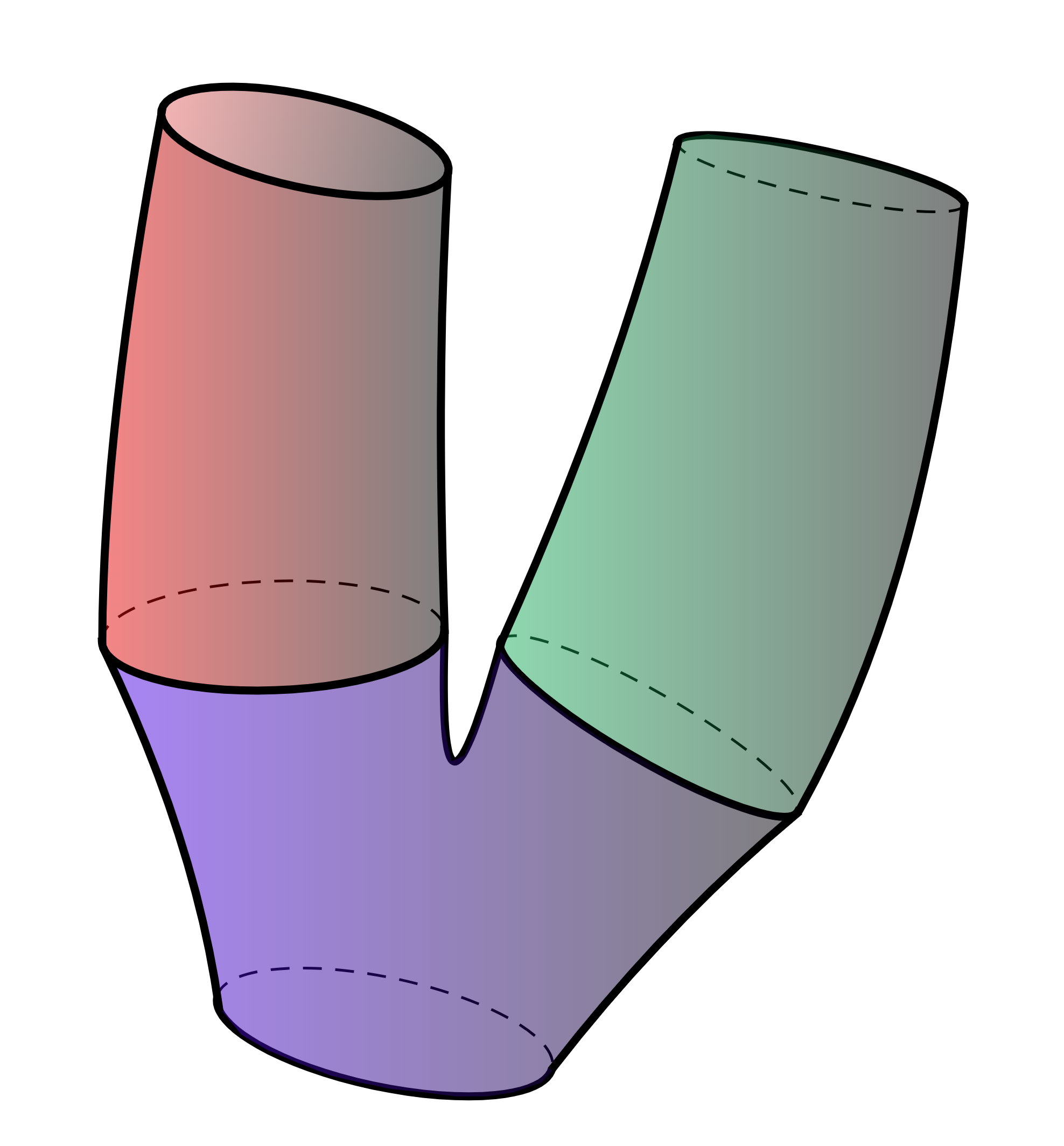}};
        \node at (5.4cm,0.5cm) {\textcolor{white}{$\domain^\is{i}$}};
        \node at (6.1cm,-1.4cm) {\textcolor{white}{$\domain^\is{j}$}};
        \node at (5.3cm,-0.3cm) (center) {};
        \draw[->,bend right = 10,thick] (BFS) edge (center);
        \node at (3.8cm,-0.9cm) {$\Theta^{\iin{ij}}$};
    \end{tikzpicture}
    \caption{Basis functions $P_k^n$ on the unit disk $\mathcal D$, for $n \leq 3$, mapped onto the target interface $\Gamma^{\iin{ij}}$. Minimum and maximum values are plotted in blue and yellow, respectively. The basis function $P_0^0$ is constant.}
    \label{fig:bbfs}
\end{figure}

We follow the same procedure presented in Section~\ref{subsec:numdiscns} for the discretization of the local variables $\uvec^{\is{j}}$ and $p^{\is{j}}$, which become $\uvec^{\is{j},h} = \sum_{i = 1}^{N_u^{\is{z_j},h}} u_i^{\is{j},h} \boldsymbol \varphi^{\is{j},h}_i$ and $p^{\is{j},h} = \sum_{i=1}^{{N_p^{\is{z_j},h}}} p^{\is{j},h}_i \psi^{\is{j},h}_i$. We remark that---since the meshes $\mathcal T^{j,h}$ are generated as transformations of the reference ones $\widehat{\mathcal T}^{z_j,h}$---the number of nodes $N_u^{\is{z_j},h}$ and $N_p^{\is{z_j},h}$ are indexed with respect to the reference building blocks; moreover, we have that $\boldsymbol \varphi^{\is{j},h}_i := \widehat{\boldsymbol \varphi}^{\irs{z_j},h} \circ (\Phi^\irs{z_j}(\param^j))^{-1}$, i.e. the velocity FE basis functions on the subdomains are obtained as composition of the reference velocity FE basis functions $\widehat{\boldsymbol \varphi}^{\irs{z_j},h}$ with the inverse of the map $\Phi^\irs{z_j}(\param^j)$ (the same holds true for the pressure FE basis functions $\psi^{\is{j},h}_i$). For the Lagrange multiplier $\boldsymbol \lambda^{\iin{jm}}$ we consider the approximation $\boldsymbol \lambda^{\iin{jm},\delta} = \sum_{i = 1}^{N_{\lambda}} \lambda_i^{\iin{jm},\delta} \bm  \xi_i^{\iin{jm},\delta}$, where $\boldsymbol \xi_i^{\iin{jm},\delta} = c_{jm}\widehat{\bm \xi}_i^\delta \circ (\Theta^{\iin{jm}})^{-1}\in [L^2(\Gamma^\iin{jm})]^d$, $c_{jm} \in \reals$, $\widehat{\bm \xi}_i^\delta \in [L^2(\mathcal D)]^d$ are a set of orthogonal functions on the unit disk $\mathcal D \in \reals^{d-1}$ and $\Theta^{\iin{jm}}: \mathcal D \mapsto \Gamma^{\iin{jm}}$ is a bijective map from the unit disk $\mathcal D$ to the target interface $\Gamma^{\iin{jm}}$. Similarly, we define the maps $\Theta^{\iin{m}}_\text{in}: \mathcal D \mapsto \Gamma^{\iin{m}}_\text{in}$ from the unit disk $\mathcal D$ to the inlet interfaces $\Gamma^{\iin{m}}_\text{in}$, the basis functions $\boldsymbol \xi_{i,\text{in}}^{\iin{m},\delta} = \widehat{\bm \xi}_i^\delta \circ (\Theta^{\iin{m}}_\text{in})^{-1}\in [L^2(\Gamma^\iin{m}_\text{in})]^d$, and $\boldsymbol \lambda^{\iin{m},\delta}_\text{in} = \sum_{i = 1}^{N_{\lambda}} \lambda_{i,\text{in}}^{\iin{m},\delta} \bm  \xi_{i,\text{in}}^{\iin{m},\delta}$. These functions are used to set the inlet velocity profiles, allowing one to effortlessly transition from the FE model to the reduced one through the process described in Section~\ref{sec:nsmodularrb}. The use of Lagrange multipliers is a classical way to weakly impose Dirichlet boundary conditions (see e.g. \cite{babuvska1973finite}) and, compared to other popular approaches such as penalty methods \cite{zhu1998modified,babuvska1973penalty}, it has the advantage of being variationally consistent. We remark that (i) we introduce the discretization parameter $\delta$ for the Lagrange multipliers to indicate that the degree of refinement is in fact independent of the mesh size in $\mathcal{T}^{\is{j},h}$ or $\mathcal{T}^{\is{m},h}$, and (ii) we consider for simplicity the same number $N_\lambda$ of basis functions at each interface.

In this paper, we construct $\{\widehat{\bm \xi}^\delta_i\}_{i = 1}^{N_\lambda}$ as follows. Let us consider Chebyshev polynomials of the second kind $U_n$, which are defined through the recurrence relation $U_0(x) = 1$, $U_1(x) = 2x$, $U_{n+1}(x) = 2x U_{n}(x) - U_{n-1}(x)$. Then, for $0 \leq k \leq n$,
\begin{equation}
P^n_k(x,y) = \dfrac{1}{\sqrt \pi }U_n \left (x \cos (\omega x) + y \sin (\omega y) \right ),\quad \omega = \dfrac{k}{n+1} \pi,
\label{eq:chebyshev}
\end{equation}
are orthonormal polynomials on the unit disk $\mathcal D$ with respect to the weight function $W(x,y) = 1 / \sqrt{\pi}$ \cite{dunkl2014orthogonal}. Given $n \geq 0$, we set
\begin{equation}
\widehat{\Xi}_n := \{\widehat{\bm \xi^\delta}_i\}_{i = 1}^{N_\lambda} = \bigcup_{i = 1}^d \bigcup_{\widetilde{n} = 0}^n \bigcup_{k = 0}^{\widetilde{n}} P^{\widetilde{n}}_k \mathbf e_i,
\end{equation}
where $\mathbf e_i$ is the $i^{\text{th}}$ canonic vector. It is trivial to find that $N_\lambda = d(n+1)(n+2)/2$.

Let us now address the discretization of the individual elements of \weakref{weak:navierstokes_dd}, which is local to subdomain $\domain^\is{j}$, the assembly of the global block system, and the discretization in time.

\paragraph{Discretization of the coupled momentum and continuity equations} For every $j = 1,\ldots,N_\domain$ and $m \in N(j)$, let us define $\algvec{\boldsymbol{\lambda}}^{\iin{jm},\delta} = [\lambda^{\iin{jm},\delta}_1,\ldots,\lambda^{\iin{jm},\delta}_{N_\lambda}]$ and the coupling matrix $\mat{B}_{pq}^{\iin{jm},h\delta} = \int_{\Gamma^\iin{jm}} \bm \xi_p^{\iin{jm},\delta} \cdot \boldsymbol \varphi^{\is{j},h}_q$. Then, Eq.~\eqref{eq:ns_dd} can be rewritten in algebraic form as
\begin{equation}
    \begin{bmatrix}
        \mat{M}^{\is{j},h} & & \\
        & & \\
        & & \\
    \end{bmatrix}
    \begin{bmatrix}
        \dot{\algvec{u}}^\is{j,h}\\
        \dot{\algvec{p}}^\is{j,h}\\
        \dot{\algvec{\lambda}}^\is{j,\delta}
    \end{bmatrix}
    +
    \begin{bmatrix}
        \mat{K}^\is{j,h} + \mat{C}^\is{j,h}\left (\algvec u^\is{j,h}\right ) & (\mat{D}^\is{j,h})^\text{T} & (\mat{B}^\is{j,h\delta})^\text{T}\\
        \mat{D}^\is{j,h} & \\
    \end{bmatrix}
    \begin{bmatrix}
        \algvec{u}^\is{j,h}\\
        \algvec{p}^\is{j,h}\\
        \algvec{\lambda}^\is{j,\delta}
    \end{bmatrix}
    =
    \begin{bmatrix}
        \algvec{f}^\is{j,h} \\
        \algvec{0}
    \end{bmatrix}.
    \label{eq:ns_dd_discr}
\end{equation}
In Eq.~\eqref{eq:ns_dd_discr}, we denoted by $\algvec{\lambda}^{\is{j},\delta}$ the vector containing all $\algvec{\lambda}^{\iin{jm},\delta}$ and by $\mat{B}^{\is{j},h\delta}$ the matrix obtained by stacking the various $\mat{B}^{\iin{jm},h\delta}$, for every $m \in N(j)$.

\paragraph{Equality of the Lagrange multipliers at the same interface} A natural way to enforce the constraint $\boldsymbol \lambda^{\iin{jm}} = - \boldsymbol \lambda^{\iin{mj}}$ in \weakref{weak:navierstokes_dd} is to choose $c_{jm} = 1 = -c_{mj}$ such that $\bm \xi_i^{\iin{jm},\delta}  = - \bm \xi_i^{\iin{mj},\delta}$ and $\lambda^{{\iin{jm}},\delta}_i = \lambda^{{\iin{mj}},\delta}_i$ for all $m \in N(j)$ and for $i = 1,\ldots,N_\lambda$. As a consequence, it is legitimate to introduce a numbering of the interfaces $\Gamma^\iin{1},\ldots,\Gamma^\iin{N_\Gamma}$ and to denote the corresponding (unique) vectors of degrees of freedom of the Lagrange multipliers by $\algvec{\lambda}^{\iin{1},\delta}$,\ldots,$\algvec{\lambda}^{\iin{N_\Gamma},\delta}$. Furthermore, it holds that $\mat{B}_{pq}^{\iin{mj},h\delta} = -\int_{\Gamma^\iin{mj}} \boldsymbol \xi_p^{\iin{jm},\delta} \cdot \boldsymbol \varphi^{\is{m},h}_q$. In the following, the coupling matrix between the $i^\text{th}$ interface and the velocity in the $j^\text{th}$ subdomain is also written $\mat{B}^{\iin{i}\is{j},h\delta}$ (to be intended as null matrix, if $\Gamma^\iin{i} \cap \overline{\domain}^\is{j} = \emptyset$).

\paragraph{Discretization of the weak continuity equation} The numerical discretization of Eq.~\eqref{eq:continuity_dd} entails the definition of a preferential side for $\Gamma^{\iin{jm}}$ which determines the computational mesh and the FE basis functions to be used in the evaluation of the integral. However, the direct approximation of Eq.~\eqref{eq:continuity_dd} is problematic because it requires the projection or interpolation of the velocity from one side of the interface to the other (such operation is required in methods such as the mortar method or INTERNODES, for instance). At the continuous level, it is evidently possible to write
\begin{equation}
    \int_{\Gamma^{\iin{jm}}} \boldsymbol \eta \cdot \left (\uvec^{\is{j}} - \uvec^{\is{m}} \right ) = \int_{\Gamma^{\iin{jm}}} \boldsymbol \eta \cdot \uvec^{\is{j}} - \int_{\Gamma^{\iin{mj}}} \boldsymbol \eta \cdot \uvec^{\is{m}}  = 0.
    \label{eq:continuity_separated}
\end{equation}
The advantage of separating the integral into the two contributions on $\Gamma^\iin{jm}$ and $\Gamma^\iin{mj}$ is apparent at the discrete level. Indeed, by substituing $\boldsymbol \eta$ with $\bm \xi_i^{\iin{jm},\delta}$ and $-\bm \xi_i^{\iin{jm},\delta}$ in the two integrals on the right hand side of Eq.~\eqref{eq:continuity_separated}, Eq.~\eqref{eq:continuity_dd} is discretized as
\begin{equation}
    \mat{B}^{\iin{jm},h\delta} \algvec u^\is{j} + \mat{B}^{\iin{mj},h\delta} \algvec u^\is{m} = \algvec{0}.
    \label{eq:continuity_dd_discr}
\end{equation}

\paragraph{Weak imposition of Dirichlet boundary conditions} We recall that we assume that every inlet interface $\Gamma^\iin{m}_\text{in}$ correspond to a single subdomain $\Omega^\is{j}$. The condition $\mathbf u^{\is{j},h} = \mathbf g_m$ on $\Gamma_\text{in}^\iin{m}$ is weakly imposed as
\begin{equation}
\int_{\Gamma_\text{in}^\iin{m}} \boldsymbol \xi_{i,\text{in}}^{\iin{m},\delta} \cdot \left ( \mathbf u^{\is{j},h} - \mathbf g_m \right ) = 0, \quad i = 1,\ldots,N_\lambda.
\label{eq:weakbcs}
\end{equation}
By introducing the interpolation of the boundary data $\mathbf g_m$ onto the FE space spanned by $\{\boldsymbol \varphi_i^{\is{j},h} \}_{i = 1}^{N_u^{\is{z_j},h}}$, the corresponding vector of degrees of freedom $\algvec{g}_m^h$, and $(\mat{B}^{\iin{m}\is{j},h\delta}_\text{in})_{pq} = \int_{\Gamma^\iin{m}_\text{in}} \boldsymbol \xi_{p,\text{in}}^{\iin{m},\delta} \cdot \boldsymbol \varphi_q^{\is{j},h}$, Eq.~\eqref{eq:weakbcs} is approximated as
\begin{equation}
    \mat{B}^{\iin{m}\is{j},h\delta}_\text{in} \algvec{u}^{\is{j},h} - \mat{B}^{\iin{m}\is{j},h\delta}_\text{in} \algvec{g}_m^h = \algvec{0},
\end{equation}
and the first equation in Eq.~\eqref{eq:ns_dd_discr} is accordingly modified as
\begin{equation}
\mat{M}^{\is{j},h} \dot{\algvec{u}}^{\is{j},h} + \left (\mat{K}^{\is{j},h} + \mat{C}^{\is{j},h}(\algvec{u}^{\is{j},h}) \right ) \algvec{u}^{\is{j},h} + (\mat D^{\is{j},h})^\text{T} \algvec{p}^{\is{j},h} + (\mat B^{\is{j},h\delta})^\text{T} \algvec{\lambda}^{\is{j},h} + (\mat B^{\iin{m}\is{j},h\delta}_\text{in})^\text{T} \boldsymbol \lambda^{\iin{m},\delta}_\text{in} = \algvec{f}
\end{equation}
to account for the Lagrange multiplier $\boldsymbol \lambda^{\iin{m},\delta}_\text{in}$.

\paragraph{Assembly of global system} It is possible to arrange the local systems corresponding to the subdomains in the form of a global block system as
\begin{equation}
    \begin{bmatrix}
        \pazocal M^h & \\
         & \\
    \end{bmatrix}
    \begin{bmatrix}
        \dot{\algvec{W}}^h \\
        \dot{\algvec{\Lambda}}^\delta
    \end{bmatrix}
    +
    \begin{bmatrix}
        \pazocal A^h(\algvec{W}^h) & (\pazocal B^{h\delta})^\text{T} \\
        \pazocal B^{h\delta} &
    \end{bmatrix}
    \begin{bmatrix}
        \algvec{W}^h \\
        \algvec{\Lambda}^\delta
    \end{bmatrix}
    =
    \begin{bmatrix}
        \algvec{F}^h \\
        \algvec{G}^{h\delta}
    \end{bmatrix},
    \label{eq:globalsys}
\end{equation}
where
\begin{equation}
\pazocal M^h := \text{diag}\left (
\begin{bmatrix}
\mat{M}^{\is{j},h} & \\
   & \\
\end{bmatrix}
\right )_{j = 1,\ldots,N_\domain},
\quad
\pazocal A^h(\algvec{W}^h) := \text{diag}\left (
\begin{bmatrix}
    \mat{K}^{\is{j},h} + \mat{C}^{\is{j},h}(\algvec u^{\is{j},h}) & (\mat{D}^{\is{j},h})^\text{T} \\
    \mat{D}^{\is{j},h} & \\
\end{bmatrix}
\right )_{j = 1,\ldots,N_\domain},
\end{equation}
$\pazocal B^{h\delta}$ is a block matrix such that $(\pazocal B^{h\delta})_{ij} = [\mat{B}^{\iin{i}\is{j},h\delta}, \mat{O}]$ ($\mat{O}$ being the null matrix) if $i \leq N_\Gamma$ and $(\pazocal B^{h\delta})_{ij} = [\mat{B}^{\iin{i - N_\Gamma}\is{j},h\delta}_\text{in}, \mat{O}]$ otherwise,
\begin{equation}
\algvec{W}^h = [\algvec{w}^{\is{1},h},\ldots,\algvec{w}^{\is{N_\domain},h}], \quad \algvec{\Lambda}^\delta = [\algvec{\lambda}^{\iin{1},\delta},\ldots,\algvec{\lambda}^{\iin{N_\Gamma},\delta},\algvec{\lambda}^{\iin{1},\delta}_\text{in},\ldots,\algvec{\lambda}^{\iin{N_\text{in}},\delta}_\text{in}]
\end{equation}
and $\algvec{F}^h$ and $\algvec{G}^{h\delta}$ are block vectors accounting for the forcing terms and Dirichlet boundary conditions, respectively.

\paragraph{Discretization in time and global nonlinear residual} The discretization in time is performed along the same lines of the discussion in Section~\ref{subsec:numdiscns}. We define $\algvec{Y}^{h\delta} = [\algvec{W}^h, \algvec{\Lambda}^\delta]$ and
\begin{equation}
\pazocal H^h :=
\begin{bmatrix}
    \pazocal M^h & \\
     & \\
\end{bmatrix},
\quad
\mathring{\algvec{F}}^{h\delta}\left (t,\algvec{Y}\right )
:=
\begin{bmatrix}
    \algvec{F}^h(t) \\
    \algvec{G}^{h\delta}(t)
\end{bmatrix}
-
\begin{bmatrix}
\pazocal A^h \left (\algvec{W}^h\right ) & (\pazocal B^{h\delta})^\text{T} \\
\pazocal B^{h\delta} & \\
\end{bmatrix}
\begin{bmatrix}
    \algvec{W}^h \\
    \algvec{\Lambda}^\delta
\end{bmatrix}.
\end{equation}
Then, given $\algvec{Y}^{h\delta}_{k-j+1}$ for $j = 1,\ldots,\sigma$, the solution at timestep $t_{k+1}$ is found by solving
\begin{equation}
    \algvec{R}\left (\algvec{Y}^{h\delta}_{k+1} \right ) := \pazocal{H}^h \algvec{Y}^{h\delta}_{k+1} - \sum_{j = 1}^{\sigma} \alpha_j \pazocal{H}^h \algvec{Y}^{h\delta}_{k-j+1} - \Delta t \beta \mathring{\algvec{F}}^{h\delta}\left (t_{k+1},\algvec{Y}^{h\delta}_{k+1}\right ) = \algvec{0}.
    \label{eq:residualbdf_dd}
\end{equation}

%
\subsection{Efficient solution of the global nonlinear system}
\label{subsec:preconditioner}
Eq.~\eqref{eq:residualbdf_dd} is nonlinear and hence solved using the Newton--Raphson algorithm. In particular, given an initial guess $\algvec{Y}^{(0)}$, the $(l+1)^\text{th}$ iteration of the algorithm for the solution of $\algvec{R}(\algvec{Y}) = \algvec{0}$ is
\begin{equation}
    \algvec{Y}^{(l+1)} = \algvec{Y}^{(l)} - \left(\pazocal J_{\algvec{R}} (\algvec{Y}^{(l)} )\right)^{-1}\algvec{R} (\algvec{Y}^{(l)} ),
    \label{eq:residual_block}
\end{equation}
where $\pazocal J_{\algvec{R}}$ is the tangent matrix of $\algvec{R}$. The stopping criterion is based on a user-provided tolerance $\tau_{\text{NR}}$ and reads $\Vert \algvec{R}(\algvec{Y}^{(l)}) \Vert_2 / \Vert \algvec{R}(\algvec{Y}^{(0)}) \Vert_2 < \tau_{\text{NR}}$.

In order to efficiently solve the linear system in Eq.~\eqref{eq:residual_block} via iterative methods such as GMRES \cite{saad1986gmres}, we need to develop a preconditioner for the tangent matrix $\pazocal J_{\algvec{R}}$. Differentiating Eq.~\eqref{eq:residualbdf_dd} with respect to its only argument yields
\begin{equation}
\pazocal J_{\algvec{R}}(\algvec{Y}) =
\begin{bmatrix}
\pazocal M + \Delta t \beta \pazocal A^h (\algvec{W}) & \Delta t \beta (\pazocal B^{h\delta})^\text{T} \\
\Delta t \beta \pazocal B^{h\delta} &
\end{bmatrix}
=
\begin{bmatrix}
\widetilde{\pazocal{A}}(\algvec{W}) & \widetilde{\pazocal{B}}^\text{T} \\
\widetilde{\pazocal{B}} & \\
\end{bmatrix}.
\label{eq:jacobian}
\end{equation}
Thus, the tangent matrix features a saddle-point structure stemming directly from the original differential problem \weakref{weak:navierstokes_dd}---which is in fact a saddle-point problem. In the remainder of this section, we omit the explicit dependence of $\pazocal J_{\algvec{R}}$ and $\widetilde{\pazocal{A}}$ on $\algvec{Y}$ and $\algvec{W}$, respectively, for the sake of clarity of notation. A possible strategy to design a preconditioner is based on the (exact) decomposition
\begin{equation}
\begin{bmatrix}
\widetilde{\pazocal{A}} & \widetilde{\pazocal{B}}^\text{T} \\
\widetilde{\pazocal{B}} & \\
\end{bmatrix}
=
\begin{bmatrix}
\pazocal{I} & \\
\widetilde{\pazocal B}\widetilde{\pazocal{A}}^{-1} & \pazocal{I}\\
\end{bmatrix}
\begin{bmatrix}
\widetilde{\pazocal{A}} & \\
 & \pazocal{S} \\
\end{bmatrix}
\begin{bmatrix}
\pazocal{I} & \widetilde{\pazocal{A}}^{-1} \widetilde{\pazocal{B}}^\text{T}\\
 & \pazocal{I}\\
\end{bmatrix},
\label{eq:decomposition}
\end{equation}
where $\pazocal S = -\widetilde{\pazocal{B}} \widetilde{\pazocal{A}}^{-1} \widetilde{\pazocal{B}}^\text{T}$ is the Schur complement of $\pazocal J_{\algvec{R}}$. This decomposition is the foundation of several preconditioners for saddle-point systems, such as SIMPLE \cite{segal2010preconditioners} and the nested block preconditioner for blood flow simulations proposed in \cite{liu2019nested}. The solution of a linear system of the form $\pazocal J_{\algvec{R}} \algvec{X} = \algvec{B}$, with $\algvec{X} = [\algvec{X}_w, \algvec{X}_\lambda]$ and $\algvec{B} = [\algvec{B}_w, \algvec{B}_\lambda]$, amounts to solving
\begin{equation}
\begin{bmatrix}
    \widetilde{\pazocal A} & \widetilde{\pazocal B}^\text{T} \\
     & \pazocal{S} \\
\end{bmatrix}
\begin{bmatrix}
    \algvec{X}_w \\
    \algvec{X}_\lambda
\end{bmatrix}
=
\begin{bmatrix}
\widetilde{\pazocal{A}} & \\
 & \pazocal{S} \\
\end{bmatrix}
\begin{bmatrix}
\pazocal{I} & \widetilde{\pazocal{A}}^{-1} \widetilde{\pazocal{B}}^\text{T}\\
 & \pazocal{I}\\
\end{bmatrix}
\begin{bmatrix}
    \algvec{X}_w \\
    \algvec{X}_\lambda
\end{bmatrix}
=
\begin{bmatrix}
    \algvec{B}_w \\
    \algvec{B}_\lambda - \widetilde{\pazocal B} \widetilde{\pazocal{A}}^{-1}\algvec{B}_w\\
\end{bmatrix},
\label{eq:solution_jac}
\end{equation}
hence $\algvec{X}_\lambda = \pazocal S^{-1}(\algvec{B}_\lambda - \widetilde{\pazocal B} \algvec{Z}_w)$ and $\algvec{X}_w = \algvec{Z}_w - \widetilde{\pazocal{A}}^{-1} \widetilde{\pazocal{B}}^\text{T} \algvec{X}_\lambda$, where $\algvec{Z}_w = \widetilde{\pazocal{A}}^{-1}\algvec{B}_w$ can be computed only once for efficiency.

The bottlenecks of the algorithm presented above are evidently the computation of the Schur complement $\pazocal{S}$ and the inversion of matrices $\widetilde{\pazocal{A}}$ and $\pazocal{S}$.

Block $(i,j)$ of the Schur complement explicitly takes the form
\begin{equation}
(\pazocal S)_{ij} = - \Delta t^ 3 \beta^3 \sum_{k = 1}^{N_\domain}
\begin{bmatrix}
\mat B^{\iin{i}\is{k},h\delta} & \mat{O}\\
\end{bmatrix}
\begin{bmatrix}
\mat M^{\is{k},h} / (\Delta t \beta) + \mat K^{\is{k},h} + \mat C^{\is{k},h}(\algvec{u}^{\is{k},h}) & (\mat D^{\is{k},h})^\text{T} \\ \mat D^{\is{k},h} &  \\
\end{bmatrix}^{-1}
\begin{bmatrix}
\mat (\mat B^{\iin{j}\is{k},h\delta})^\text{T} \\
\mat{O}
\end{bmatrix},
\label{eq:schur}
\end{equation}
which is a sum of the contributions of each subdomain $\domain^\is{j}$; we recall that $\mat B^{\iin{i}\is{j},h\delta}$ is not null only if $\Gamma^\iin{i}$ is an interface of $\domain^\is{j}$. As we discussed in Section~\ref{subsec:discretization_dd}, the Lagrange multiplier basis functions required to achieve $h$-convergence are often very few; hence, the Navier--Stokes matrix inverse (or an approximation thereof) can be efficiently applied to every column of $[\mat B^{\iin{j}\is{k},h\delta}, \mat{O}]^\text{T}$, whose number is typically small.
\begin{remark}
    The special structure of the Schur complement in Eq.~\eqref{eq:schur} makes its computation particularly attractive in the context of high-performance computing, assuming that each subdomain be assigned to one or few processors. An effective partition strategy should take into account the number of degrees of freedom associated with each geometrical building block in order to preserve the load balancing among the set of computational nodes. At the same time, the fact that the coupling conditions do not require to explicitly transfer information from one side of the interface to the other can be exploited to limit the amount of processor communications. In this paper, we only consider a serial implementation, but future directions of the current work include investigations on the parallel performance of the preconditioner.
    \label{remark:schur_parallel}
\end{remark}
\begin{remark}
    The computation of the Schur complement $\pazocal{S}$, whose cost scales linearly with the number of columns of $\widetilde{\pazocal B}$ (i.e. the total number of Lagrange multiplier basis functions), is in fact the most expensive step of the application of the preconditioner discussed in this section. In our numerical experiments, however, we observed that significant performance improvements are achieved by reusing the same Schur complement for $n$ consecutive applications of the preconditioner. This strategy does not significantly affect the number of GMRES iterations, if $n$ is small enough (in our case, $n \sim 20$).
\end{remark}

Let us now address the application of the inverses of $\widetilde{\pazocal{A}}$ and $\pazocal{S}$. Since $\pazocal{S} \in \reals^{N_\Gamma N_\lambda \times N_\Gamma N_\lambda}$, the Schur complement is typically small and is inexpensively inverted by solving the linear system either directly or via iterative methods (complemented with standard preconditioners such as multigrid or ILU). Matrix $\widetilde{\pazocal{A}}$ features a block diagonal structure in which each diagonal block is itself a saddle-point system and inverting $\widetilde{\pazocal{A}}$ is therefore equivalent to solving linear systems which are local to each subdomain $\domain^\is{j}$. Employing, instead of $\widetilde{\pazocal{A}}^{-1}$, a suitable approximation thereof, gives rise to different suitable preconditioners.

\begin{figure}
    \tikzsetnextfilename{preconditioner}
    \centering
    \setlength
    \figureheight{0.4\textwidth}
    \setlength
    \figurewidth{\textwidth}
    \input{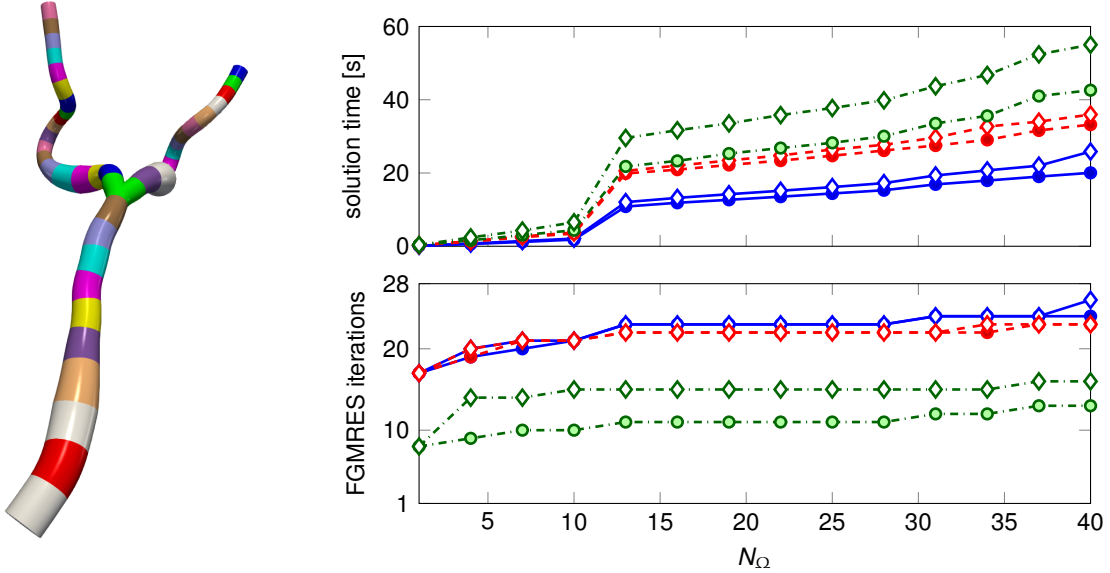}
    \caption{Solution time of a single linear system preconditioned as discussed in Section~\ref{subsec:preconditioner} (right-top) and corresponding number of FGMRES iterations (right-bottom), as functions of number of subdomains $N_\domain$ in the geometry of the aorta and iliac arteries (left). Blue solid lines: $\widetilde{\pazocal{A}}^{-1}$ approximated with a single application of SIMPLE for each subdomain; red dashed lines and blue dash-dotted lines: $\widetilde{\pazocal{A}}^{-1}$ solved with GMRES and tolerance 5e--1 and 1e--2, respectively; circles: $N_\lambda = 3$; diamonds = $N_\lambda = 84$.}
    \label{fig:preconditioner}
\end{figure}

In this work, we choose to approximate every block diagonal matrix of $\widetilde{\pazocal{A}}^{-1}$ by considering a single application of the SIMPLE preconditioner, both in the computation of the Schur complement shown in Eq.~\eqref{eq:schur} and in the application of the preconditioner. We recall that SIMPLE is based on Eq.~\eqref{eq:decomposition} applied to the Navier--Stokes equations and that the approximated inverse of the top left diagonal block is performed by extracting and inverting its diagonal \cite{segal2010preconditioners}. In Fig.~\ref{fig:preconditioner}, we show the robustness of our preconditioner with respect to the number of blocks and to the number of basis functions for the Lagrange multipliers per interface $N_\lambda$. The considered geometry is that of the aorta and the illiac arteries in Fig.~\ref{fig:preconditioner} (left). The blocks are sequentially added starting from the inlet (for this reason, we remark that the size of the system increases proportionally with the number of blocks). We compare the preconditioner performance with that achieved by inverting every block in $\widetilde{\pazocal{A}}$ with GMRES and relatively large tolerances (5e--1 and 1e--2). We remark that, as the preconditioner in the latter approach varies at each iteration, we are obliged to employ flexible GMRES (FGMRES) \cite{saad1993flexible}. If the local systems are solved exactly, the preconditioner is in fact the original global matrix, as Eq.\eqref{eq:decomposition} is an exact decomposition. For this reason, solving the local linear systems with GMRES leads to a better performance in terms of number of iterations. However, approximating each local inverse with SIMPLE is more efficient in terms of solution time, as each FGMRES iteration is less computational expensive. We conclude by observing that the increase in solution time occurring at $N_\Omega = 13$ is due to the introduction of the bifurcation---which is composed of a larger number of elements than the other blocks---in the set of considered subdomains.

\section{The reduced Navier--Stokes equations on modular domain-decompositions of arteries}
\label{sec:reducedrb}
\begin{figure}
    \tikzsetnextfilename{offline}
    \centering
    \setlength
    \figureheight{\textwidth}
    \setlength
    \figurewidth{\textwidth}
    \input{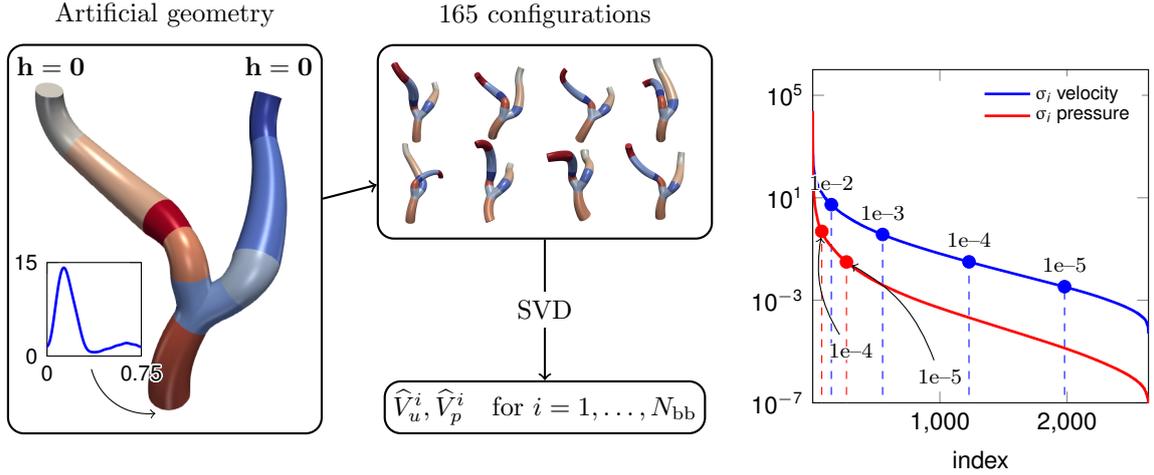}
    \caption{Offline phase on an artificial geometry featuring $N_\domain = 9$ subdomains and $N_\text{bb} = 4$ building blocks (left) and singular values decay for the velocity and pressure reduced basis built on the bifurcation (right). The snapshots are constructed by solving---on 165 deformed configurations---the flow problem obtained by imposing the flow rate $Q$ depicted in the box on the left ($y$-axis: $Q$ [cc/s], $x$-axis: $t$ [s]) at the inlet, and $\mathbf{h} = \boldsymbol{0}$ on both outflows, in the range $t = (0,0.3)$. The colored dots on the plot on the right refer to different values of POD tolerances $\epsilon_u$ and $\epsilon_p$.}
    \label{fig:offline}
\end{figure}

\label{sec:nsmodularrb}
As discussed in Section~\ref{sec:modular}, the subdomains in the target geometry $\domain_\text{m}$ are obtained from the parametrized geometrical deformation of a number of building blocks $\bbref^\is{i}$, $i = 1,\ldots,N_\text{bb}$. In this paper, these are a model symmetric bifurcation (B), and straight tubes with aspect ratios length/diameter 1:1 (T1), 1:2 (T2) and 1:3 (T3). The offline phase of our reduced order model algorithm consists of defining reduced basis functions in each of these building blocks $\bbref^\is{i}$. The snapshots are collected from a single decomposed ``artificial'' geometry $\domain_\text{m} = \bigcup_{j = 1}^{N_{\domain}} \domain^\is{j}$ by sampling the geometrical parameters $\param^1,\ldots,\param^{N_\domain}$ describing each subdomain from uniform distributions centered on the values characterizing the original configuration, as depicted in Fig.~\ref{fig:offline}. The snapshots are found by solving a flow problem with $\rho_\text{f} = 1.06$ $\text{gr}\,\text{cm}^{-3}$, $\mu_\text{f} = 0.04$ $\text{gr}\,\text{cm}^{-1}\,\text{s}^{-1}$, the imposed inflow flow rate shown in Fig.~\ref{fig:offline} (in the box on the left) with a parabolic profile and homogeneous Neumann conditions on the outlets on 165 random configurations of the artificial geometry. There exist other equally valid possibilities to generate the database of snapshots. For example, these could be taken by solving flow problems on a collection of target geometries. This approach allows us to avoid issues related to the random sampling of the geometrical parameters---e.g. physiological feasibility of the resulting global geometry---but requires the aid of an automatic algorithm for the decomposition to be efficient. The development of such an algorithm is one of the possible future extensions of the present work. The simulations are run from $t_0 = 0$~s to $T = 0.3$~s with a BDF scheme of order $\sigma = 2$ and $\Delta t = 2.5 \times 10^{-3}$~s. The initial condition at $t_0$ is computed by gradually increasing the inflow flow rate profile at the inlet by the law $Q(t) = Q_0[1 - \cos((t - t_0^\text{ramp}) \pi / (t_0 - t_0^\text{ramp}))] / 2$, $Q_0$ being the desired flow rate at time $t_0$, from $t = t_0^\text{ramp} = -2 \times 10^{-2}$ s to $t = t_0$.  For the discretization of the Lagrange multipliers on each interface, we employ the set of basis functions $\widehat{\Xi}_n$ with $n = 5$, which corresponds to $N_\lambda = 63$ basis functions. We remark that the artificial geometry is not included in the configurations used for the snapshots generation and is considered in Section~\ref{subsec:numres_artificial} to assess the performance of the method.

\begin{figure}
\centering
\tikzsetnextfilename{basisfunctions}
\begin{tikzpicture}
    \node[inner sep=0pt] (bif) at (0,0)
    {\includegraphics[width = 0.12\textwidth,trim={4cm 0cm 3.4cm 16cm},clip]{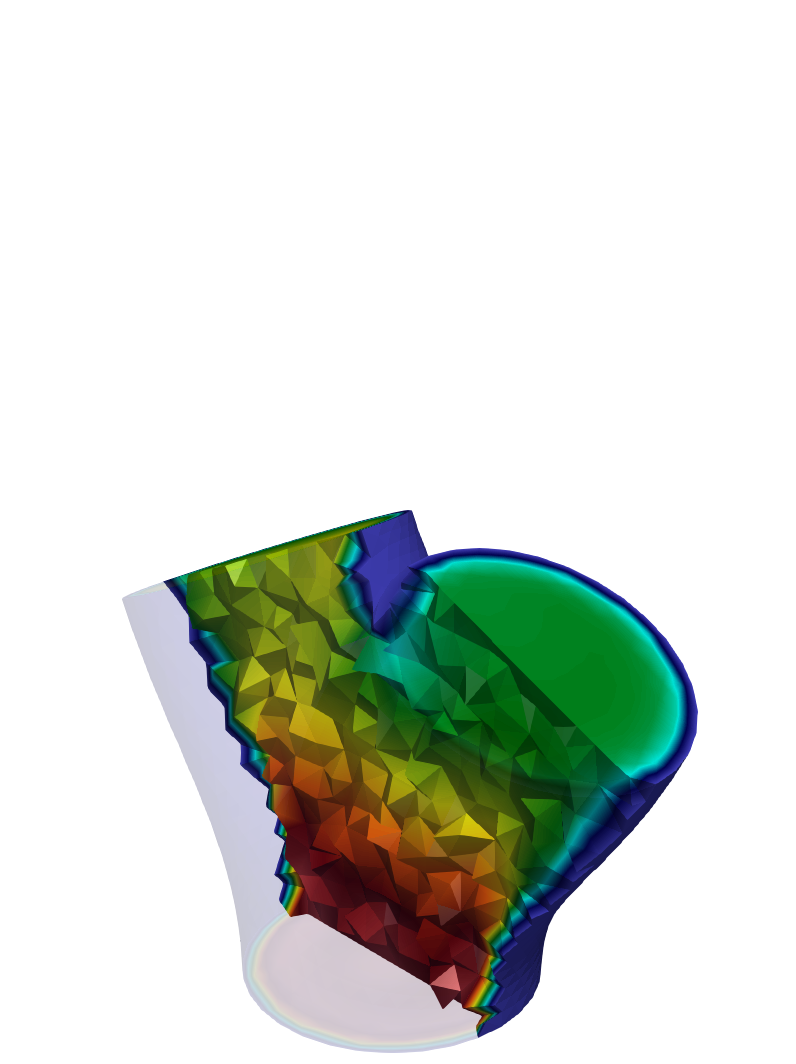}
     \includegraphics[width = 0.12\textwidth,trim={4cm 0cm 3.4cm 16cm},clip]{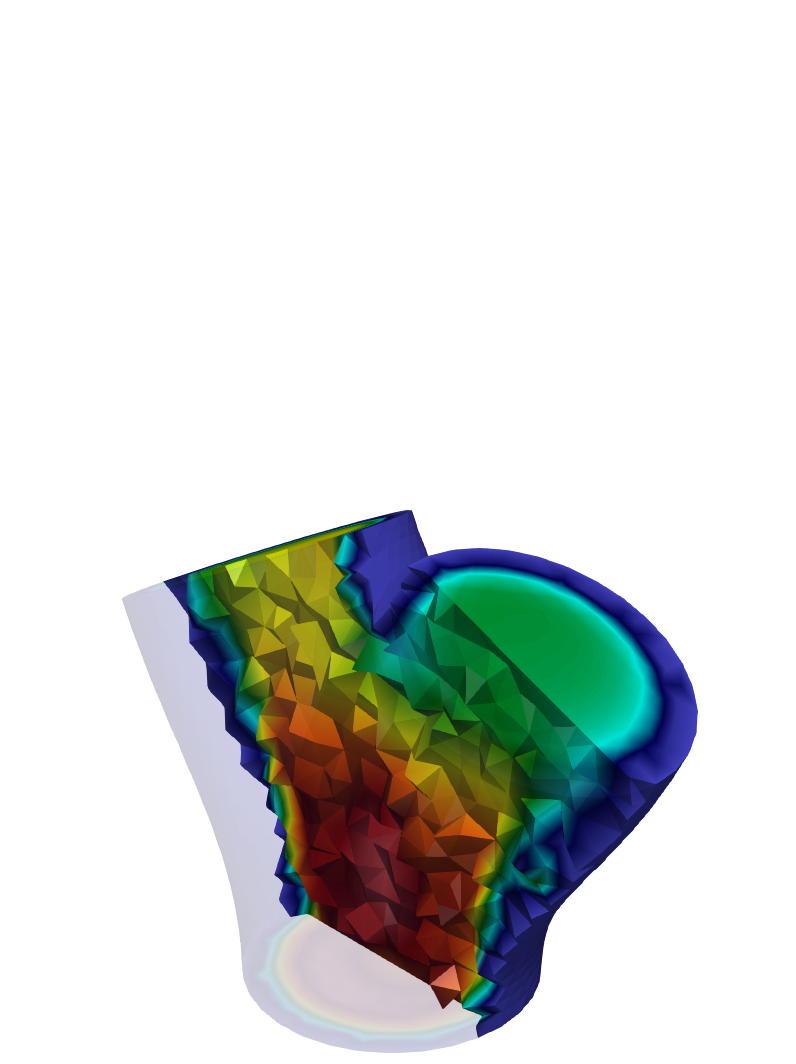}
     \includegraphics[width = 0.12\textwidth,trim={4cm 0cm 3.4cm 16cm},clip]{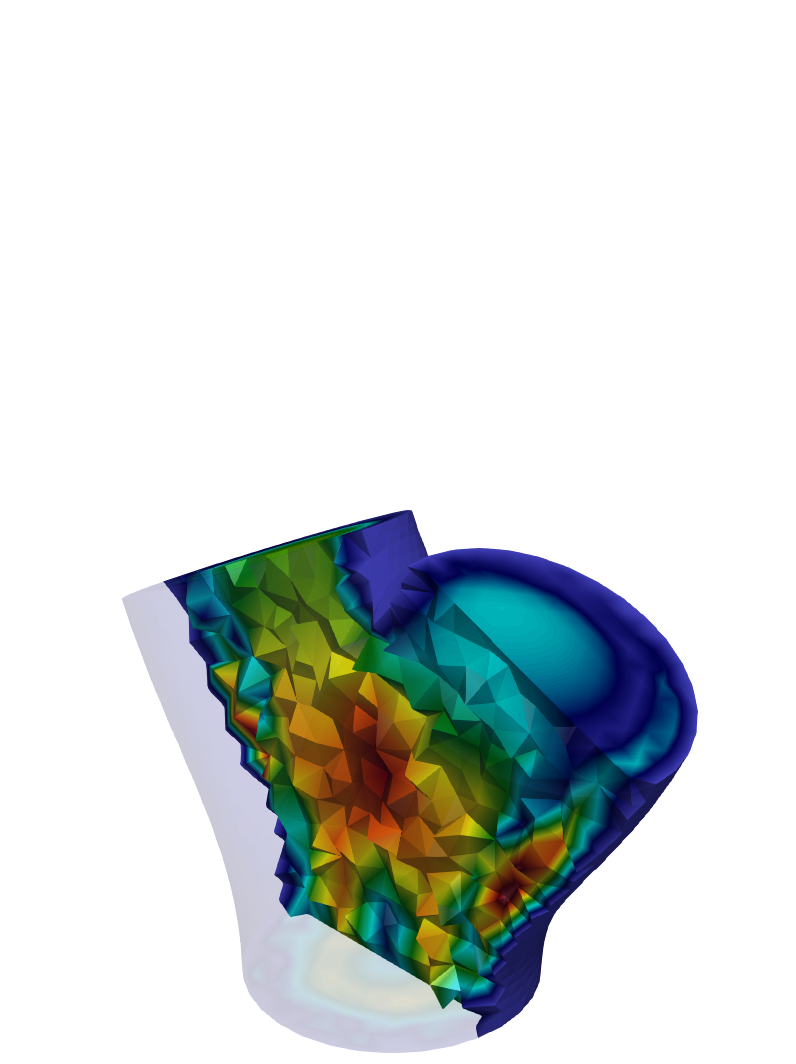}
     \includegraphics[width = 0.12\textwidth,trim={4cm 0cm 3.4cm 16cm},clip]{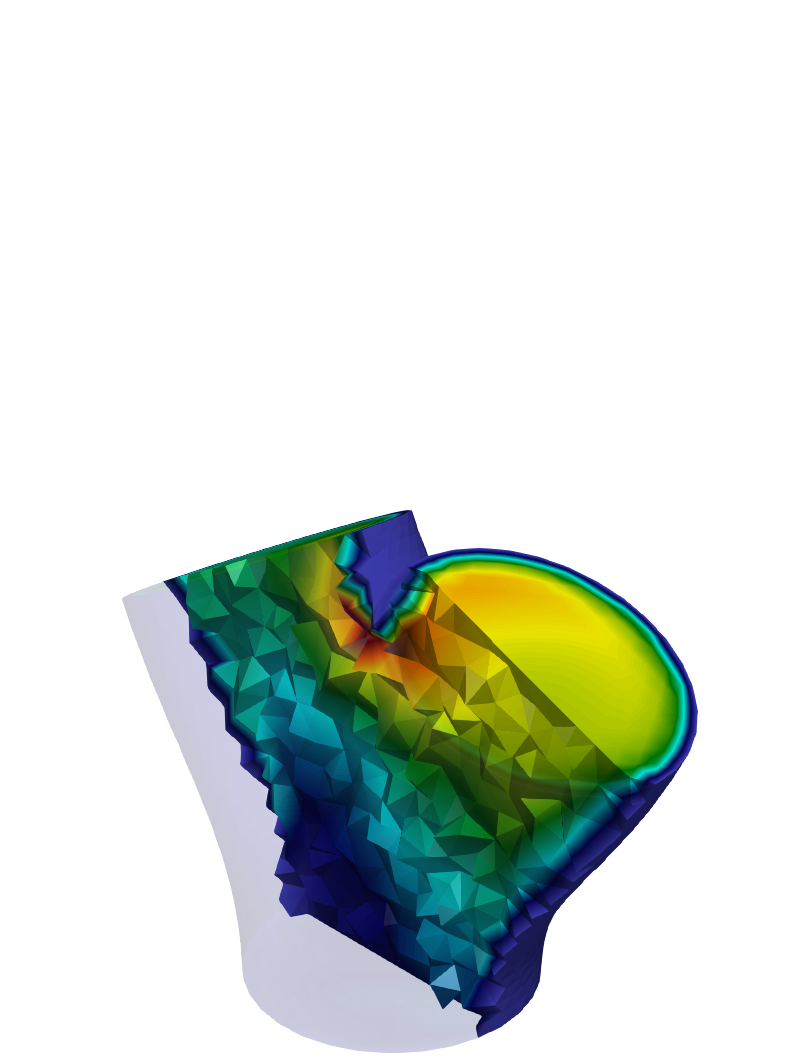}
     \includegraphics[width = 0.12\textwidth,trim={4cm 0cm 3.4cm 16cm},clip]{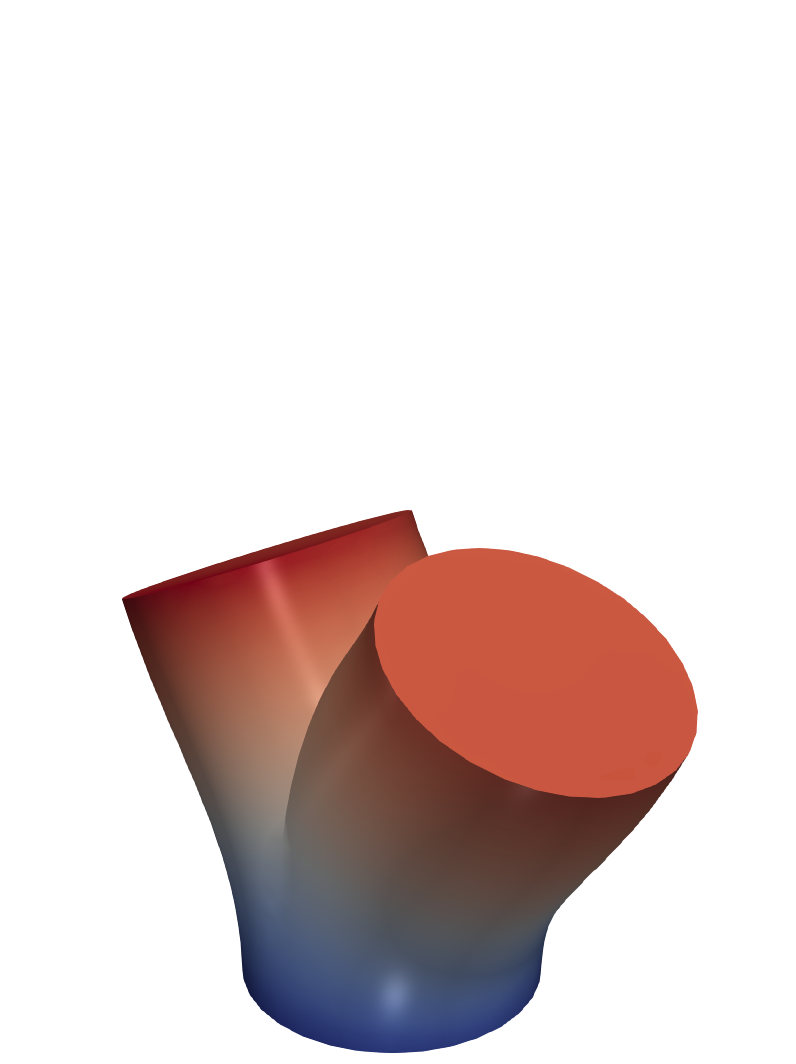}
     \includegraphics[width = 0.12\textwidth,trim={4cm 0cm 3.4cm 16cm},clip]{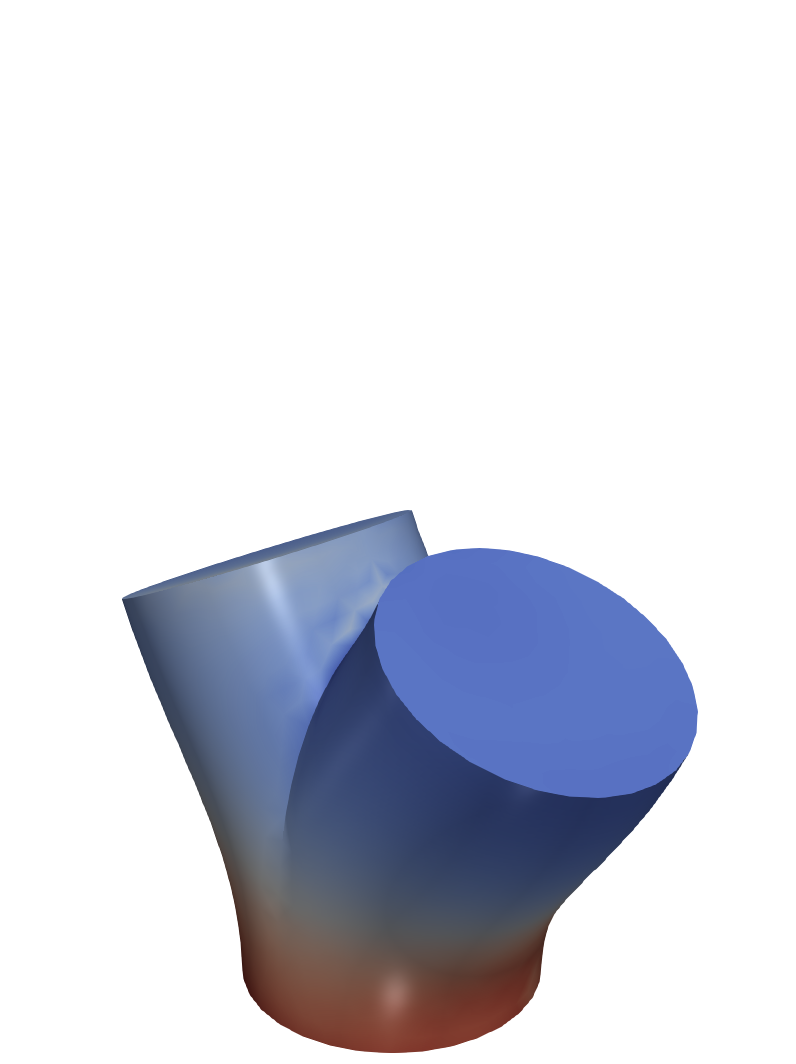}
     \includegraphics[width = 0.12\textwidth,trim={4cm 0cm 3.4cm 16cm},clip]{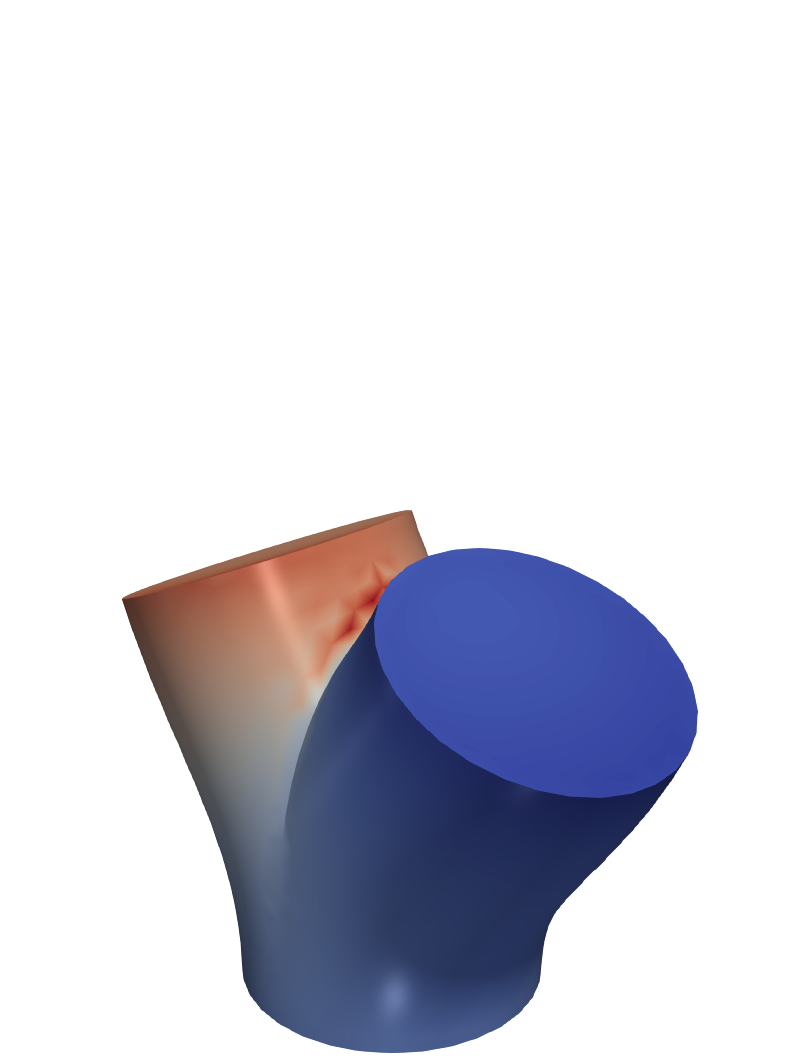}
     \includegraphics[width = 0.12\textwidth,trim={4cm 0cm 3.4cm 16cm},clip]{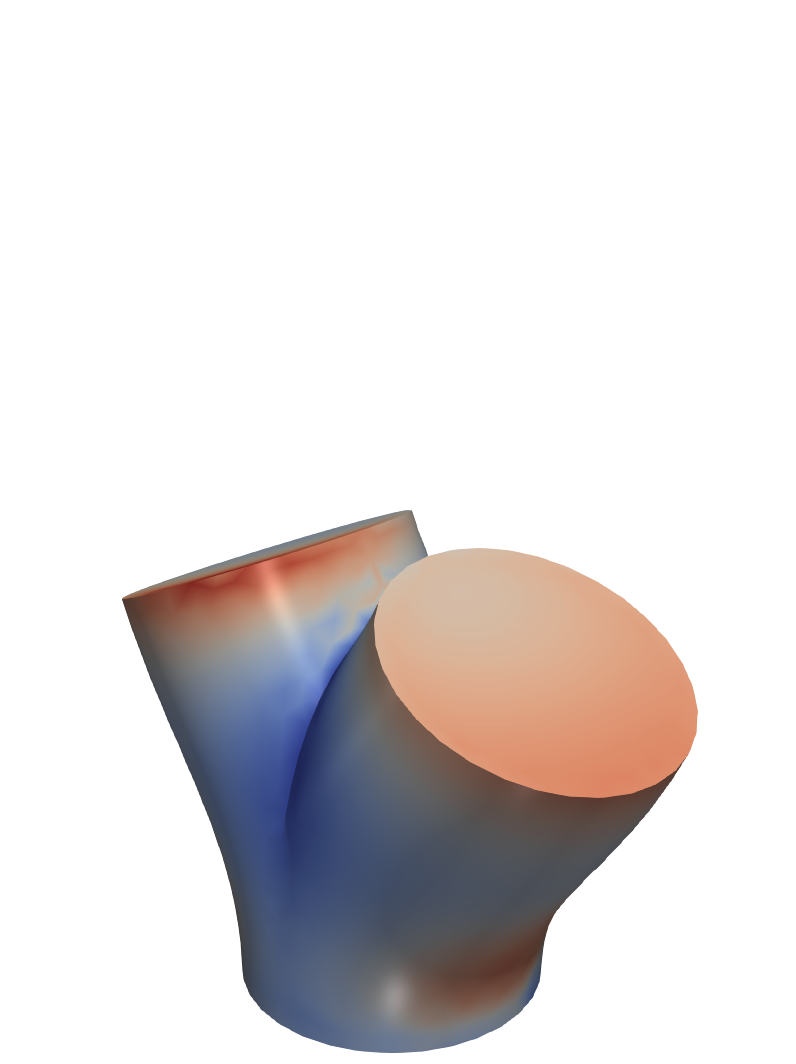}};
    \node[inner sep=0pt, below=0.5cm of bif] (t1x1)
    {\includegraphics[width = 0.12\textwidth,trim={4cm 0cm 3.4cm 21cm},clip]{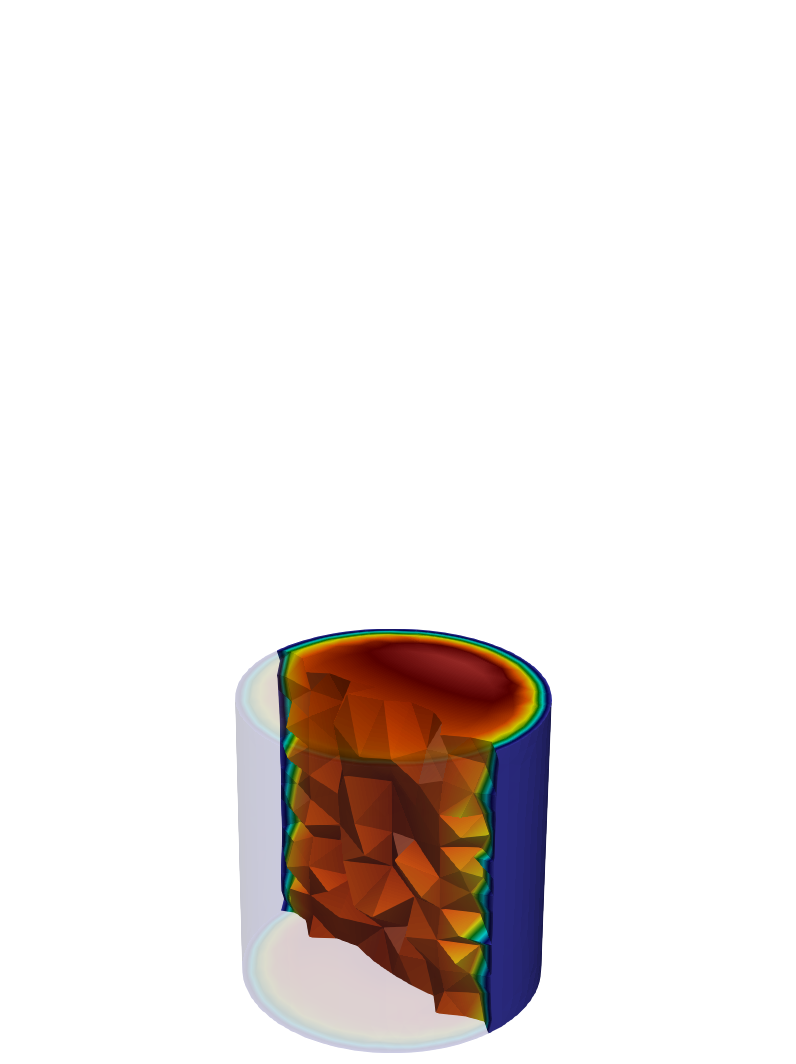}
     \includegraphics[width = 0.12\textwidth,trim={4cm 0cm 3.4cm 21cm},clip]{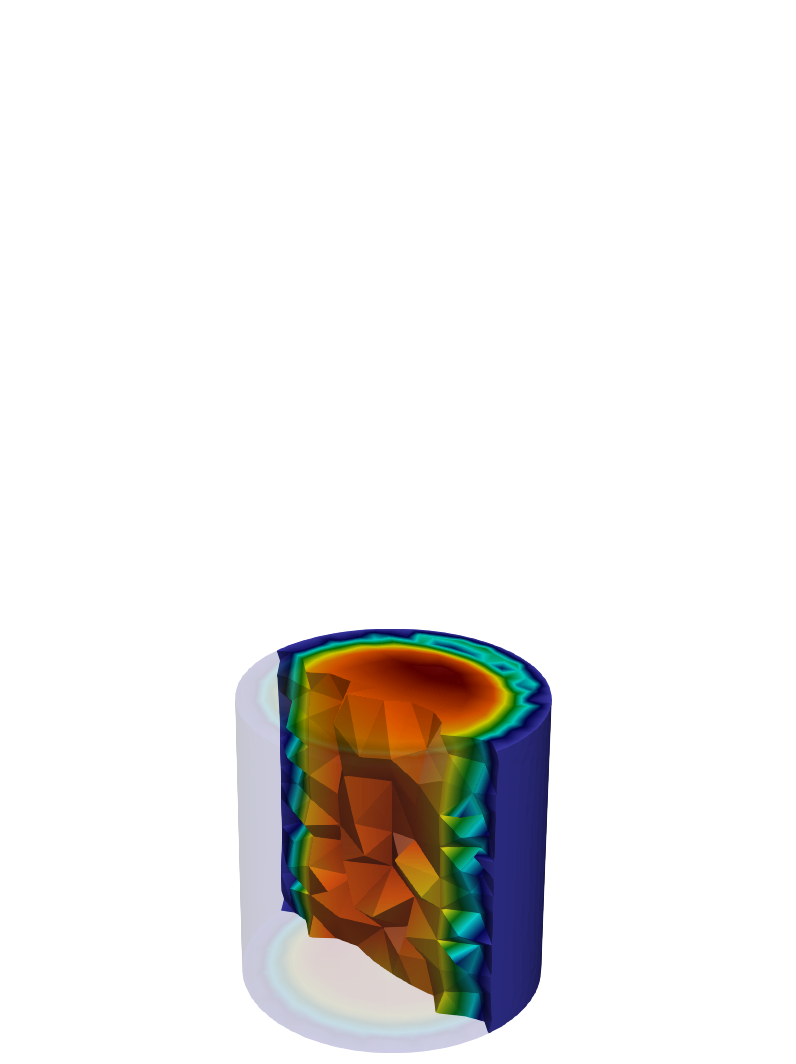}
     \includegraphics[width = 0.12\textwidth,trim={4cm 0cm 3.4cm 21cm},clip]{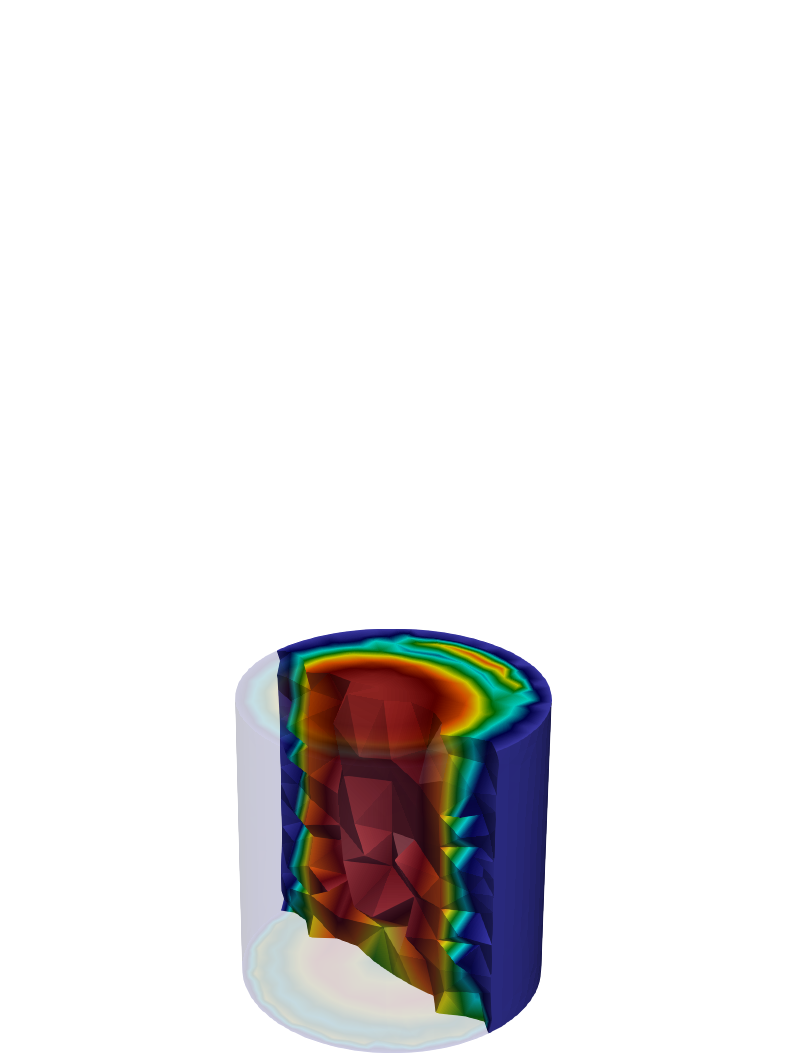}
     \includegraphics[width = 0.12\textwidth,trim={4cm 0cm 3.4cm 21cm},clip]{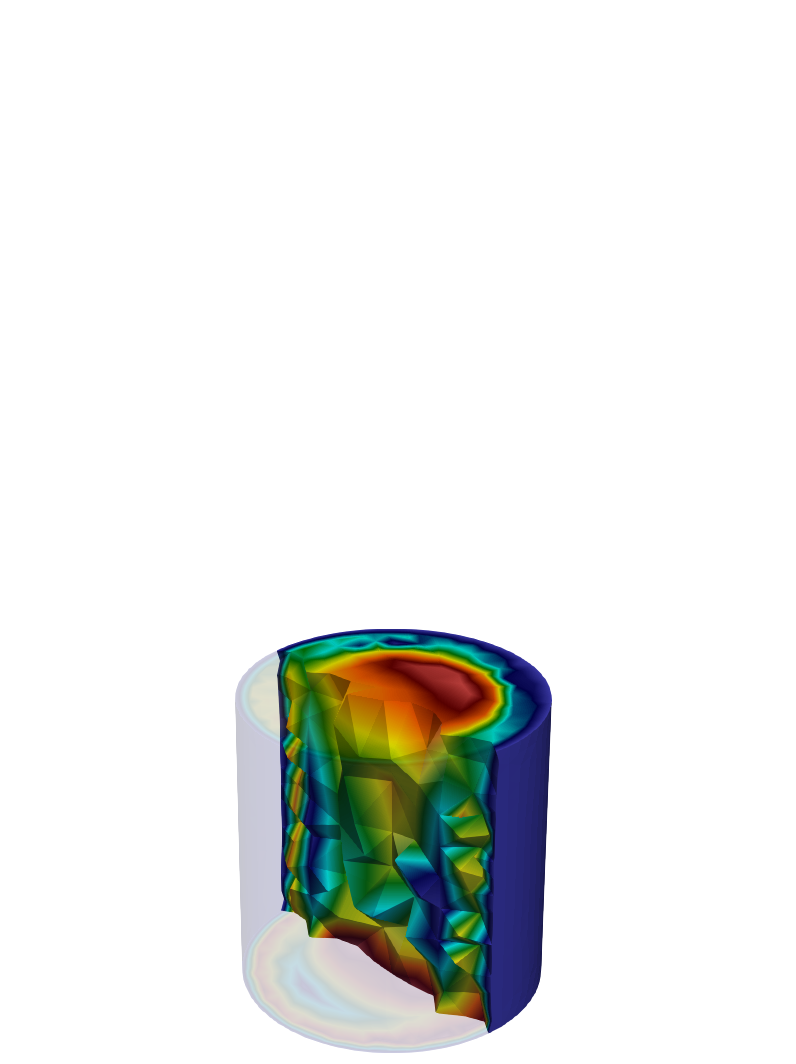}
     \includegraphics[width = 0.12\textwidth,trim={4cm 0cm 3.4cm 21cm},clip]{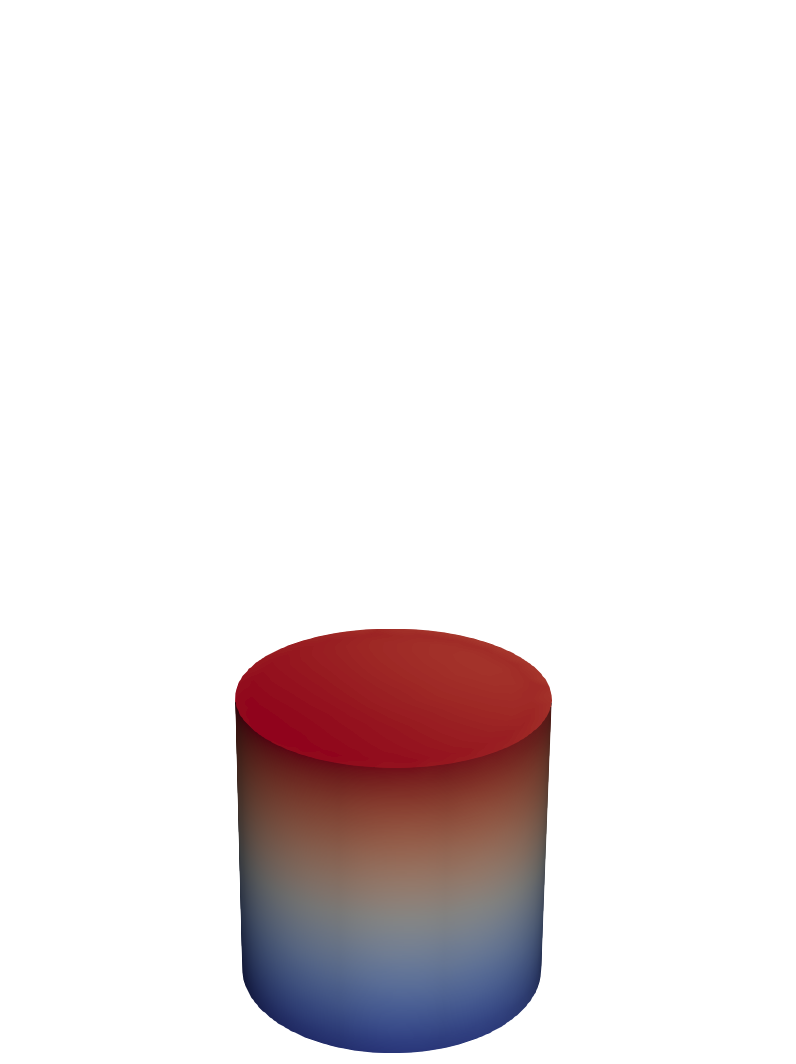}
     \includegraphics[width = 0.12\textwidth,trim={4cm 0cm 3.4cm 21cm},clip]{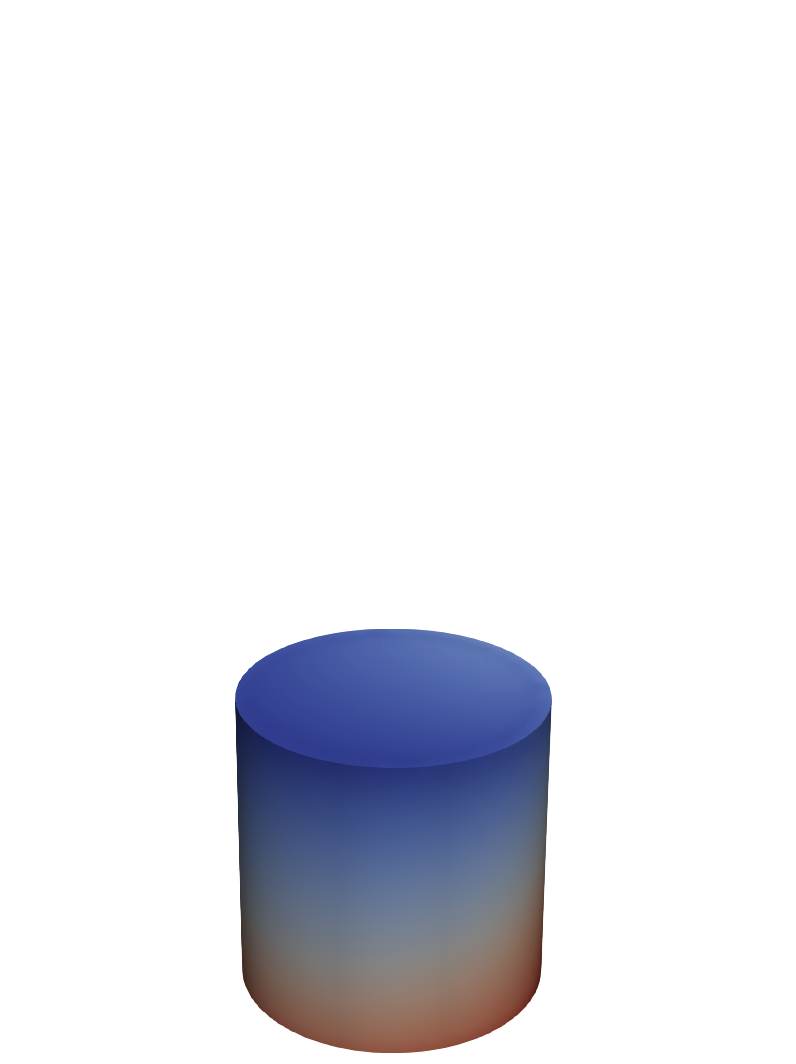}
     \includegraphics[width = 0.12\textwidth,trim={4cm 0cm 3.4cm 21cm},clip]{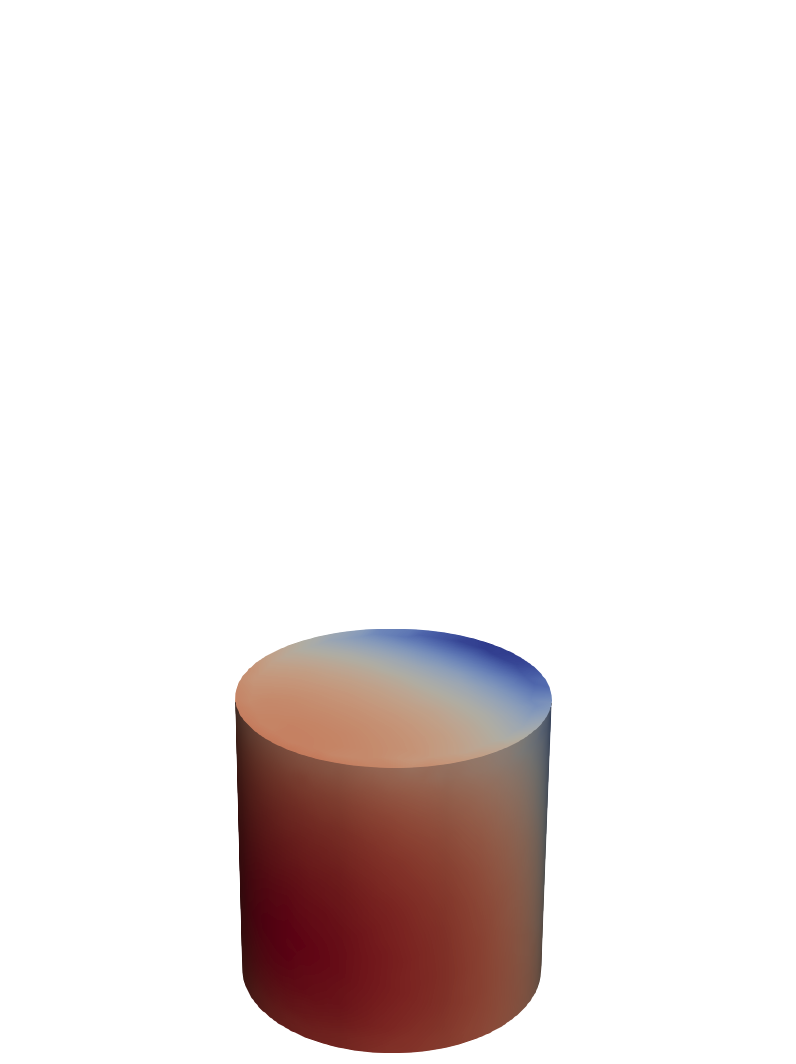}
     \includegraphics[width = 0.12\textwidth,trim={4cm 0cm 3.4cm 21cm},clip]{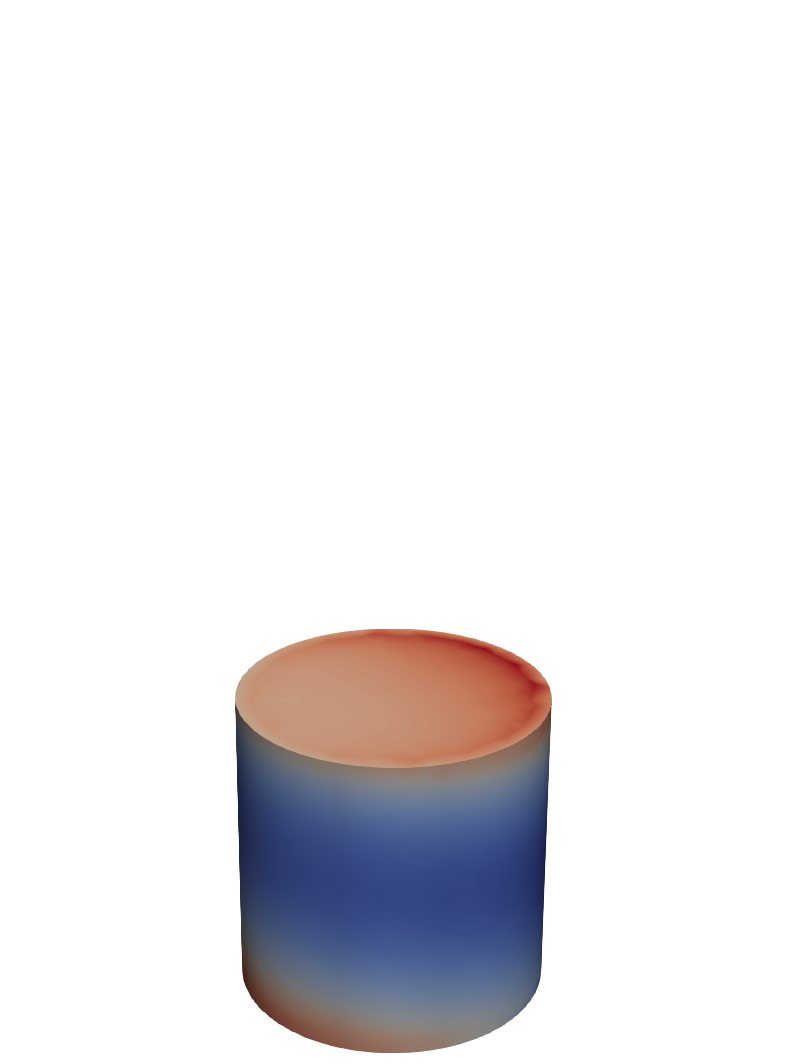}};
     \node[inner sep=0pt,below=0.5cm of t1x1] (t1x2)
     {\includegraphics[width = 0.12\textwidth,trim={4cm 0cm 3.4cm 10cm},clip]{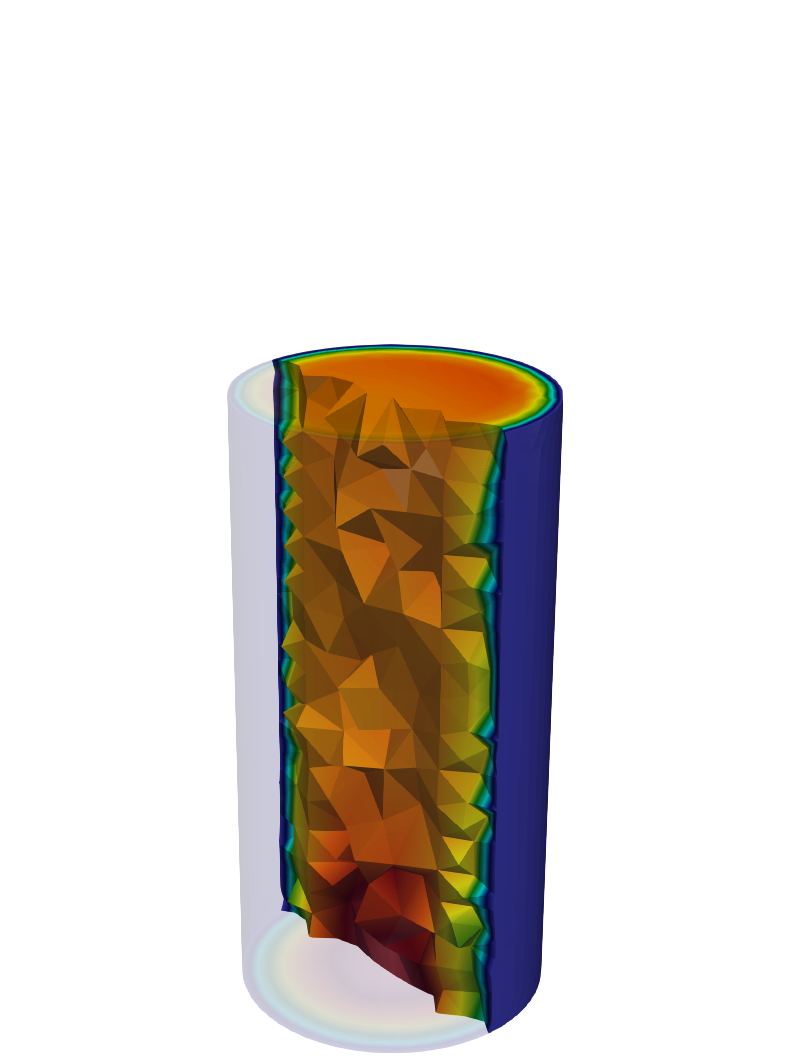}
      \includegraphics[width = 0.12\textwidth,trim={4cm 0cm 3.4cm 10cm},clip]{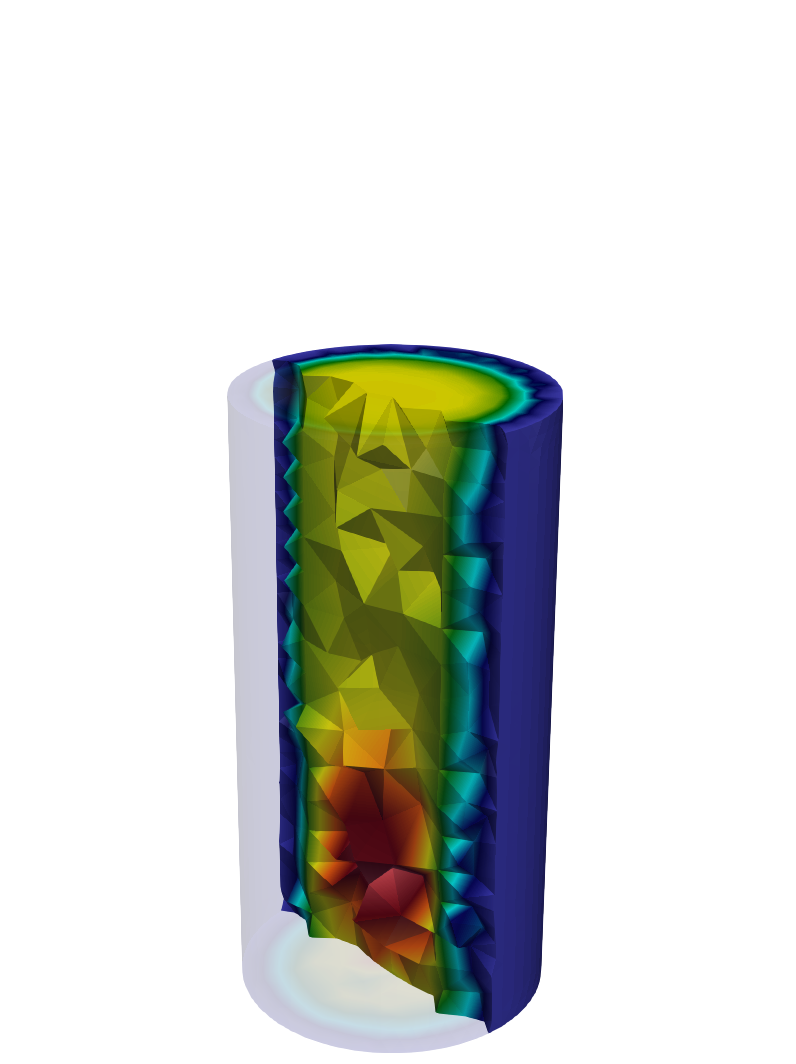}
      \includegraphics[width = 0.12\textwidth,trim={4cm 0cm 3.4cm 10cm},clip]{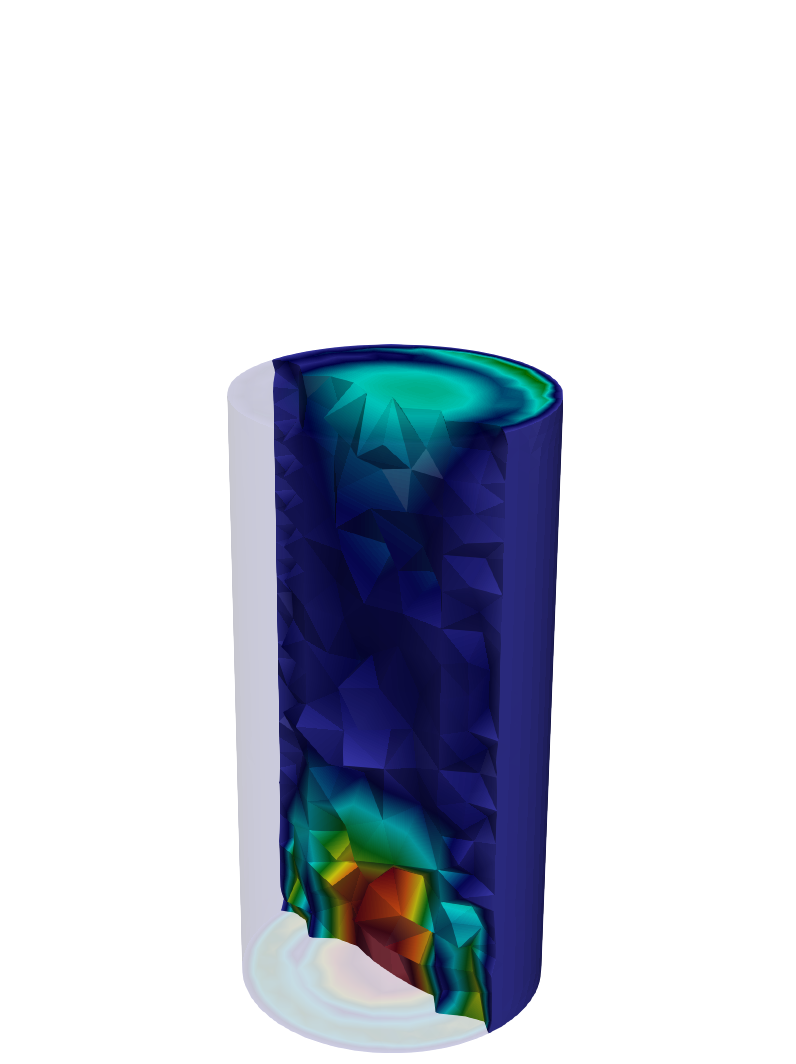}
      \includegraphics[width = 0.12\textwidth,trim={4cm 0cm 3.4cm 10cm},clip]{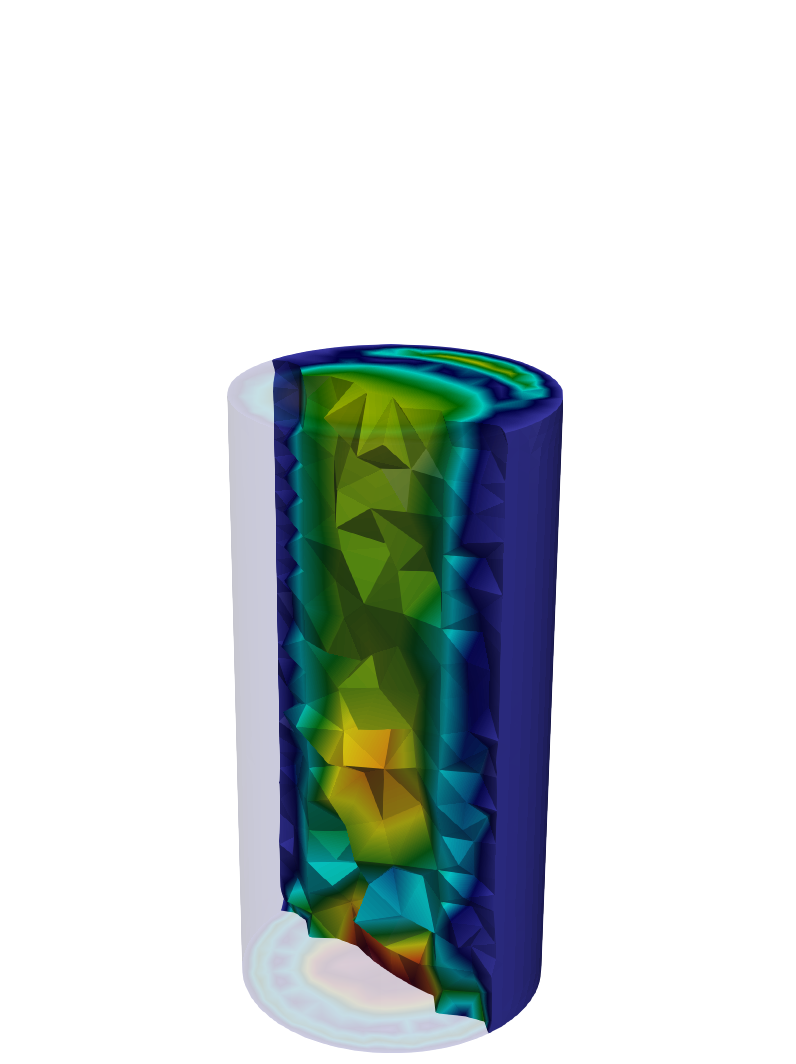}
      \includegraphics[width = 0.12\textwidth,trim={4cm 0cm 3.4cm 10cm},clip]{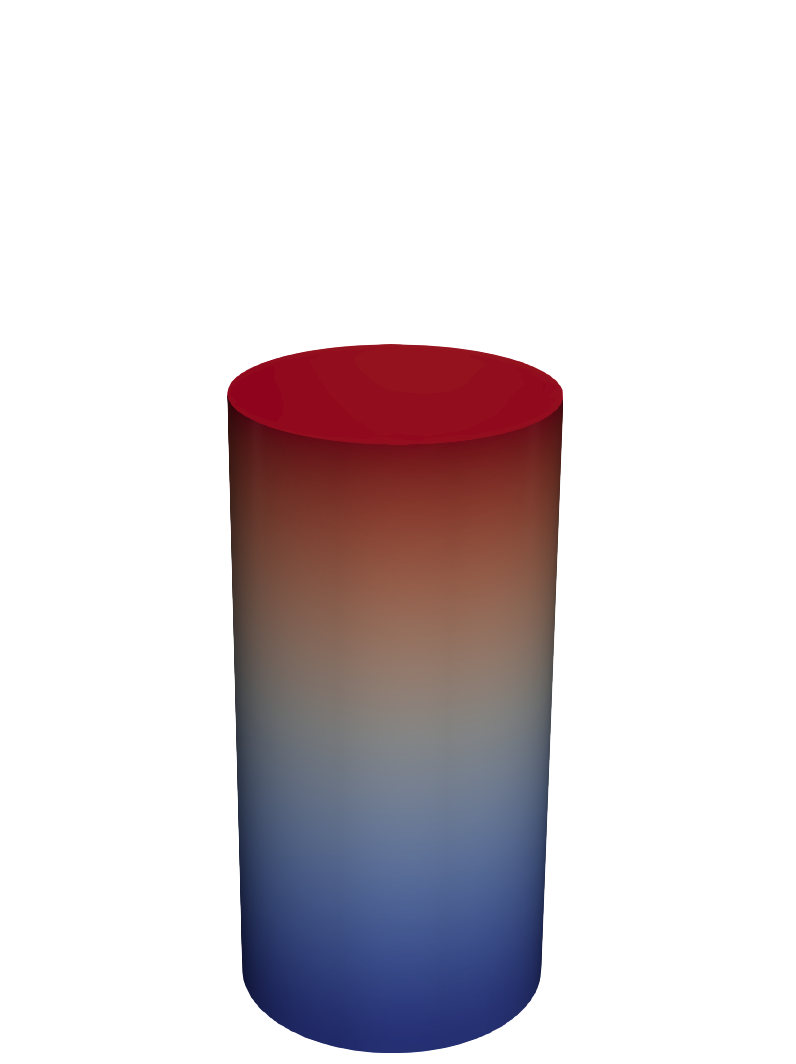}
      \includegraphics[width = 0.12\textwidth,trim={4cm 0cm 3.4cm 10cm},clip]{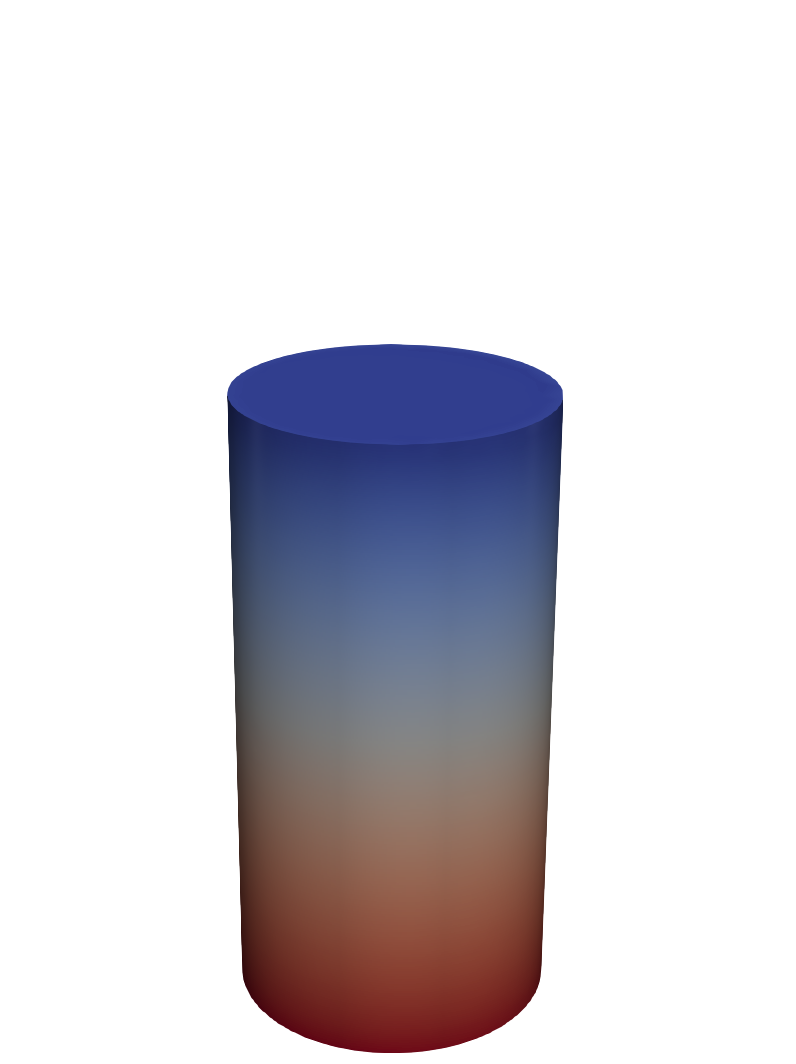}
      \includegraphics[width = 0.12\textwidth,trim={4cm 0cm 3.4cm 10cm},clip]{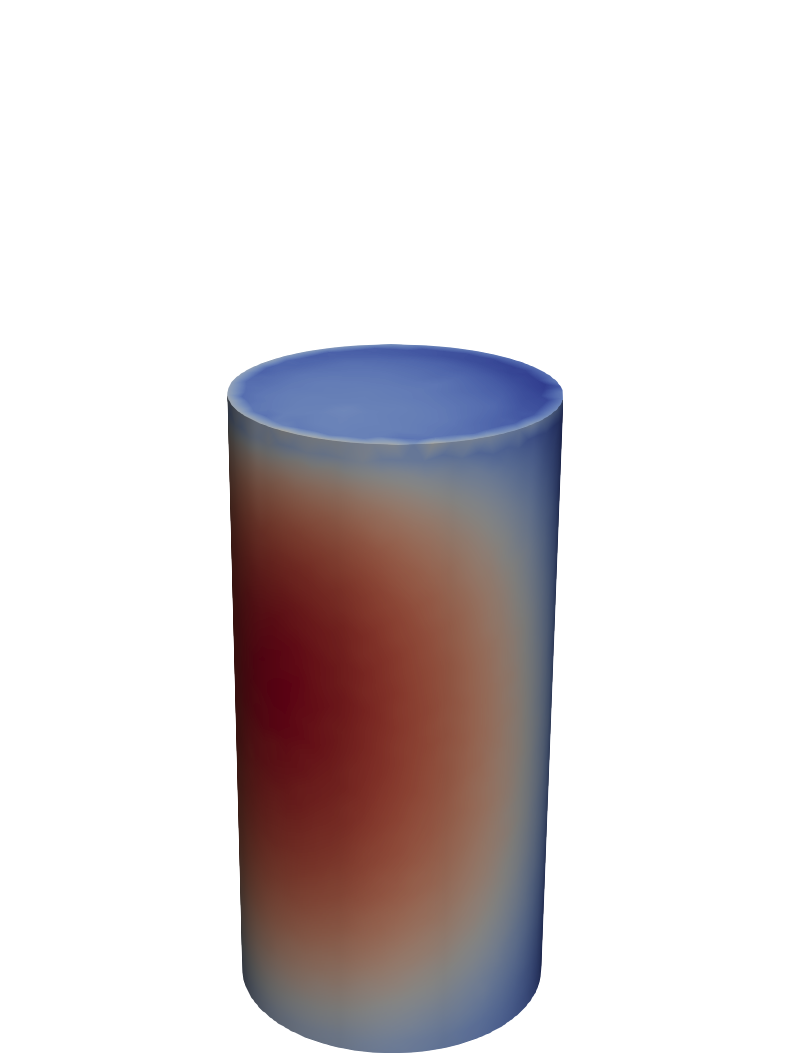}
      \includegraphics[width = 0.12\textwidth,trim={4cm 0cm 3.4cm 10cm},clip]{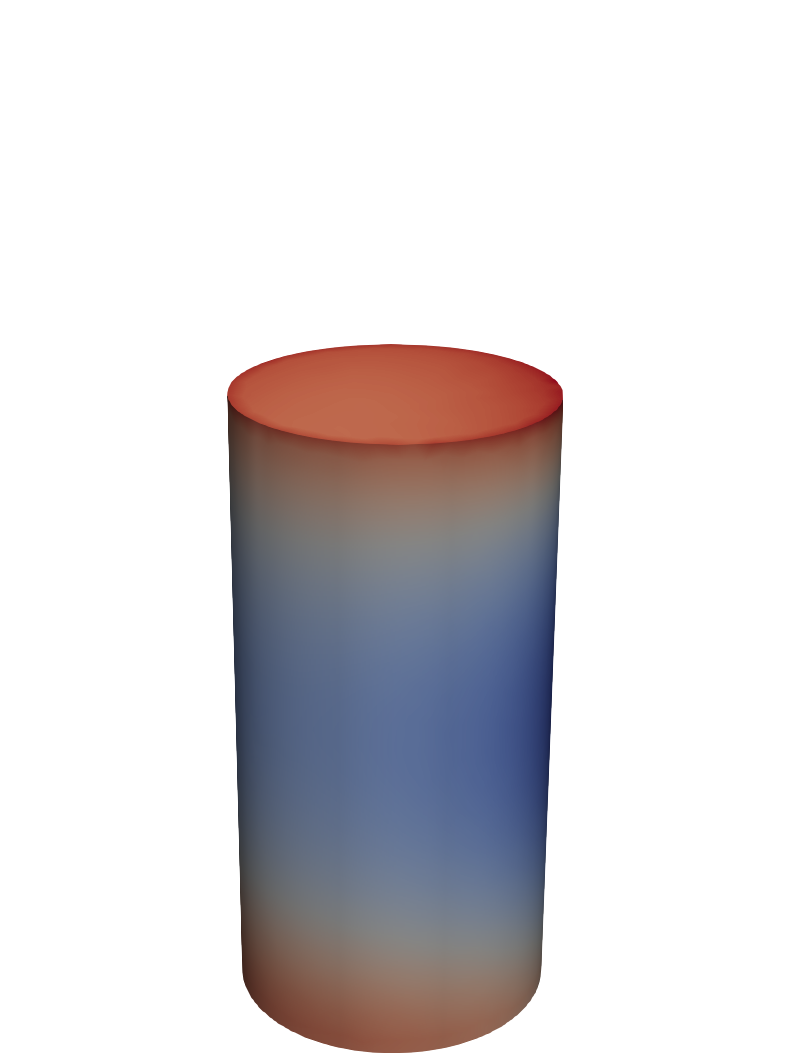}};
      \node[inner sep=0pt,below=0.5cm of t1x2] (t1x3)
      {\includegraphics[width = 0.12\textwidth,trim={4cm 0cm 3.4cm 0cm},clip]{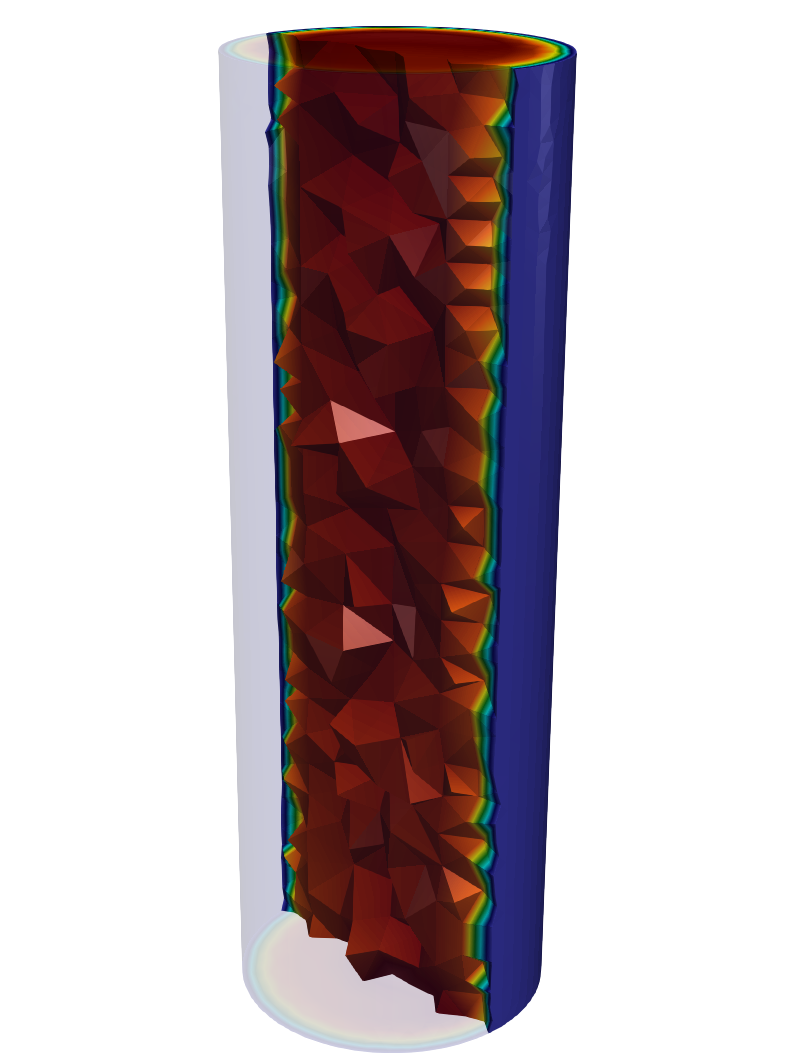}
       \includegraphics[width = 0.12\textwidth,trim={4cm 0cm 3.4cm 0cm},clip]{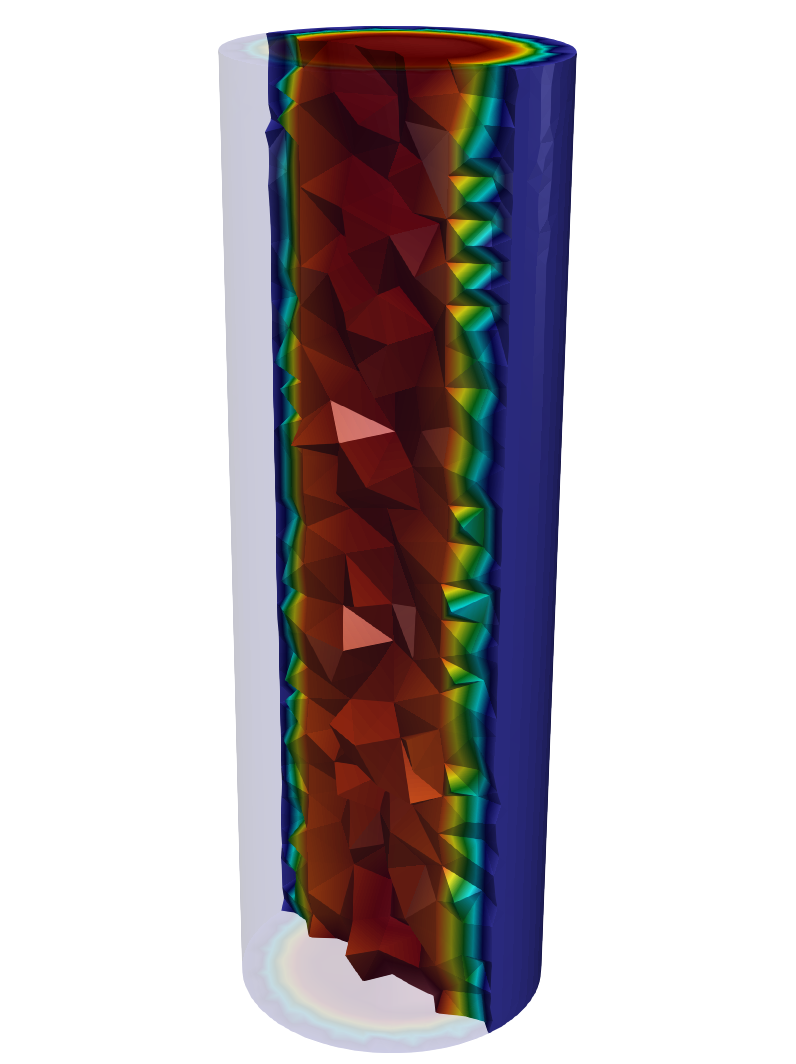}
       \includegraphics[width = 0.12\textwidth,trim={4cm 0cm 3.4cm 0cm},clip]{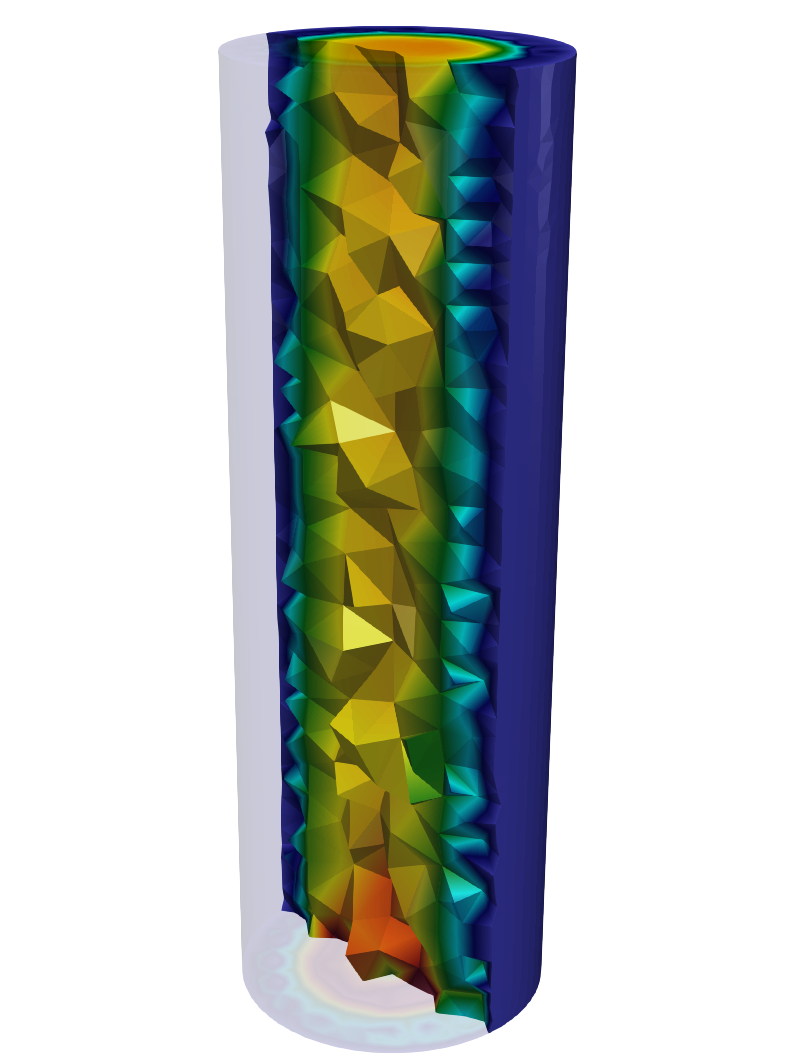}
       \includegraphics[width = 0.12\textwidth,trim={4cm 0cm 3.4cm 0cm},clip]{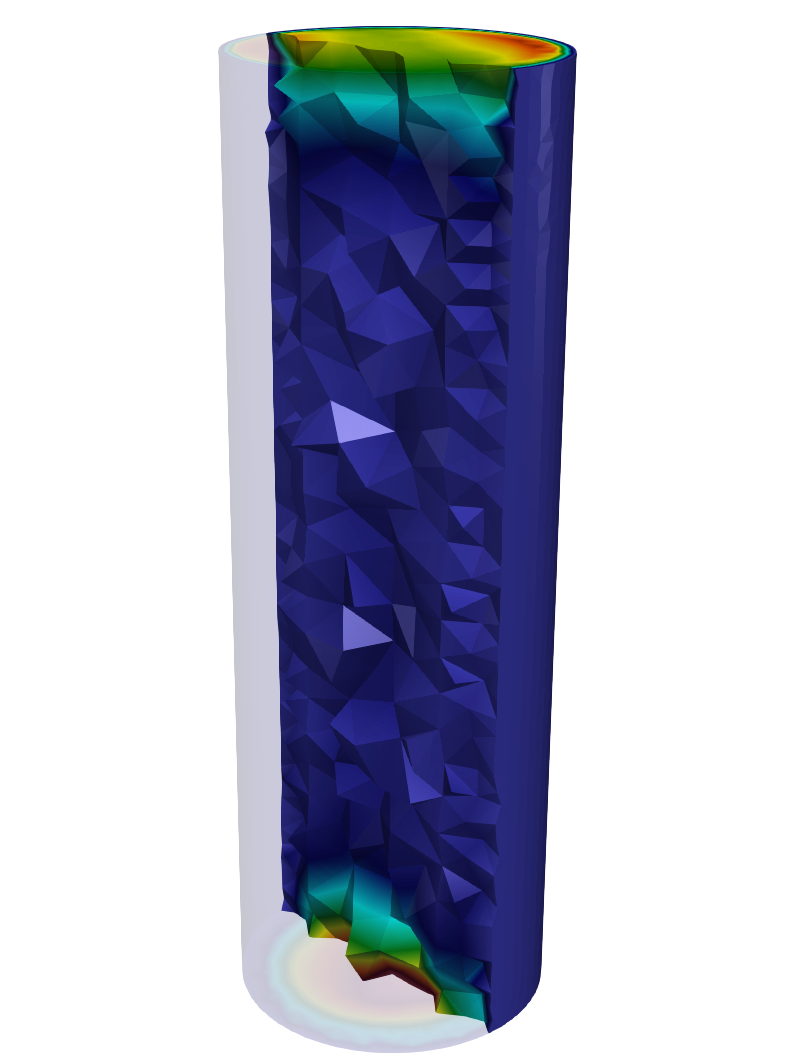}
       \includegraphics[width = 0.12\textwidth,trim={4cm 0cm 3.4cm 0cm},clip]{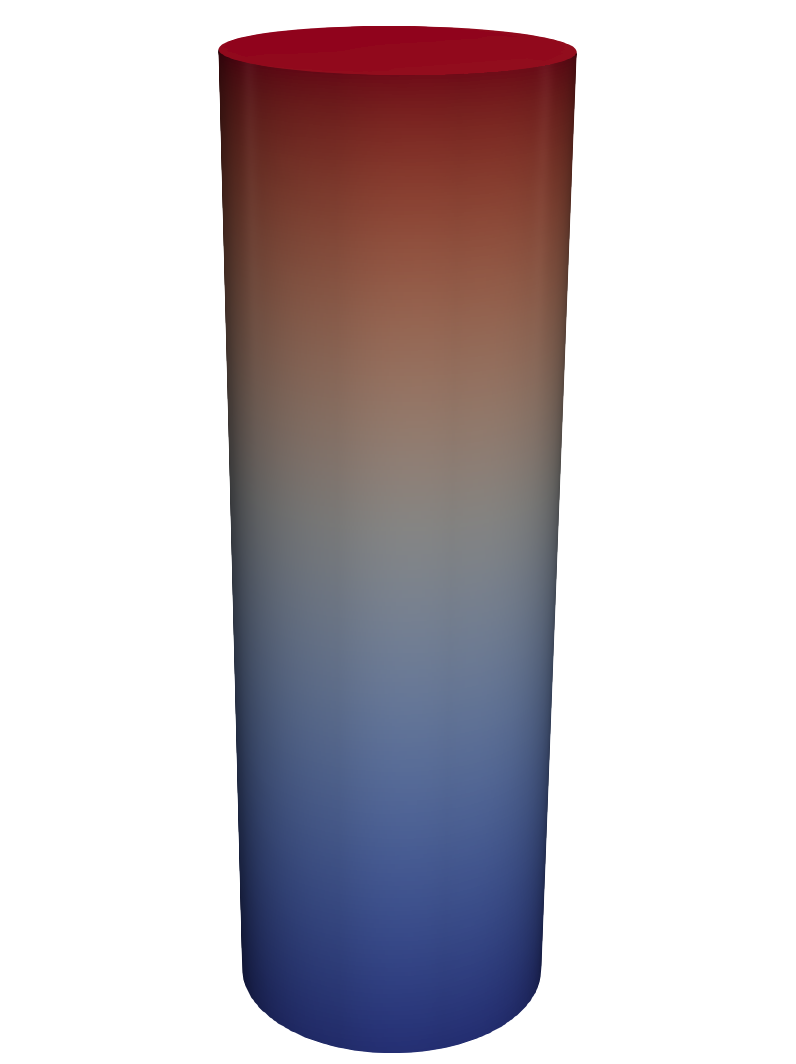}
       \includegraphics[width = 0.12\textwidth,trim={4cm 0cm 3.4cm 0cm},clip]{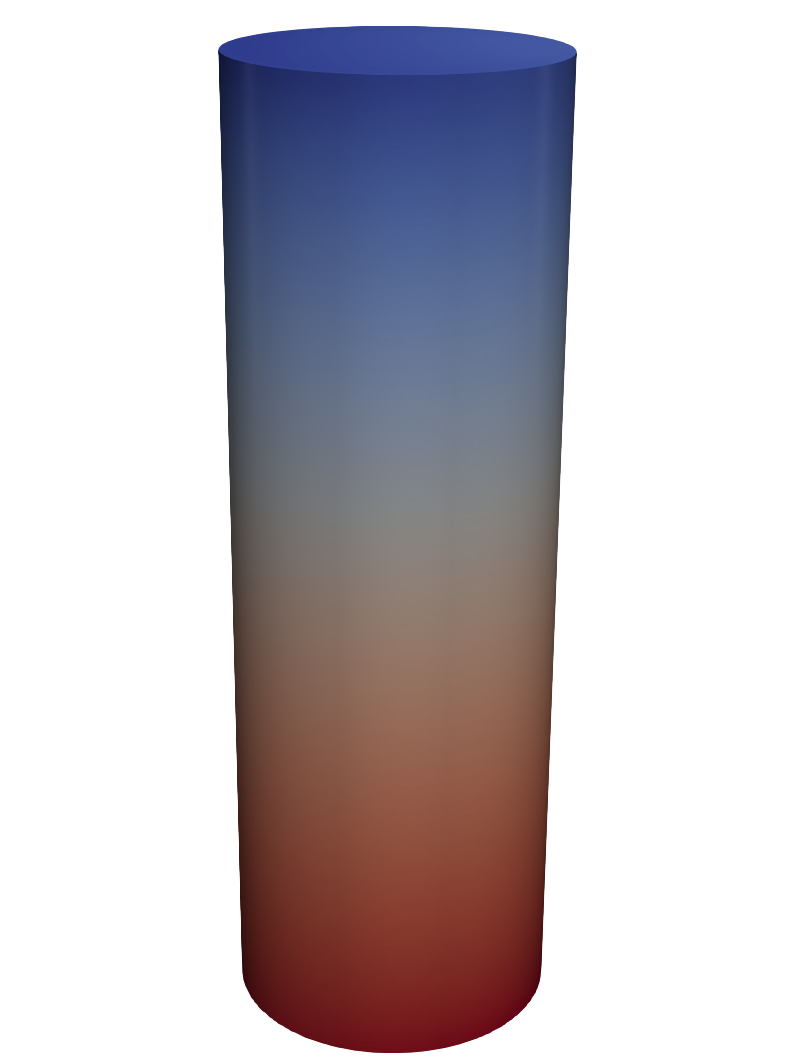}
       \includegraphics[width = 0.12\textwidth,trim={4cm 0cm 3.4cm 0cm},clip]{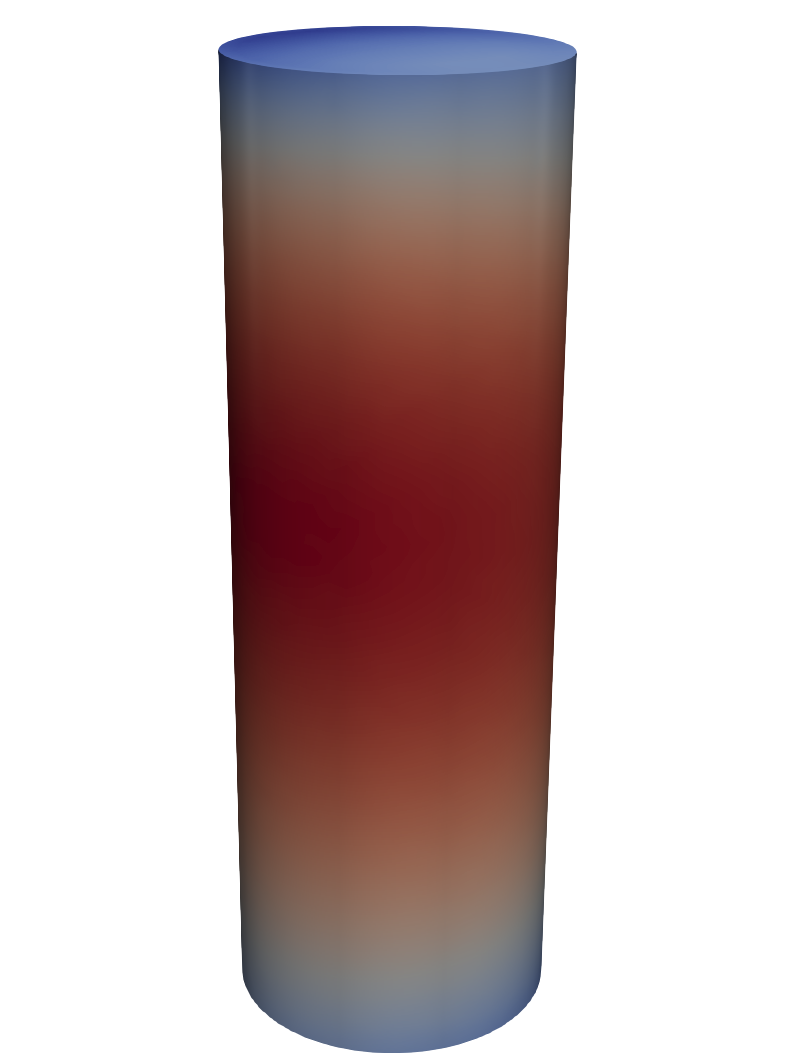}
       \includegraphics[width = 0.12\textwidth,trim={4cm 0cm 3.4cm 0cm},clip]{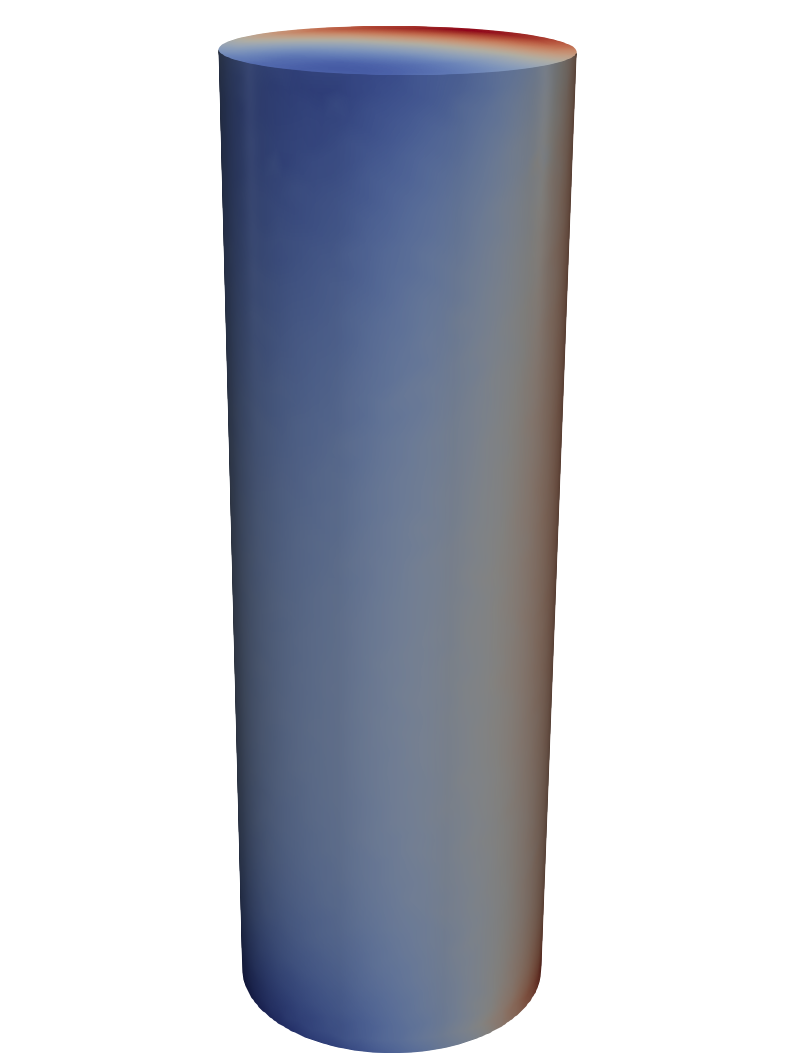}};
\end{tikzpicture}
\caption{First four velocity modes of velocity (left) and pressure (right) for the four building blocks. From top to bottom: bifurcation (B), tubes with aspect ratio length/diameter 1:1 (T1), 1:2 (T2) and 1:3 (T3).}
\label{fig:rbfunctions}
\end{figure}
\begin{table}[htbp]
\centering
\begin{tabular} {c c c c c c c c c c c c c c}
\toprule
& & & & \multicolumn{6}{c}{$N_u^i$ w.r.t. $\epsilon$} & \multicolumn{4}{c}{$N_p^i$ w.r.t. $\epsilon$} \\
\cmidrule(r{4pt}){5-10} \cmidrule(l){11-14}
 & $N_s^i$ & $N^{i,h}_u$ & $N^{i,h}_p$ & 1.6e-2 & 8e-3 & 4e-3 & 2e-3 & 1e-3 & 5e-4 & 8e-5 & 4e-5 & 2e-5 & 1e-5 \\
\midrule
B & 2640 & 76974 & 3552 & 99 & 172 & 270 & 394 & 549 & 732 & 82 & 131 & 189 & 265 \\
T1& 7920& 42708& 2162 & 69 & 134 & 233 & 379 & 582 & 848 & 54 & 86 & 133 & 198 \\
T2& 7920& 76416& 3830 & 60 & 103 & 162 & 243 & 354 & 503 & 59 & 92 & 131 & 181 \\
T3& 5280& 103728& 5211 & 23 & 45 & 78 & 131 & 211 & 324 & 29 & 45 & 66 & 95 \\
\bottomrule
\end{tabular}
\caption{Number of snapshots $N_s^\irs{i}$, velocity and pressure FE basis sizes $N^{\irs{i},h}_u$ and $N^{\irs{i},h}_p$, and velocity and pressure RB sizes $N^\irs{i}_u$ and $N^\irs{i}_p$ with respect to different POD tolerances $\epsilon_u$ and $\epsilon_p$, for 4 different building blocks (B: bifurcation, T1, T2 and T3: tubes with aspect ratios diameter/length 1:1, 1:2, 1:3, respectively).}
\label{table:redbasis}
\end{table}

The $N_s^\irs{i}$ snapshots for velocity and pressure for the $i^\text{th}$ building block are collected in matrices $\widehat{\mat{S}}_u^\irs{i} \in \reals^{N_{u}^{\irs{i},h} \times N_{s}^\irs{i}}$ and $\widehat{\mat{S}}_p^\irs{i} \in \reals^{N_{p}^{\irs{i},h} \times N_{s}^\irs{i}}$. We remark that, since we are dealing with unsteady problems, these matrices collect snapshots sampled at different timesteps for different values of the geometrical parameters. It is worth noting that each velocity snapshot, which is divergence free in the deformed configuration, does not retain such property on the reference building block. In order to consider snapshots which are divergence-free in the reference configuration, the columns of $\widehat{\mat{S}}_u^\irs{i}$ are scaled by means of the divergence-preserving Piola transformation, which is defined as follows. Given a vector field $\mathbf v$ defined in $\domain^\is{j}(\boldsymbol \param^\is{j})$ and such that $\nabla_{\mathbf x} \cdot \mathbf v = 0$, the field
\begin{equation}
    \widehat{\mathbf v}(\mathbf{\widehat{x}}) = |\mat{J}_{\Phi^{z_j}}(\widehat{\mathbf{x}};\boldsymbol \param^\is{j})| \mat{J}_{\Phi^{z_j}}^{-1}(\widehat{\mathbf{x}};\boldsymbol \param^\is{j}) \mathbf{v}(\Phi^{z_j}(\widehat{\mathbf{x}};\param^j))
    \label{eq:piola}
\end{equation}
is such that $\nabla_{\widehat{\mathbf x}} \cdot \widehat{\mathbf v} = 0$ in $\bbref^\irs{z_j}$. Matrix $\mat{J}_{\Phi^{z_j}}(\widehat{\mathbf{x}}; \param^\is{j})$ is the Jacobian of transformation $\Phi^{z_j}$ defined in Eq.~\eqref{eq:transformation}. This explicitly takes the form $\mat{J}_{\Phi^{z_j}}(\widehat{\mathbf{x}}; \param^\is{j}) =  (Q(\param^\is{j}) \mat{J}_{\varphi^{z_j}}(\widehat{\mathbf{x}}; \param^\is{j}))^{-1}$.

\begin{remark}
    The nonaffine deformation $\boldsymbol \varphi^{z_j}(\cdot;\param^\is{j})$ is defined in analytic form for tubes T1, T2 and T3; therefore, for those building blocks the Jacobian is computed exactly. The nonaffine deformation of the bifurcation B is performed by prescribing the position of the outlets in the physical configuration and by solving a linear elasticity problem such that the displacement field operates the desired rotation of such interfaces. Due to the complications of the evaluation of the Jacobian at the mesh nodes, in the bifurcation we consider $\mat{J}_{(\varphi^{z_j})^{-1}}(\mathbf{x}; \param^\is{j}) \approx I$. This simplification is also justified by the fact that for this building block we restrict ourselves to small deformations.
\end{remark}

The basis matrices $\widehat{\mat{V}}_u^\irs{i} = [\widehat{\algvec{\zeta}}_1^{\irs{i},h}|\ldots|\widehat{\algvec{\zeta}}^{\irs{i},h}_{N^\irs{i}_u}]\in \reals^{N_u^{\irs{i},h} \times N^\irs{i}_u}$ and $\widehat{\mat{V}}_p^\irs{i} = [\widehat{\algvec{\eta}}_1^{\irs{i},h}|\ldots|\widehat{\algvec{\eta}}^{\irs{i},h}_{N^\irs{i}_p}] \in \reals^{N_p^{\irs{i},h} \times N^\irs{i}_p}$ are constructed by POD as described in Section~\ref{subsec:rboffline} by considering two tolerances $\varepsilon_u$ and $\varepsilon_p$ for every building block. Matrices $\widehat{\mat{V}}_u^\irs{i}$ and $\widehat{\mat{V}}_p^\irs{i}$ are made orthonormal with respect to $\widehat{\mat{X}}_u^\irs{i}$ (matrix discretization of the H1 norm on the $i^\text{th}$ reference building block) and $\widehat{\mat{X}}_p^\irs{i}$ (matrix discretization of the L2 norm on the $i^\text{th}$ reference building block) respectively, by following the procedure presented in Remark~\ref{remark:ortho}. The first four modes of velocity and pressure for each building block are depicted in Fig.~\ref{fig:rbfunctions}. We also introduce the local basis matrices $\mat{V}_u^\is{j}$, which are obtained by applying to each column of $\widehat{\mat{V}}_u^\irs{z_j}$ the Piola transformation from the reference configuration $\bbref^{z_j}$ to the physical one $\domain^j$ (being dependent on the geometrical parameters $\pazocal M$, these must be computed during the online phase). Table~\ref{table:redbasis} reports data about the $N_\text{bb} = 4$ building blocks composing the artificial geometry used for the data generation, such as number of snapshots $N_s^\irs{i}$ and the sizes of FE and RB spaces for velocity and pressure. Although the RB sizes are considerably smaller than  the FE ones, the number of basis functions needed to achieve low POD tolerances is substantial. This indicates that the amount of information carried by the snapshots impedes the reduction of the problem. The basis size could be decreased by considering narrower sampling intervals for the geometrical parameters describing each building block. However, as in Section~\ref{subsec:numres_aorta} we test the ability of the same reduced basis to generalize to the case of a geometry which is not considered in the offline phase, here we decide to allow for significant deviations of the configurations from the original geometry during the snapshots generation. It is worth noting that, in order to decrease the already high computational burden of the offline phase, we settle for a number of configurations (165) that is possibly too limited to capture the geometrical variability we consider in the dataset (see Fig.~\ref{fig:offline} for examples of some the configurations). As we verify in Section~\ref{subsec:numres_artificial}, the errors that we obtain in the online phase are---although sufficiently low for most cardiovascular applications---considerably larger than the POD tolerance as a consequence of Remark~\ref{remark:quality}.

\subsection{Supremizers enrichment for pressure and coupling Lagrange multipliers}
\label{subsec:supremizers}
In Section~\ref{subsec:numdiscns} we recall that the Navier--Stokes equations represent an example of saddle-point equations and that this class of problems is associated with stability issues related with the discretization spaces employed for the primal and dual fields (velocity and pressure, respectively). Unfortunately, even if a stable discretization is considered during the reduced basis generation, the stability is in general not preserved in the reduced system. Furthermore, the global system obtained from the nonconforming method introduced in Section~\ref{subsec:discretization_dd} is also a saddle-point problem where the velocity and the Lagrange multipliers play the role of the primal and dual fields, respectively. Among the ways to deal with the loss of stability in the reduced system are the use of least squares Petrov--Galerkin approaches for the solution of the minimization problem associated to the nonlinear residual of the reduced equations \cite{carlberg2017galerkin,dal2019algebraic} and the supremizers enrichment \cite{ballarin2015supremizer,dalsanto2018hyper,rozza2005optimization}. Here we follow the latter approach.

The stability condition is often called inf-sup condition \cite{brezzi1974existence} and must be satisfied both at the continuous and discrete level in order to ensure the existence and uniqueness of the respective solutions. Let us consider \weakref{weak:navierstokes_dd}, which we assume to be well-posed in the continuous setting. We first address the stability with respect to the constraint imposed by the pressure (divergence free velocity). At the FE level, the inf-sup condition requires the existence of $\beta^{\is{j},h}_p \in \reals$ such that, for all $j = 1,\ldots,N_\domain$,
\begin{equation}
    \beta^{\is{j},h}_p = \inf_{\algvec{q} \neq \algvec{0}} \sup_{\algvec{v} \neq \algvec{0}} \dfrac{\algvec{q}^\text{T} \mat D^{\is{j},h} \algvec{v}}{\Vert \algvec{v} \Vert_{\pazocal V^{\is{j},h}} \Vert \algvec{q} \Vert_{\pazocal Q^{\is{j},h}}} > 0,
\end{equation}
where we used the notation $\Vert \algvec{v} \Vert_{\pazocal V^{\is{j},h}} = \algvec{v}^\text{T} \mat X_u^{\is{j},h} \algvec{v} = \Vert \mathbf{v}^h \Vert_{\pazocal V^\is{j}}$ and $\Vert \algvec{p} \Vert_{\pazocal Q^{\is{j},h}} = \algvec{p}^\text{T} \mat X_p^{\is{j},h} \algvec{p} = \Vert p^h \Vert_{\pazocal Q^\is{j}}$.
Taylor--Hood elements \cite{hood1974navier} are an example of stable choice of elements ensuring that $\beta^{\is{j},h}_p > 0$, as mentioned in Section~\ref{subsec:numdiscns}. In the RB context, we introduce $\mat X_u^{\is{j},N} = (\mat{V}_u^\is{j})^\text{T} \mat X_u^{\is{j},h} \mat{V}_u^\is{j}$, $\Vert \algvec{v}^N \Vert_{\pazocal V^{\is{j},N}} = \algvec{v}^\text{T} \mat X_u^{\is{j},N} \algvec{v}$, the pressure counterparts $\mat X_p^{\is{j},N}$ and $\Vert \algvec{p}^N \Vert_{\pazocal Q^{\is{j},N}}$, and $\mat{D}^{\is{j},N} = (\widehat{\mat V}_p^\irs{z_j})^\text{T} \mat{D}^\is{j} \mat V^\is{j}_u$. The inf-sup condition becomes
\begin{equation}
    \beta^{\is{j},N}_p = \inf_{\algvec{q} \neq \algvec{0}} \sup_{\algvec{v} \neq \algvec{0}} \dfrac{\algvec{q}^\text{T} \mat D^{\is{j},N} \algvec{v}}{\Vert \algvec{v} \Vert_{\pazocal V^{\is{j},N}} \Vert \algvec{q} \Vert_{\pazocal Q^{\is{j},N}}} > 0.
\end{equation}
The main idea of supremizers enrichment is to augment the reduced basis for the velocity with vectors (the supremizers) specifically computed from the pressure modes to ensure the positivity of the inf-sup constant. Formally, let us consider the supremizers $\algvec{s}_1^\is{j},\ldots,\algvec{s}^\is{j}_{N_p^\is{z_j}}$ which are obtained by solving for $j = 1,\ldots,N_\domain$ and for $l = 1,\ldots,N_p^\is{z_j}$ the problem $\mat{X}_u^{\is{j},h} \algvec{s}_l^\is{j} = (\mat{D}^{\is{j},h})^\text{T} \algvec{\eta}_l^{z_j}$. It can be shown \cite{ballarin2015supremizer} that substituting $\mat{V}_u^\is{j}$ with the enriched basis $\widetilde{\mat{V}}_u^{\is{j},+} = [\mat{V}^\is{z_j}_u|\algvec{s}_1^\is{j}|\ldots|\algvec{s}_{N_p^\is{z_j}}^\is{j}]$ leads to $\beta^{\is{j},N}_p > \beta^{\is{j},h}_p > 0$. A major drawback of this type of (exact) supremizers enrichment is that the supremizers for subdomain $\domain^\is{j}$ are dependent on the geometrical parameter $\param^\is{j}$, which implies that they must be computed in the online phase. For this reason, in this paper we follow an approximate approach similar to the one considered in \cite{ballarin2015supremizer} in the case of geometrical parameters. For every reference building block $\bbref^\irs{i}$ and for $l = 1,\ldots,N_p^\irs{i}$, we introduce the problems $\widehat{\mat{X}}_u^{\irs{i},h} \widehat{\algvec{s}}_l^\irs{i} = (\widehat{\mat{D}}^{\irs{i},h})^\text{T} \widehat{\algvec{\eta}}_l^\irs{i}$, where $\widehat{\mat{X}}_u^\irs{i}$ and $\widehat{\mat{D}}^{\irs{i},h}$ are velociy norm and divergence matrices assembled on the $i^\text{th}$ reference building block. Then, the enriched velocity matrix becomes $\widehat{\mat{V}}_u^{\irs{i},+} = [\widehat{\mat{V}}^\irs{i}_u|\widehat{\algvec{s}}_1^\irs{i}|\ldots|\widehat{\algvec{s}}_{N_p^\irs{i}}^\irs{i}]$. The coupling stabilization is performed by following the same procedure. In particular, in this case the problems to be solved for each reference building block $\bbref^i$ read $\widehat{\mat{X}}_u^{i,h} \widehat{\algvec{z}}_l^i = (\widehat{\mat{B}}^{\is{i}\iin{m},h\delta})^\text{T} \algvec{e}_l$, for $l = 1,\ldots,N_\lambda$; $\widehat{\mat{B}}^{\is{i}\iin{m},h\delta}$ is the matrix assembled on $\bbref^i$ discretizing the coupling with the $m^\text{th}$ interface (specifically, if $\bbref^i$ is a tube $m = 1,2$, whereas if it is a bifurcation $m > 2$). In the following, we simply denote by $\widehat{V}_u^\irs{i}$ the enriched basis matrix for the velocity in $\bbref^\irs{i}$ obtained by arranging columnwise $\widehat{\algvec{\zeta}}^{\irs{i},h}_l$ and the supremizers for the pressure and coupling stabilizations $\widehat{\algvec{s}}_l^\irs{i}$ and $\widehat{\algvec{z}}_l^\irs{i}$; the enriched basis $\widehat{V}_u^\irs{i}$ is made orthonormal with respect to $\widehat{\mat X}_u^i$ with the Gram--Schmidt algorithm. The numerical results in Section~\ref{sec:numericalresults} are obtained by following this stabilization strategy.
\subsection{Assembly and solution of the global reduced system}
\label{subsec:online_rbdd}
Let us define the global basis matrix
\begin{equation}
\pazocal{W} := \text{diag}\left (
\begin{bmatrix}
\mat{V}_u^\irs{z_j} &  \\
                 & \widehat{\mat{V}}_p^\irs{z_j}
\end{bmatrix}
\right )_{j = 1,\ldots,N_\domain},
\end{equation}
the matrices
\begin{equation}
\pazocal M^N := \pazocal W^\text{T} \pazocal M^h \pazocal W,\quad
\pazocal A^N(\algvec{W}^N) :=\pazocal W^\text{T} \pazocal A^h(\pazocal W \algvec{W}^N) \pazocal W,\quad
\pazocal B^{N\delta} := \pazocal B^{h\delta} \pazocal W,
\end{equation}
the vector of reduced degrees of freedom for all the subdomains $\algvec{W}^{N} = [\algvec{w}^{\is{1},N},\ldots,\algvec{w}^{\is{N_\domain},N}]$, and the vectors encoding the data $\algvec{F}^N := \pazocal W^\text{T} \algvec{F}^h$ and $\algvec{G}^{N\delta} := \algvec{G}^{h\delta}$. Then, the reduced residual at timestep $t_{k+1}$ is obtained from Eq.~\eqref{eq:residualbdf_dd} and reads
\begin{equation}
    \algvec{R}^N(\algvec{Y}^{N\delta}_{k+1}) := \pazocal{H}^N \algvec{Y}^{N\delta}_{k+1} - \sum_{j = 1}^{\sigma} \alpha_j \pazocal{H}^N \algvec{Y}^{N\delta}_{k-j+1} - \Delta t \beta \mathring{\algvec{F}}^{N\delta}(t_{k+1},\algvec{Y}^{N\delta}_{k+1}) = \algvec{0},
\end{equation}
where
\begin{equation}
\pazocal H^N :=
\begin{bmatrix}
    \pazocal M^N & \\
     & \\
\end{bmatrix},
\quad
\mathring{\algvec{F}}^{N\delta}(t,\algvec{Y}^N)
:=
\begin{bmatrix}
    \algvec{F}^N(t) \\
    \algvec{G}^{N\delta}(t)
\end{bmatrix}
-
\begin{bmatrix}
\pazocal A^N (\algvec{W}^N) & (\pazocal B^{N\delta})^\text{T} \\
\pazocal B^{N\delta} & \\
\end{bmatrix}
\begin{bmatrix}
    \algvec{W}^N \\
    \algvec{\Lambda}^\delta
\end{bmatrix},
\end{equation}
and $\algvec{Y}^{N\delta} := [\algvec{W}^N,\algvec{\Lambda}^\delta]$.

As discussed in Section~\ref{subsec:preconditioner}, finding the root of nonlinear equations using the Newton--Raphson algorithm entails the solution of a nonlinear system in the tangent matrix of the corresponding residual. Formally, solving $\algvec{R}^N(\algvec{Y}) = \algvec{0}$ given an initial guess $\algvec{Y}^{(0)}$ leads to the iterative algorithm
\begin{equation}
    \algvec{Y}^{(l+1)} = \algvec{Y}^{(l)} - \left(\pazocal J_{\algvec{R}^N}\left (\algvec{Y}^{(l),N}\right )\right)^{-1}\algvec{R}^N\left (\algvec{Y}^{(l)}\right ),
    \label{eq:residual_block_red}
\end{equation}
which is equivalent to Eq.\eqref{eq:residual_block} in the reduced context.

The efficiency of the reduction relies on the fast assembly of the tangent matrix $\pazocal J_{\algvec{R}^N}$ and residual $\algvec{R}^N$. Regarding the former, we consider the following approximation
\begin{equation}
\widetilde{\pazocal J}_{\algvec{R}}^N =
\begin{bmatrix}
\pazocal M^N + \Delta t \beta \pazocal A^N_\text{lin} & \Delta t \beta (\pazocal B^{N\delta})^\text{T} \\
\Delta t \beta \pazocal B^{N\delta} &
\end{bmatrix},
\label{eq:reduced_jac}
\end{equation}
where
\begin{equation}
    \pazocal A^N_\text{lin} := \pazocal W^\text{T} A_\text{lin} \pazocal W,\quad
    \pazocal A_\text{lin} := \text{diag}\left (
    \begin{bmatrix}
        \mat{K}^{\is{j},h} & (\mat{D}^{\is{j},h})^\text{T} \\
        \mat{D}^{\is{j},h} & \\
    \end{bmatrix}
    \right )_{j = 1,\ldots,N_\domain}
\end{equation}
is the matrix obtained by neglecting the convective terms in the Navier--Stokes equations. We remark that the reduced tangent matrix features a saddle-point structure as its full order counterpart in Eq.~\eqref{eq:jacobian}. Therefore, system $\widetilde{\pazocal J}_{\algvec{R}^N} \algvec{X}^N = \algvec{B}^N$ can be solved directly by applying the reduced version of Eq.~\eqref{eq:solution_jac}. The advantages of this approach are: (i) the tangent matrix $\widetilde{\pazocal J}_{\algvec{R}^N}$ is never entirely allocated, because every stage for applying Eq.~\eqref{eq:solution_jac} involves operations that are local to either subdomains or interfaces (we recall that inverting $\pazocal{A}_\text{lin}$ amounts to inverting each of its diagonal blocks), (ii) as a result of approximation \eqref{eq:reduced_jac}, the tangent matrix is the same for every solution of the linear system; hence the reduced Schur complement is assembled only once and it can be factorized---along with the other local matrices to be inverted---at the start of the simulation. The linearized version of the tangent matrix \eqref{eq:reduced_jac} is nonconsistent, which implies that the Newton--Raphson algorithm is not expected to convergence quadratically. However, the reduced complexity of the assembly results in a considerable performance gain overall.

\begin{figure}
    \tikzsetnextfilename{rbcoefs}
    \centering
    \setlength
    \figureheight{0.2\textwidth}
    \setlength
    \figurewidth{0.9\textwidth}
    \input{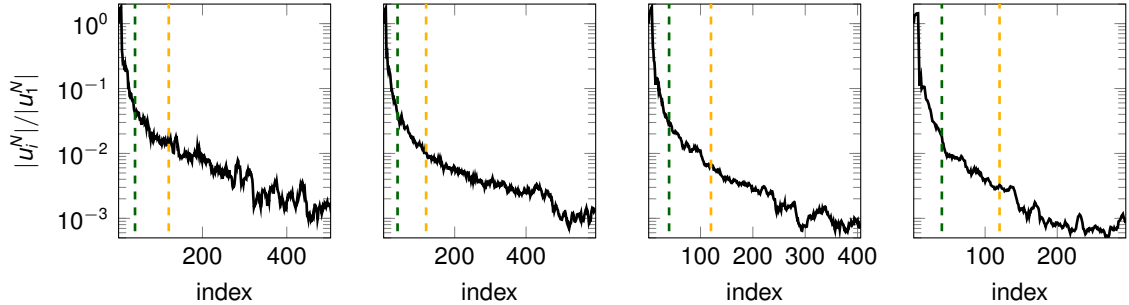}
    \caption{Average over time of the RB velocity solutions---normalized with respect to the first coefficient---in the four building blocks (from left to right: bifurcation B and tubes T1, T2 and T3) in the RB simulation considered in Section~\ref{subsec:numres_artificial} with $\epsilon_u = \text{4e--3}$ and $\epsilon_p = \text{8e--5}$. The left and right dashed lines in every plot correspond to the indices 40 and 120, which are two of the choices considered in Section~\ref{subsec:numres_artificial} for the truncation of the computation of the convective term.}
    \label{fig:rbcoefs}
\end{figure}

The problematic part of the computation of $\algvec{R}^N(\algvec{Y})$ is evidently the nonlinear term, as the matrices encoding the linear ones are computed only once and the corresponding contributions are found at each timestep by inexpensive matrix-vector multiplications. After trivial but repetitious steps, we find that the the blocks of the nonlinear part of the reduced residual read, for all $j = 1,\ldots,N_\domain$,
\begin{equation}
\algvec{c}^{\is{j},N} = \left (\mat{V}_u^\is{z_j}\right )^\text{T} \mat{C}\left ( \mat{V}_u^\is{z_j}\algvec{u}^{\is{j},N} \right ) \mat{V}_u^\is{z_j}\algvec{u}^{\is{j},N}.
\label{eq:nlin_red}
\end{equation}
One way to compute the nonlinear term for every block is then to assemble the full order nonlinear term $\mat{C} ( \mat{V}_u^\is{z_j}\algvec{u}^{\is{j},N} ) \mat{V}_u^\is{z_j}\algvec{u}^{\is{j},N}$ and to project it onto the reduced space. Another strategy is based on the decomposition
\begin{equation}
\algvec{c}^{\is{j},N}_i = \sum_{l = 1}^{N_u^{z_j}}\sum_{m = 1}^{N_u^{z_j}}u_l^{\is{j},N}u_m^{\is{j},N}\int_{\domain_j} \left [\left (\bm \zeta_m^{z_j,h} \cdot \nabla \right )\bm \zeta_l^{z_j,h} \right ]\bm \zeta_i^{z_j,h}.
\label{eq:nonlin}
\end{equation}
The integrals in Eq.\eqref{eq:nonlin} are independent of the reduced solution and can be computed as a setup step in the first stages of the simulation. However, the amount of computation increases quadratically with the size of the reduced basis and could therefore nullify the performance gain.
Using the fact that the velocity modes in $\mat{V}_u^{\is{z_j}}$ are sorted in order of significance as a consequence of Proposition~\ref{prop:optimality}, it is legitimate to consider the following approximation
\begin{equation}
\algvec{c}^{\is{j},N}_i \approx \sum_{l = 1}^{N_c^{z_j}}\sum_{m = 1}^{N_c^{z_j}}u_l^{\is{j},N}u_m^{\is{j},N}\int_{\domain_j} \left [\left (\bm{\zeta}_m^{z_j,h} \cdot \nabla \right )\bm{\zeta}_l^{z_j,h} \right ]\bm{\zeta}_i^{z_j,h},
\label{eq:nonlin_red_sym}
\end{equation}
where $0 < N_c^{j} \leq N_u^{z_j}$. In other words, the observation that the magnitude of the reduced coefficients $u_i^{\is{j},N}$ quickly decreases as $i$ increases---as we show in Fig.~\ref{fig:rbcoefs} for the case of one of the simulations presented in Section~\ref{subsec:numres_artificial}---allows to truncate the two sums in Eq.~\eqref{eq:nonlin} to the first $N_c^{j}$ terms. In the numerical results in Section~\ref{sec:numericalresults} we investigate the effects of considering both Eq.~\eqref{eq:nlin_red} and Eq.~\eqref{eq:nonlin_red_sym} for the computation of the convective part of the residual.

\section{Numerical results}
\label{sec:numericalresults}

In this section, we assess the performance of our numerical method in two applications. In Section~\ref{subsec:numres_artificial}, we consider the same modular artificial geometry employed in the offline phase presented in Section~\ref{sec:nsmodularrb} and we compare the results obtained by solving the flow problem using the RB and the FE methods. 
These are obtained with a code based on LifeV, a C++ FE library with support to high-performance computing \cite{bertagna2017lifev}. In Section~\ref{subsec:numres_aorta}, we consider a simple but more physiological geometry of an aorta and the two iliac arteries. In this case, we compare the results obtained with the RB method against the ones obtained on a reference geometry (i.e. not partitioned into approximated subdomains) with SimVascular\footnote{http://simvascular.github.io/} \cite{updegrove2017simvascular}, an open-source software for patient-specific modeling and blood flow simulations.

For all the simultations, we fix $\rho_\text{f} = 1.06$ $\text{gr}\,\text{cm}^{-3}$, $\mu_\text{f} = 0.04$ $\text{gr}\,\text{cm}^{-1}\,\text{s}^{-1}$, and we consider the same choice for the discrete Lagrange multipliers space as in the snapshot generation phase (i.e. $N_\lambda = 63$ for each interface).

\subsection{Online phase on the artificial geometry}
\label{subsec:numres_artificial}

\begin{figure}
\centering
\tikzsetnextfilename{qualitative}
\begin{tikzpicture}
    \node[inner sep=0pt] (vt15) at (0,0)
    {\includegraphics[width = 0.23\textwidth]{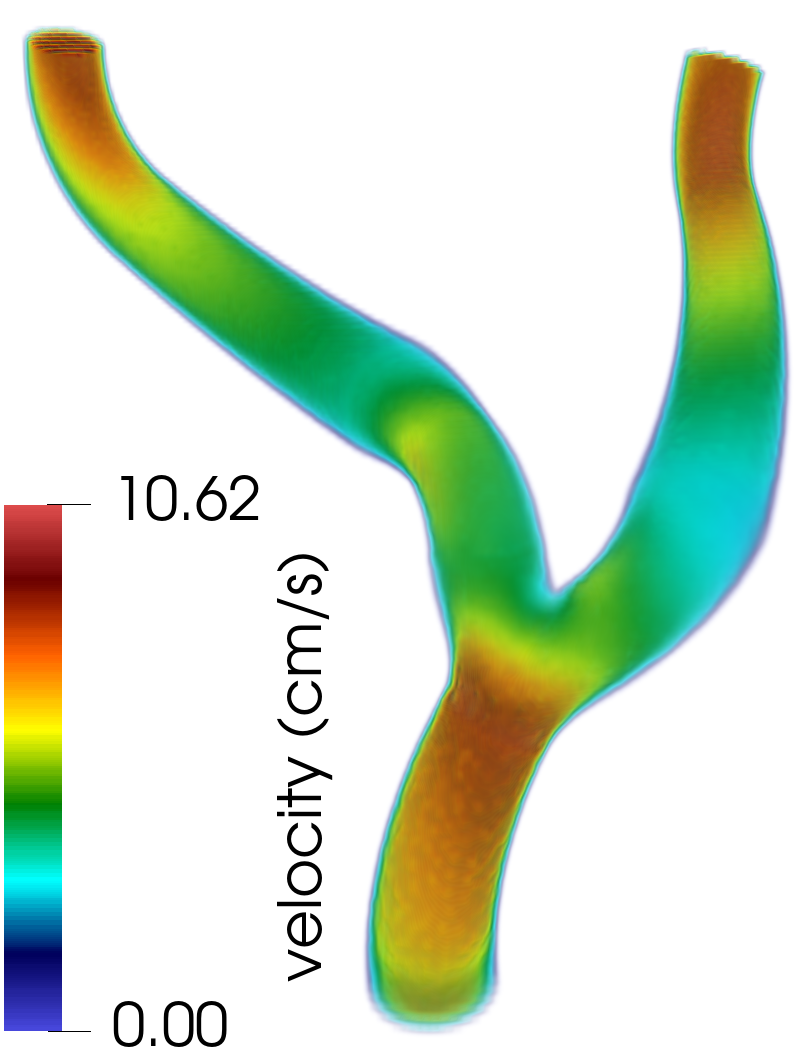}
     \includegraphics[width = 0.23\textwidth]{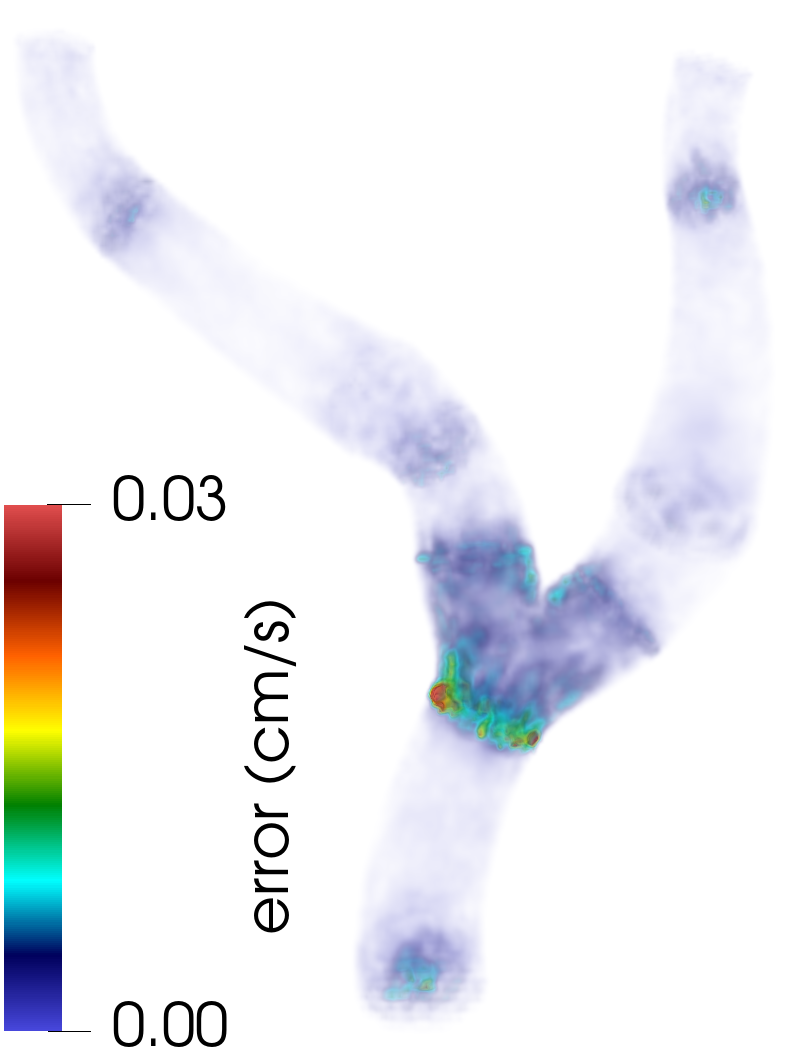}};
    \node[inner sep=0pt,right=0.6cm of vt15] (vt25)
    {\includegraphics[width = 0.23\textwidth]{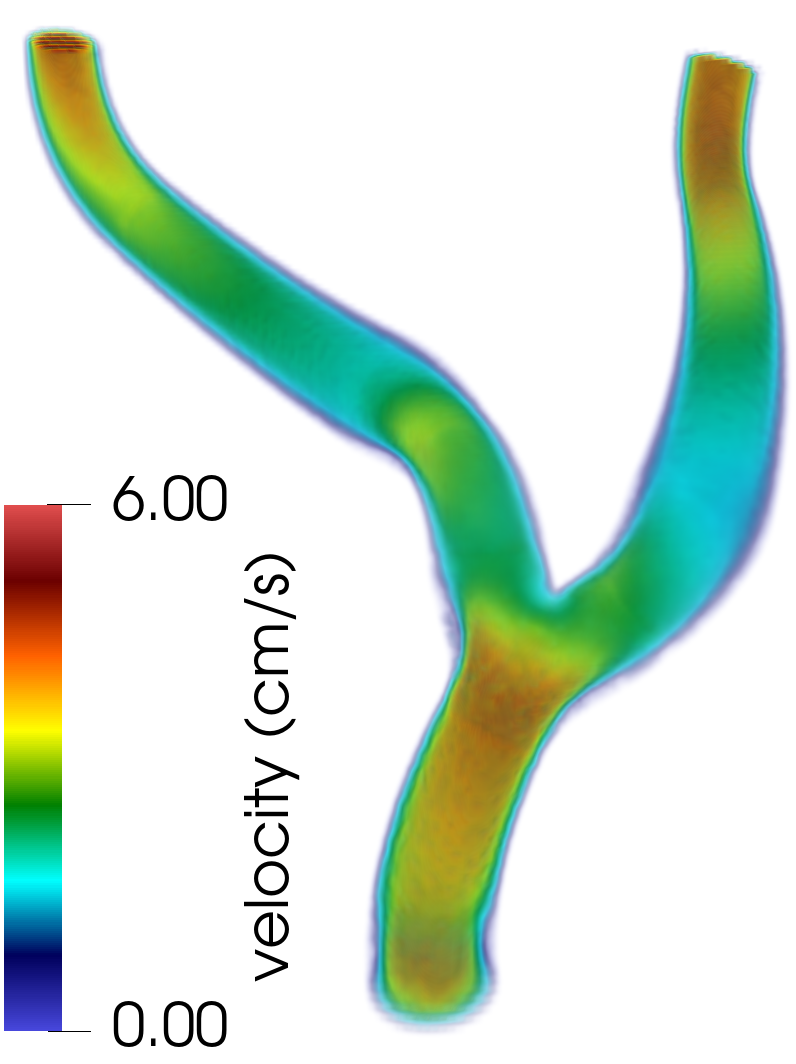}
     \includegraphics[width = 0.23\textwidth]{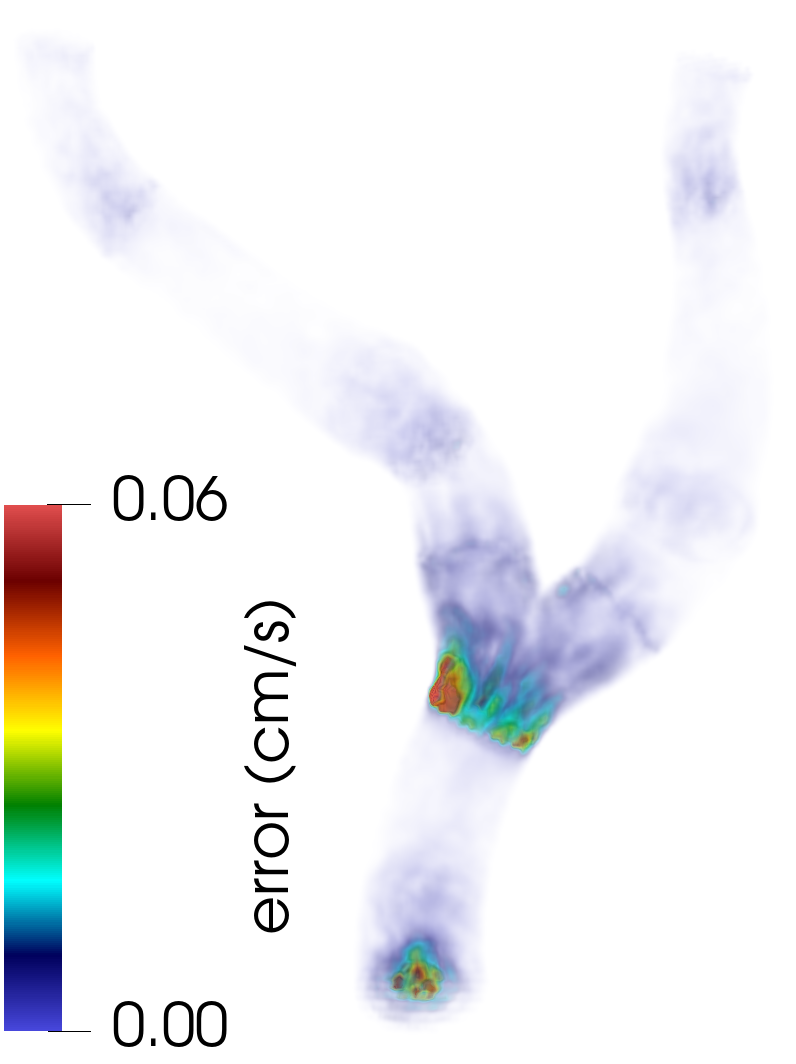}};
    \node[above=0.2cm of vt15] {$t = 0.15$ s};
    \node[above=0.2cm of vt25] {$t = 0.25$ s};
    \node[inner sep=0pt,below=0.6cm of vt15] (pt15)
    {\includegraphics[width = 0.23\textwidth]{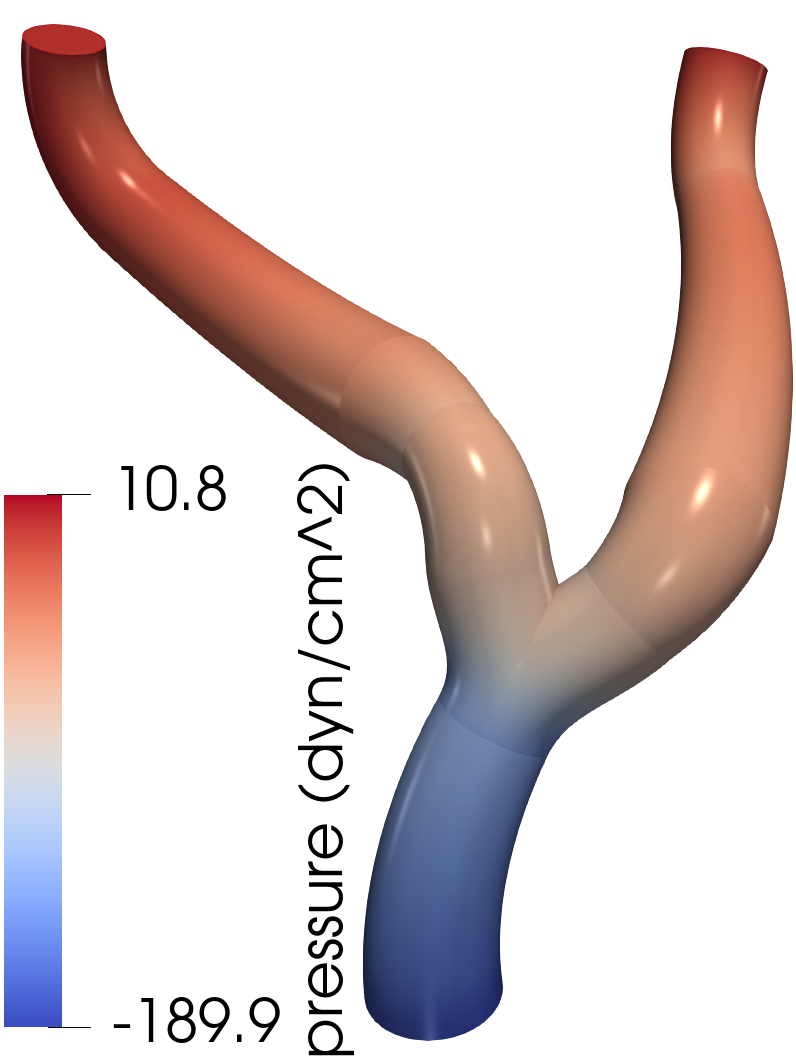}
     \includegraphics[width = 0.23\textwidth]{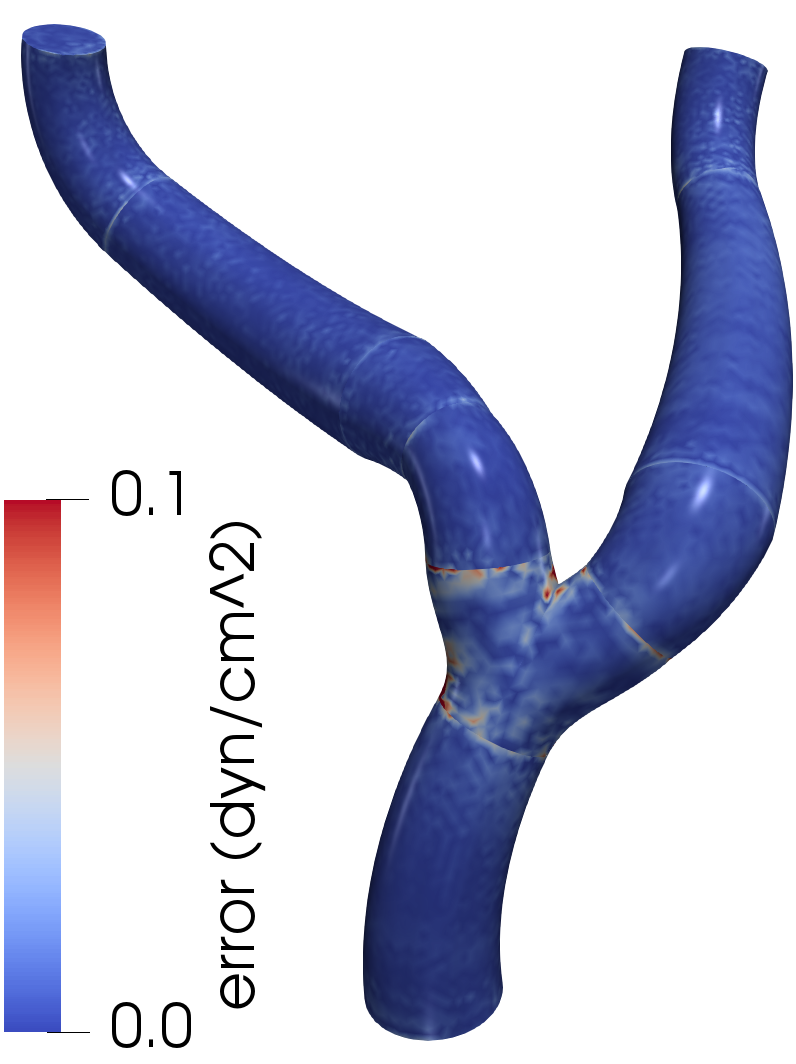}};
    \node[inner sep=0pt,right=0.6cm of pt15] (pt25)
    {\includegraphics[width = 0.23\textwidth]{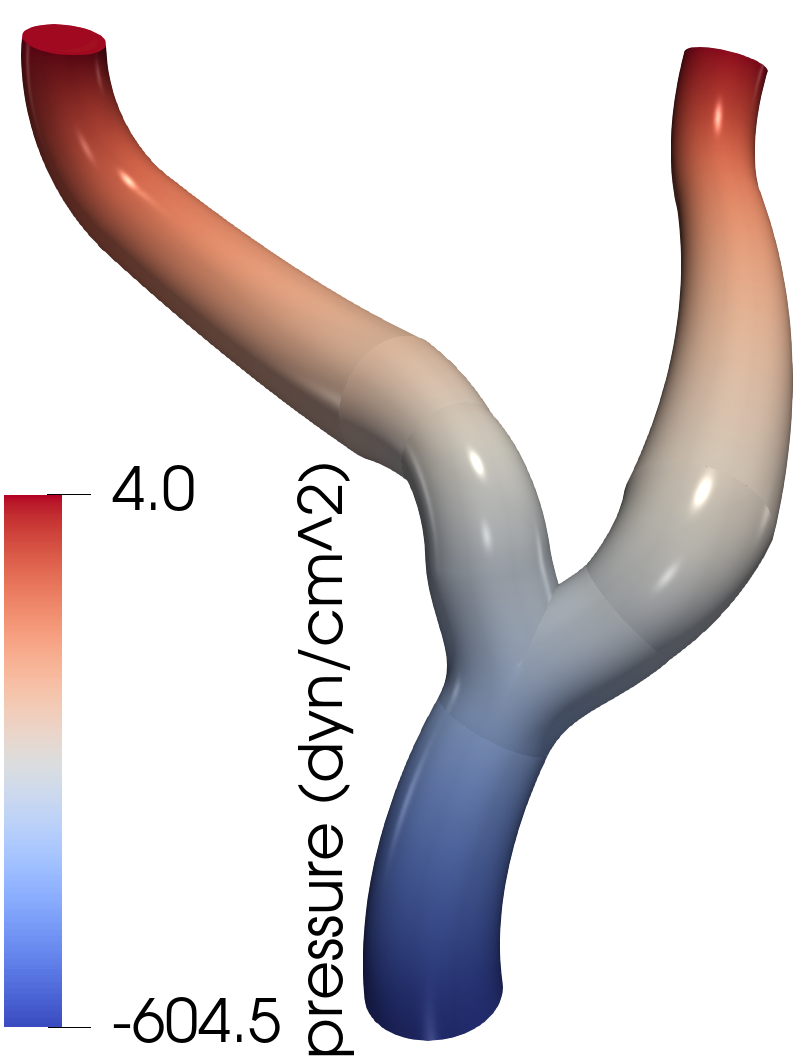}
     \includegraphics[width = 0.23\textwidth]{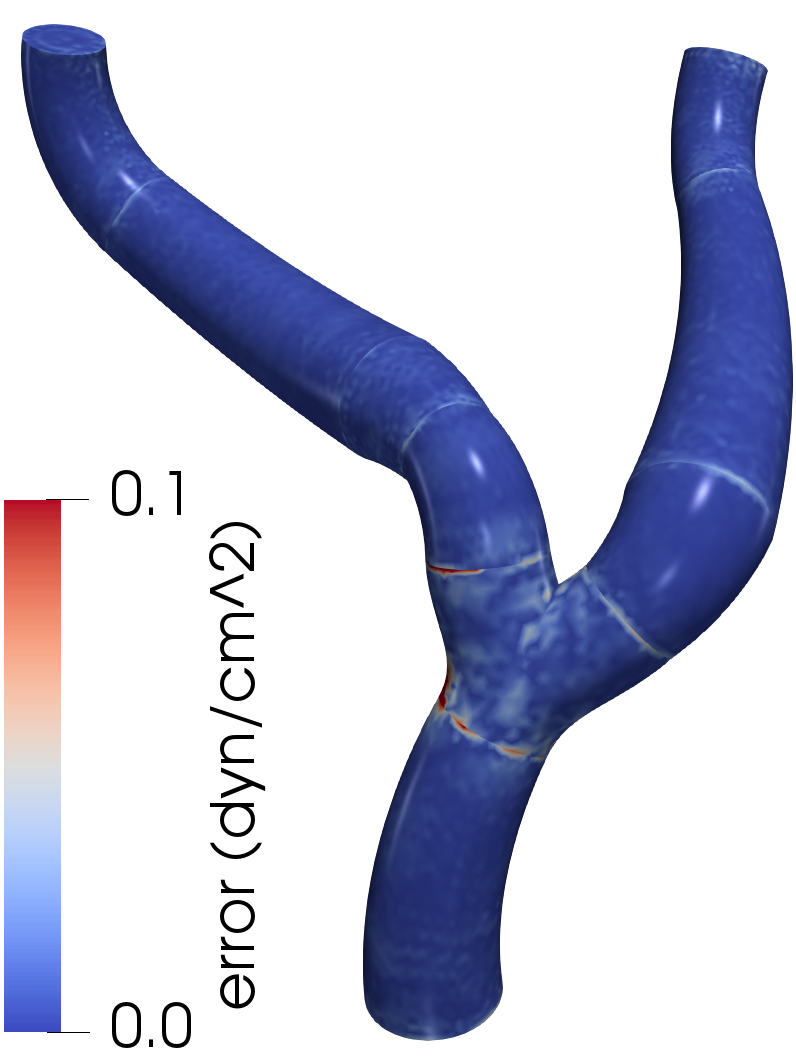}};
     \node[inner sep=0pt,below=0.6cm of pt15] (wss15)
     {\includegraphics[width = 0.23\textwidth]{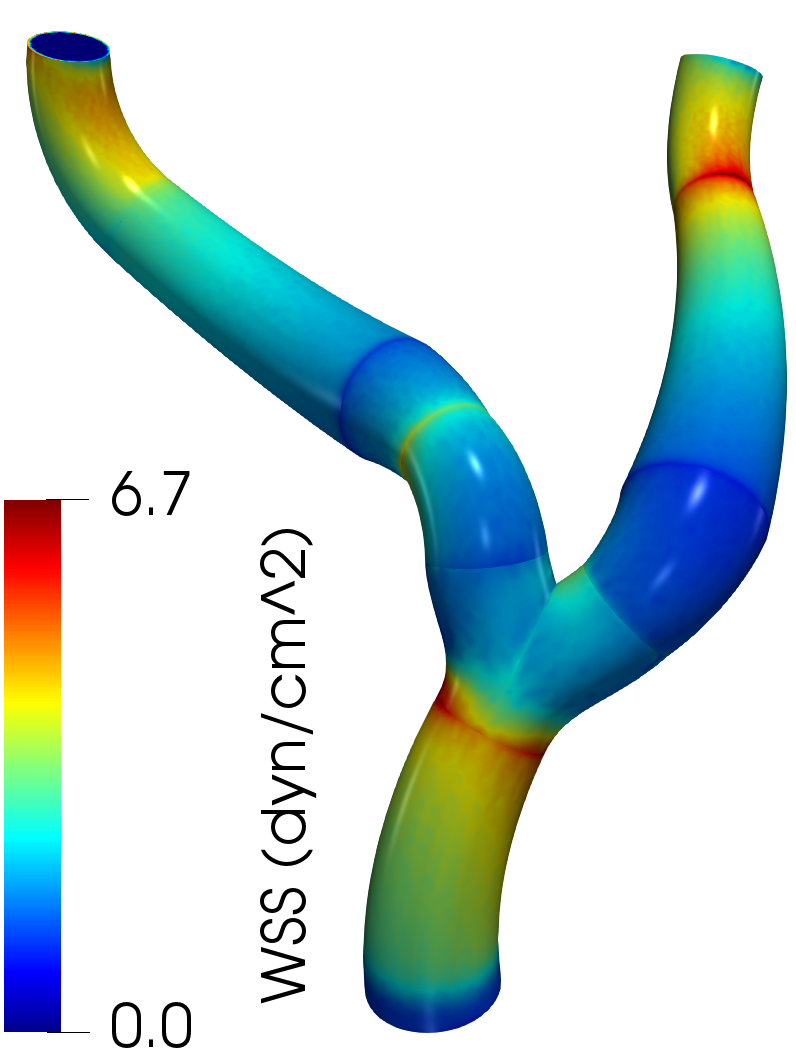}
      \includegraphics[width = 0.23\textwidth]{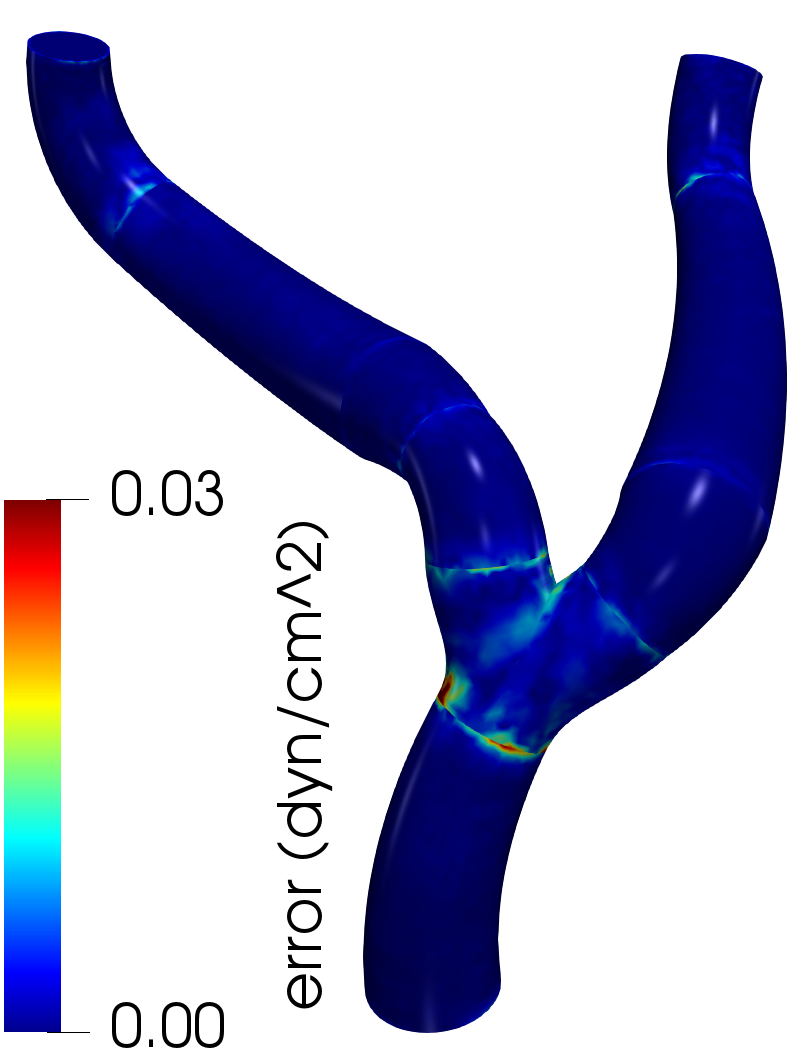}};
     \node[inner sep=0pt,right=0.6cm of wss15] (wss25)
     {\includegraphics[width = 0.23\textwidth]{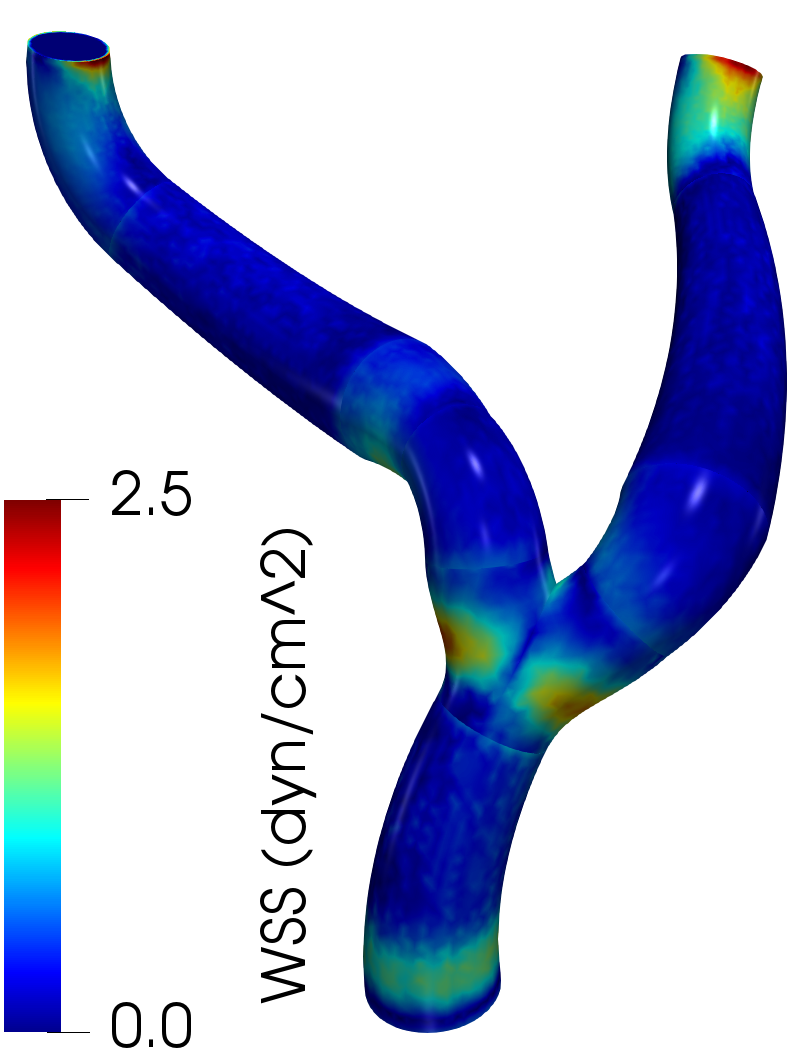}
      \includegraphics[width = 0.23\textwidth]{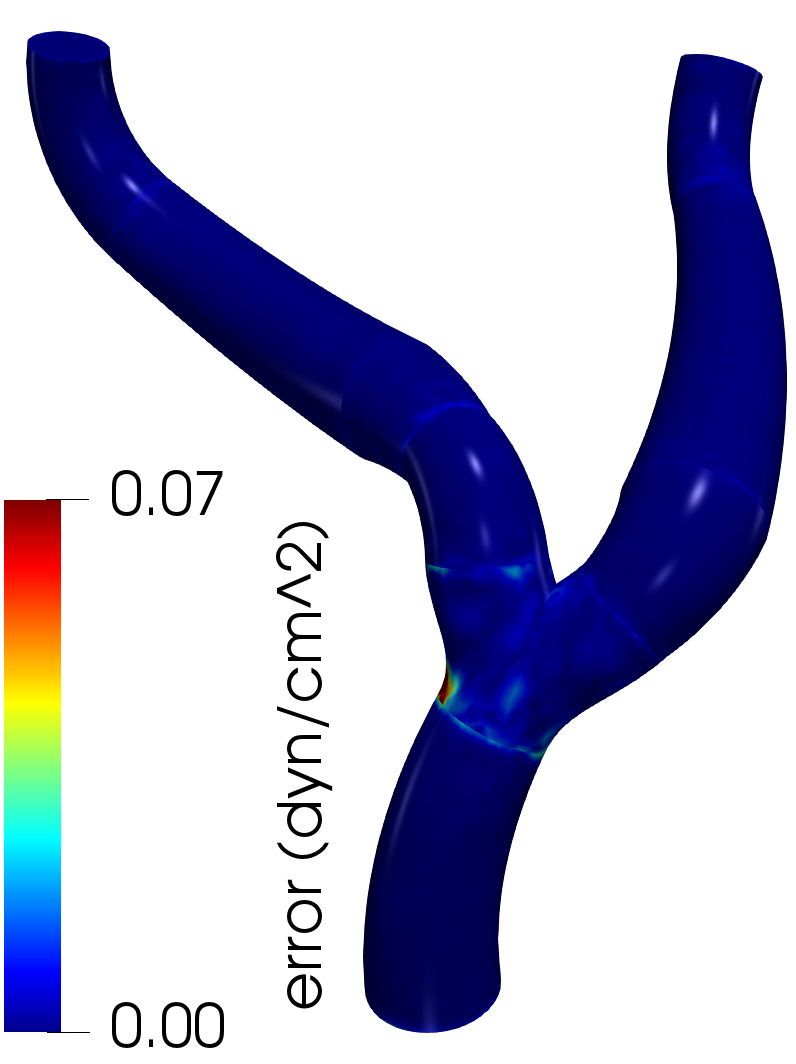}};
\end{tikzpicture}
\caption{The left and right columns---each composed of two sub-columns of plots---refer to time $t = 0.15$ s and $t = 0.25$ s, respectively. First row: velocity magnitude volume plot of the RB solution (sub-column left) and magnitude of the point-wise velocity error w.r.t. the FE solution (sub-column right). Second row: pressure plot of the RB solution (sub-column left) and absolute value of the point-wise pressure error w.r.t. the FE solution (sub-column right). Third row: magnitude of the WSS of the RB solution (sub-column left) and magnitude of the point-wise WSS error w.r.t. the FE solution (sub-column right). The RB solution corresponds to the choice $\epsilon_u = 1\text{e--3}$ and $\epsilon_p = 1\text{e--5}$.}
\label{fig:qualitative}
\end{figure}

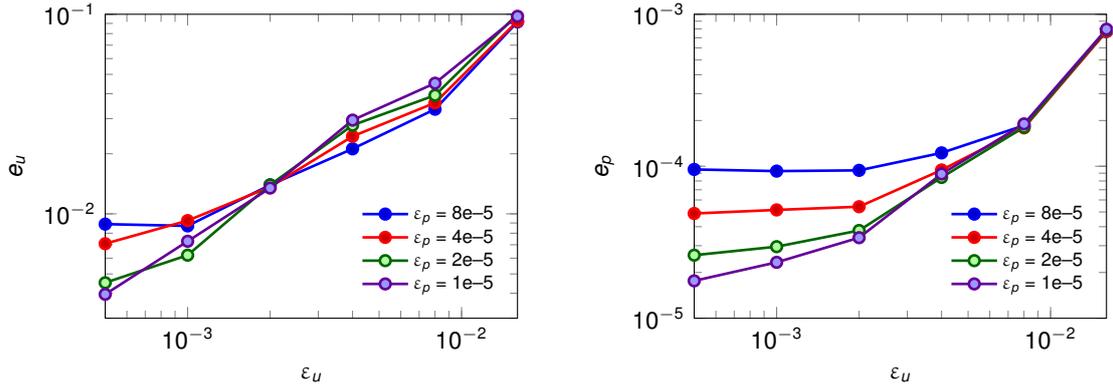
\begin{figure}
    \tikzsetnextfilename{errors_aorta_portion}
    \centering
    \setlength
    \figureheight{\textwidth}
    \setlength
    \figurewidth{\textwidth}
%
%
\definecolor{mypurple}{rgb}{0.6000,0.6000,1.0000}%
\definecolor{mycolor1}{rgb}{0.6000,0.6000,1.0000}%
\definecolor{mycolor2}{rgb}{0.6000,0.6000,1.0000}%
\definecolor{mycolor3}{rgb}{0.6000,0.6000,1.0000}%
\definecolor{mycolor4}{rgb}{0.6000,0.6000,1.0000}%

\begin{tikzpicture}

\begin{axis}[%
width=0.35\figurewidth,
height=0.26\figureheight,
at={(0\figurewidth,0\figureheight)},
scale only axis,
xmode=log,
xmin=0.0005000000,
xmax=0.0160000000,
xminorticks=true,
xlabel={$\varepsilon_u$},
ymode=log,
ymin=0.0030000000,
ymax=0.100000000,
yminorticks=true,
ylabel={$e_u$},
axis background/.style={fill=white},
legend style={at={(0.6,0.4)}, anchor=north west, legend cell align=left, align=left}
]
\addplot [blue,every mark/.append style={fill=blue!80!black},mark=*,line width = 1.0pt]
  table[row sep=crcr]
{%
0.0005000000	0.0088823724\\
0.0010000000	0.0087031637\\
0.0020000000	0.0139507596\\
0.0040000000	0.0211798275\\
0.0080000000	0.0334381219\\
0.0160000000	0.0915903056\\
};
\addlegendentry{$\varepsilon_p = 8\text{e--5}$}

\addplot [red,every mark/.append style={fill=red!80!black,solid},mark=*,line width = 1.0pt]
  table[row sep=crcr]
{%
0.0005000000	0.0070860560\\
0.0010000000	0.0092498922\\
0.0020000000	0.0136472999\\
0.0040000000	0.0244601522\\
0.0080000000	0.0360127891\\
0.0160000000	0.0924143989\\
};
\addlegendentry{$\varepsilon_p = 4\text{e--5}$}

\addplot [green!40!black,every mark/.append style={fill=mygreen,solid},mark=*,line width = 1.0pt]
  table[row sep=crcr]{%
0.0005000000	0.0045169412\\
0.0010000000	0.0062082706\\
0.0020000000	0.0139628585\\
0.0040000000	0.0278044544\\
0.0080000000	0.0392798412\\
0.0160000000	0.1001710699\\
};
\addlegendentry{$\varepsilon_p = 2\text{e--5}$}

\addplot [red!40!blue,every mark/.append style={fill=mypurple,solid},mark=*,line width = 1.0pt]
  table[row sep=crcr]{%
0.0005000000	0.0039538060\\
0.0010000000	0.0072898132\\
0.0020000000	0.0134663643\\
0.0040000000	0.0294826948\\
0.0080000000	0.0451377811\\
0.0160000000	0.0976545741\\
};
\addlegendentry{$\varepsilon_p = 1\text{e--5}$}

\end{axis}

\begin{axis}[%
width=0.35\figurewidth,
height=0.26\figureheight,
at={(0.5\figurewidth,0\figureheight)},
scale only axis,
xmode=log,
xmin=0.0005000000,
xmax=0.0160000000,
xminorticks=true,
xlabel={$\varepsilon_u$},
ymode=log,
ymin=0.000010000,
ymax=0.001000000,
yminorticks=true,
ylabel={$e_p$},
axis background/.style={fill=white},
legend style={at={(0.6,0.4)}, anchor=north west, legend cell align=left, align=left}
]
\addplot [blue,every mark/.append style={fill=blue!80!black},mark=*,line width = 1.0pt]
  table[row sep=crcr]{%
0.0005000000	0.0000952605\\
0.0010000000	0.0000928815\\
0.0020000000	0.0000939531\\
0.0040000000	0.0001224450\\
0.0080000000	0.0001855353\\
0.0160000000	0.0007721430\\
};
\addlegendentry{$\varepsilon_p = 8\text{e--5}$}

\addplot [red,every mark/.append style={fill=red!80!black,solid},mark=*,line width = 1.0pt]
  table[row sep=crcr]
{%
0.0005000000	0.0000488688\\
0.0010000000	0.0000515924\\
0.0020000000	0.0000541541\\
0.0040000000	0.0000944625\\
0.0080000000	0.0001789382\\
0.0160000000	0.0007717241\\
};
\addlegendentry{$\varepsilon_p = 4\text{e--5}$}

\addplot [green!40!black,every mark/.append style={fill=mygreen,solid},mark=*,line width = 1.0pt]
  table[row sep=crcr]{%
0.0005000000	0.0000259856\\
0.0010000000	0.0000295612\\
0.0020000000	0.0000377796\\
0.0040000000	0.0000843681\\
0.0080000000	0.0001807089\\
0.0160000000	0.0007876148\\
};
\addlegendentry{$\varepsilon_p = 2\text{e--5}$}

\addplot [red!40!blue,every mark/.append style={fill=mypurple,solid},mark=*,line width = 1.0pt]
  table[row sep=crcr]{%
0.0005000000	0.0000176157\\
0.0010000000	0.0000233247\\
0.0020000000	0.0000338799\\
0.0040000000	0.0000889907\\
0.0080000000	0.0001899401\\
0.0160000000	0.0007973689\\
};
\addlegendentry{$\varepsilon_p = 1\text{e--5}$}

\end{axis}

\end{tikzpicture}%
    \caption{Error on velocity $e_u$ (left) and error on pressure $e_p$ (right), computed as in Eq.~\eqref{eq:errors}, in function of the POD tolerances for velocity and pressure $\varepsilon_u$ and $\varepsilon_p$.}
    \label{fig:erroraortaportion}
\end{figure}

We evaluate the performance of our reduced order model on the artifical problem employed for the generation of the reduced basis. We recall that the \textit{Artificial geometry} in Fig.~\ref{fig:offline}, which is not included in the set of 165 configurations used to produce the snapshots, is a legitimate candidate to test the accuracy on geometries not ``seen'' in the offline phase. The solution by the RB method is compared to the global solution obtained by considering FE method solutions in each subdomain (with the same meshes used in the RB case) coupled with the discretization strategy presented in Section~\ref{subsec:discretization_dd}. We consider the same choice for the discrete Lagrange multipliers space as in the snapshot generation phase. The reasons for considering such comparison are the following: (i) being that the geometry is exactly the same, it is possible to easily compute H1 and L2 errors for velocity and pressure in order to verify the convergence of the RB approximation with respect to the FE one, and (ii) it is possible to fairly discuss the speedup achieved by the RB method, as the RB and FE solutions share the same computational mesh. As for the generation of the reduced basis, we consider $t_0 = 0$ s, $T = 0.3$ s, and a second order BDF scheme with $\Delta t = 2.5\times10^{-3}$ s.

\begin{table}
\centering
\begin{tabular} {c c c c c c c}
\toprule
$\varepsilon_p$\textbackslash$\varepsilon_u$ & $1.6\text{e--2}$ & $8\text{e--3}$ & $4\text{e--3}$ & $2\text{e--3}$ & $1\text{e--3}$ & $5\text{e--4}$ \\
\midrule
$8\text{e--5}$  & 33(50) & 30(48) & 26(43) & 22(39) & 17(32) & 14(31) \\
$4\text{e--5}$  & 27(46) & 25(44) & 23(38) & 21(38) & 17(35) & 13(30) \\
$2\text{e--5}$  & 28(45) & 26(43) & 23(41) & 28(38) & 15(33) & 12(28) \\
$1\text{e--5}$  & 26(45) & 25(44) & 22(41) & 17(33) & 14(30) & 12(29) \\
\bottomrule
\end{tabular}
\caption{Overall speedups w.r.t. the FE solution and, in parenthesis, speedups of the solve part of the online phase, i.e. speedup relative to the total running time excluding the setup part in which the reduced bases are loaded and the constant matrices are assembled and projected onto the reduced spaces.}
\label{table:speedup1}
\end{table}

Fig.~\ref{fig:qualitative} shows, in the first two rows, the magnitude of the velocity field and pressure distribution at times $t = 0.15$ s and $t = 0.25$ s obtained with the RB method and the corresponding point-wise errors with respect to the FE solution. The POD tolerances in every subdomain have been set to $\epsilon_u = 1\text{e--3}$ and $\epsilon_p = 1\text{e--5}$. We observe that, despite the global mesh being nonconforming, the velocity and pressure appear to be quite smooth at the interfaces. The comparison with the FE solution highlights the fact that the largest errors are committed in the region of the bifurcation. This is likely due to the fact that, as shown in Table~\ref{table:redbasis}, the reduced basis for the corresponding building block (B) is based on a smaller number of snapshots. However, the RB and the FE solutions match quite accurately overall, as the relative error is negligible in every part of the domain. The last row of Fig.~\ref{fig:qualitative} depicts the distribution of the magnitude of the WSS on the boundary of the artery in the RB and the magnitude of the error. The influence of the coupling is noticeable: indeed, it is clearly possible to spot the interfaces as regions with abnormally low or high WSS. However, this effect is not due to the RB approximation but rather to the coupling strategy: indeed, the RB and FE approximations are extremely close, as proven by the small magnitude of the error on the WSS.

Fig.~\ref{fig:erroraortaportion} shows the H1 and L2 relative errors on velocity and pressure integrated in time, defined as
\begin{equation}
e_u^2 := \dfrac{\int_{0}^\text{T} \Vert \mathbf u^{(j),h} - \mathbf u^{(j),N} \Vert^2_{\text{b},\pazocal{V}}}{\int_{0}^\text{T} \Vert \mathbf u^{(j),h} \Vert^2_{\text{b},\pazocal{V}}}, \quad e_p^2 := \dfrac{\int_{0}^\text{T} \Vert p^{(j),h} - p^{(j),N} \Vert^2_{\text{b},\pazocal{Q}}}{\int_{0}^\text{T} \Vert p^{(j),h} \Vert^2_{\text{b},\pazocal{Q}}},
\label{eq:errors}
\end{equation}
where $\Vert \mathbf u \Vert^2_{\text{b},\pazocal{V}} = \sum_{i = 1}^{N_\domain} \Vert \mathbf u^{(j)} \Vert^2_{\pazocal{V}^{(j)}}$ and $\Vert p \Vert^2_{\text{b},\pazocal{Q}} = \sum_{i = 1}^{N_\domain} \Vert p^{(j)} \Vert^2_{\pazocal{Q}^{(j)}}$ are the broken norms. The errors are plotted as functions of the velocity and pressure POD tolerances $\varepsilon_u$ and $\varepsilon_p$, highlighting the convergence of the RB solution to the FE one as the reduced basis size increases. Clearly, $\varepsilon_u$ and $\varepsilon_p$ both contribute to the errors in velocity and pressure. Indeed, for large $\varepsilon_p$, $e_u$ and $e_p$ set on a plateau as $\varepsilon_u$ decreases, indicated that the error in the pressure is dominating the error in the velocity. For each data point in Fig.~\ref{fig:erroraortaportion}, the corresponding speedup is reported in Table~\ref{table:speedup1}. The runtime of the reference FE solution---which is composed of 641502 degrees of freedom for velocity and pressure and 567 degrees of freedom for the Lagrange multipliers---is 66892 s ($\sim$ 18.5 hours). The speedups are relative to the total runtime and, in parentheses, to the part of the online phase after the initial setup (which include the loading of the reduced basis, the assembly of the constant matrices and their projection onto the reduced spaces). The motivations to consider both speedups are twofold. Firstly, in this paper we do not focus on the optimization of the assembly part of the system, which could considerably increase the total speedup; such optimization could be carried out, for example, by employing (M)DEIM, as mentioned in Section~\ref{subsec:rbonline}. Secondly, the setup part of the RB algorithm is particular to the geometry we are interested in. As a matter of fact, should we be interested in solving flow problems corresponding to different boundary conditions and/or fluid properties but on the same geometry, the setup phase can be executed only once, and for each solution of the reduced system we take advantage of the speedups relative to the only solve phase. The gain in performance is in all cases quite substantial (at least one order of magnitude with respect to the full order solution), and we observe, as expected, the trend of increasing speedup as the size of the reduced system decreases. However, the careful profiling of the simulation highlights that most of the time of the solve phase is spent in the assembly of the reduced convective term rather than in the actual solution of the reduced system. This is because, as discussed in Section~\ref{subsec:online_rbdd}, the exact assembly of the reduced convective terms entails two projections and the construction of the full order convective term.

\begin{figure}
    \tikzsetnextfilename{errors_aorta_convectiveterm}
    \centering
    \setlength
    \figureheight{\textwidth}
    \setlength
    \figurewidth{\textwidth}
    \input{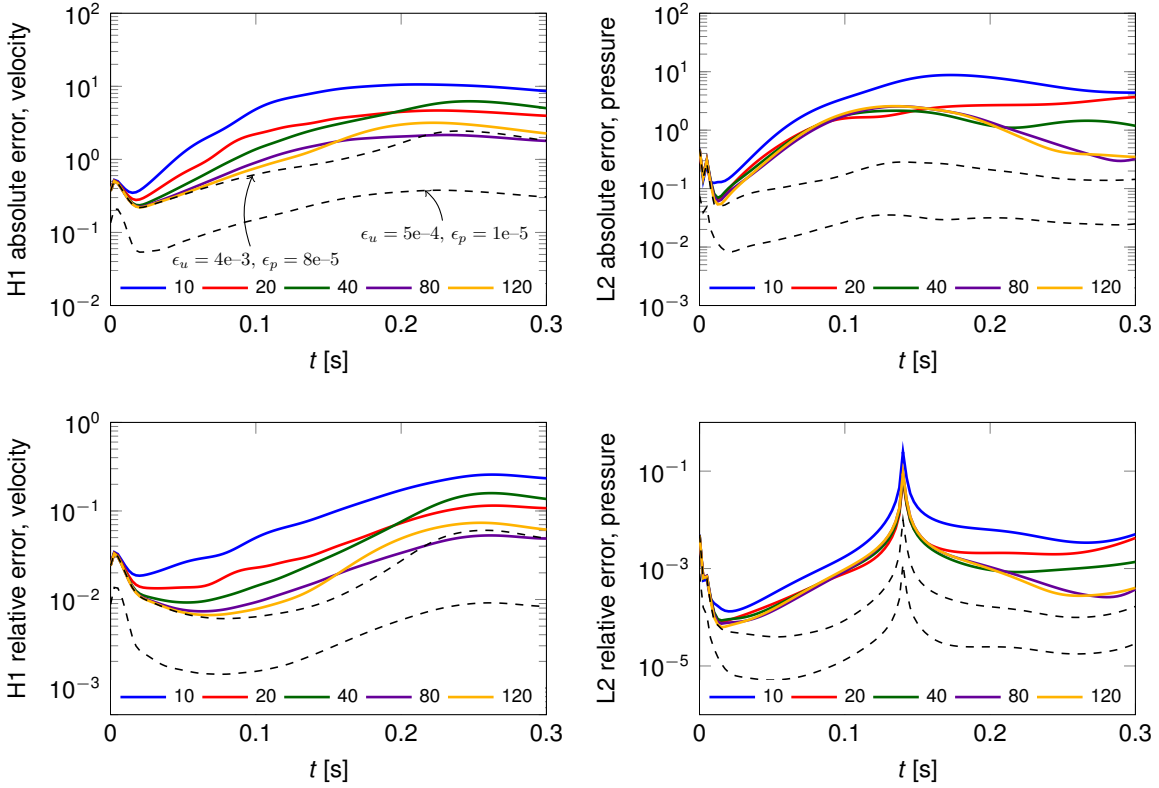}
    \caption{Errors of reduced velocity and pressure against the FE solution vs time. The colored lines refer to different choices of $N_c$ for the approximation of the nonlinear convective term. The black dashed lines show the reference errors obtained without approximation of the convective term with $\epsilon_u = 4\text{e--3}$ and $\epsilon_p = 8\text{e--5}$ (which are the same tolerances used in the simulations corresponding to the colored lines) and $\epsilon_u = 1\text{e--3}$ and $\epsilon_p = 1\text{e--5}$.}
    \label{fig:erroraortaportionappr}
\end{figure}

With the purpose of achieving higher speedups during the solution time, we consider the approximation of the convective term given in Eq.~\eqref{eq:nonlin_red_sym}. Fig.~\ref{fig:erroraortaportionappr} shows the absolute and relative H1 and L2 errors on velocity and pressure over time in function of different degrees of truncation of the convective term (i.e. different values of $N_c^{z_j}$, which we set equal to $N_c^{z_j} = N_c$ for every subdomain). The achieved speedups are, from $N_c = 10$ to $N_c =120$ and using the same notation adopted in Table~\ref{table:speedup1}, 56(998), 36(620), 22(464), 10(313), 5(215). The POD tolerances are constant and take the values $\epsilon_u = 4\text{e--3}$ and $\epsilon_p = 8\text{e--5}$.
We remark that the values of $N_c$  are to be considered in relation with the decay of the RB solutions shown in Fig.~\ref{fig:rbcoefs}, which refer to the velocity coefficients of the reference solution corresponding to $\epsilon_u = 4\text{e--3}$ and $\epsilon_p = 8\text{e--5}$.
As expected, the runtime of the solve phase is greatly decreased with respect to both the FE solution, against which the speedup achieved is always higher than 200, and with respect to the RB solution with the convective term computed as in Eq.~\eqref{eq:nlin_red}. As $N_c$ increases, the total speedup rapidly dicreases due to the quadratic dependence on that parameter of the number of integrals computed during the setup phase. Nevertheless, we believe that this strategy for approximating the convective term is of great benefit whenever it is required to run multiple simulations on the same geometry, as in this scenario the setup phase is only performed once.

\subsection{Online phase on the aorta and iliac arteries}
\label{subsec:numres_aorta}
\begin{figure}
\centering
\includegraphics[height = 0.4\textwidth,trim = {5cm 0cm 8cm 0cm},clip]{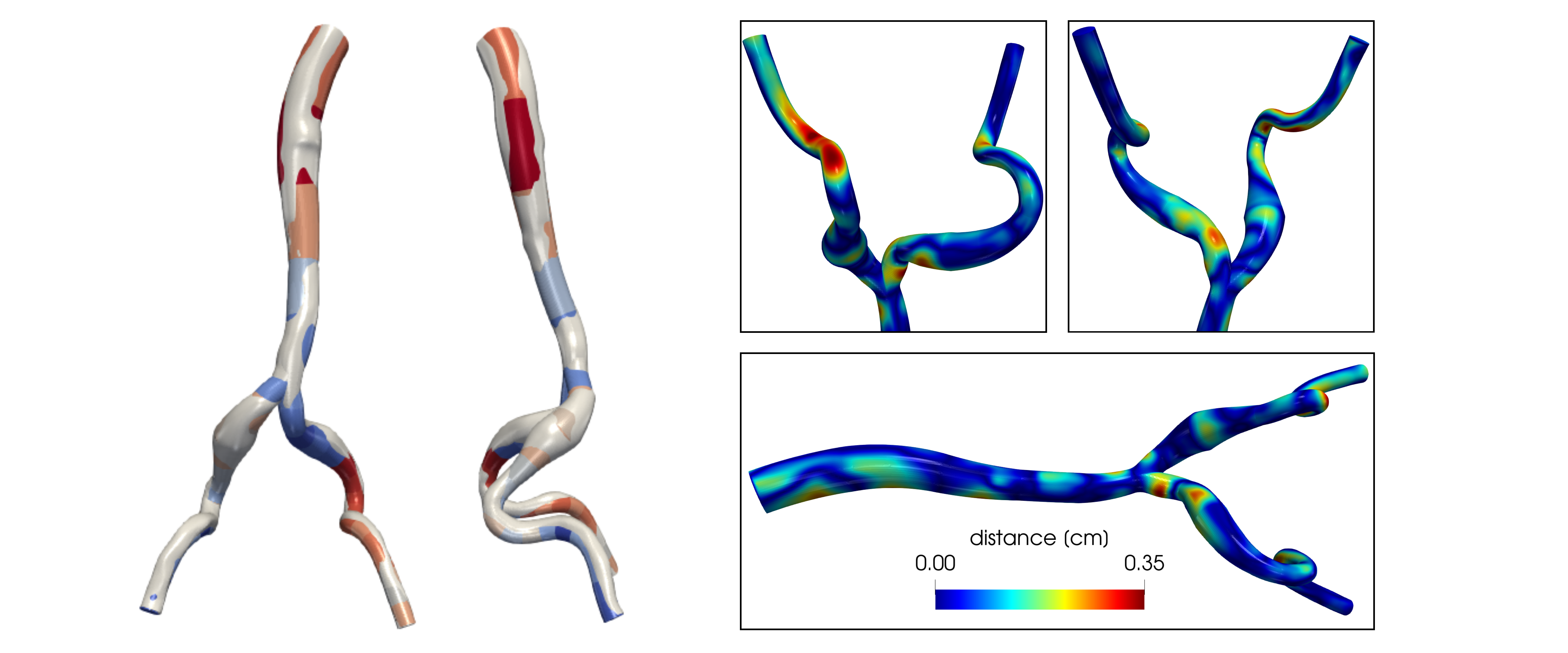}
\caption{On the left, qualitative comparison of the reference mesh with the decomposed one. On the right, quantitative estimation of the distance between the two.}
\label{fig:geometric_dif}
\end{figure}

In this section we consider a physiological geometry of an aorta with the two iliac arteries\footnote{A SimVascular tutorial based on the same geometry considered here is available on the software website (http://simvascular.github.io/docsQuickGuide.html).}. Our goal is to evaluate the effects of the geometrical approximation on the solution given by our ROM. In order to do so, we employ the geometries depicted in Fig.~\ref{fig:geometric_dif}. Specifically, on the left we show the decomposed geometry along with the ``exact'' one. On the right, we provide a quantitative analysis of the difference between the two. The algorithm to generate the decomposed geometry from the target one is out of the scope of this paper. However, the development of efficient and accurate reconstruction strategies is a topic of interest and will be addressed in future works. Importantly, we choose to employ the same set of RB basis functions computed in the offline phase described in Section~\ref{sec:reducedrb} (which, we recall, is built upon modifications of the same artificial geometry we use to test the accuracy of the ROM in the previous section). This is motivated by the perspective of employing the method in realistic scenarios which may be considerably different from the ones explored during the offline phase.

The flow problem consists of imposing the same inflow profile shown in Fig.~\ref{fig:offline} at the inlet (the aorta) and homogeneous Neumann conditions to the outlets (the iliac arteries). We take $T = 1.5\,\text{s}$ (i.e. two heartbeats) and $\Delta t = 1.25\times 10^{-3}\,\text{s}$. As already anticipated, the reference simulation is computed with the SimVascular solver svSolver. This software is based on the FE method with P1-P1 elements and VMS-SUPG stabilization; we refer to \cite{bazilevs2007variational} for more information regarding this numerical approach. It is therefore challenging to devise a \textit{fair} comparison between the ROM---which, we recall, is built upon a P2-P1 discretization---and the reference solution in terms of efficiency and accuracy. Nevertheless, we provide for sake of completeness some data regarding the reference solution. This is computed on a fine mesh (composed of 1823827 nodes) which is selected by studying the convergence of the WSS on the boundary; the simulation using SimVascular took 46457 s ($\sim$ 13 hours) by employing 28 cores.

\begin{figure}
\centering
\tikzsetnextfilename{qualitative_aorta_iliac}
\begin{tikzpicture}
    \node[inner sep=0pt] (vt15) at (0,0)
    {\includegraphics[width = 0.23\textwidth]{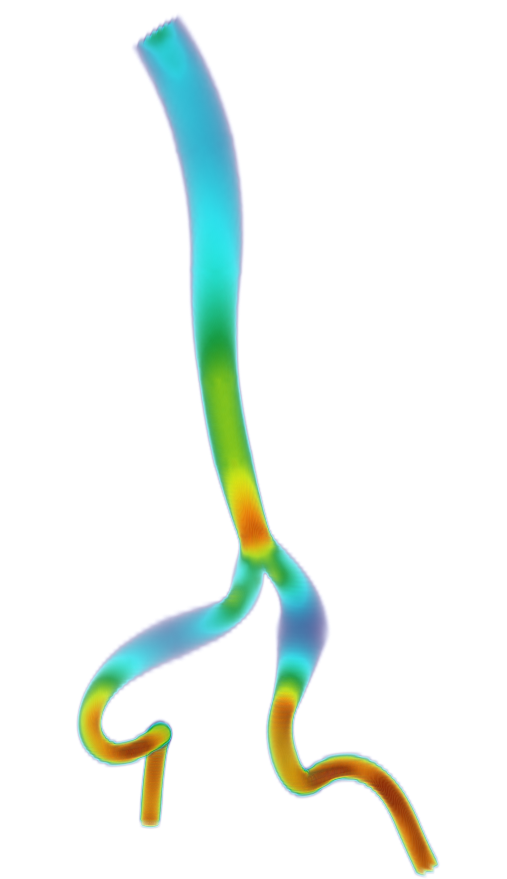}
     \includegraphics[width = 0.23\textwidth]{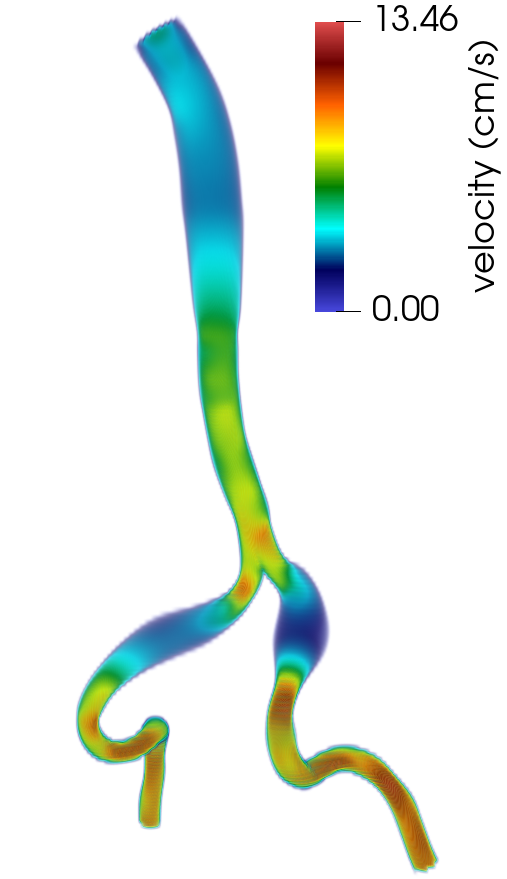}};
    \node[inner sep=0pt,right=0.6cm of vt15] (vt25)
    {\includegraphics[width = 0.23\textwidth]{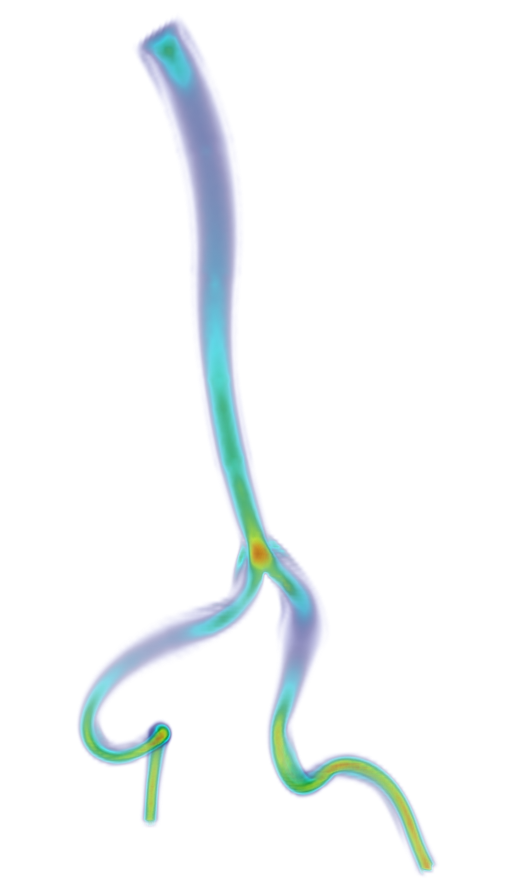}
     \includegraphics[width = 0.23\textwidth]{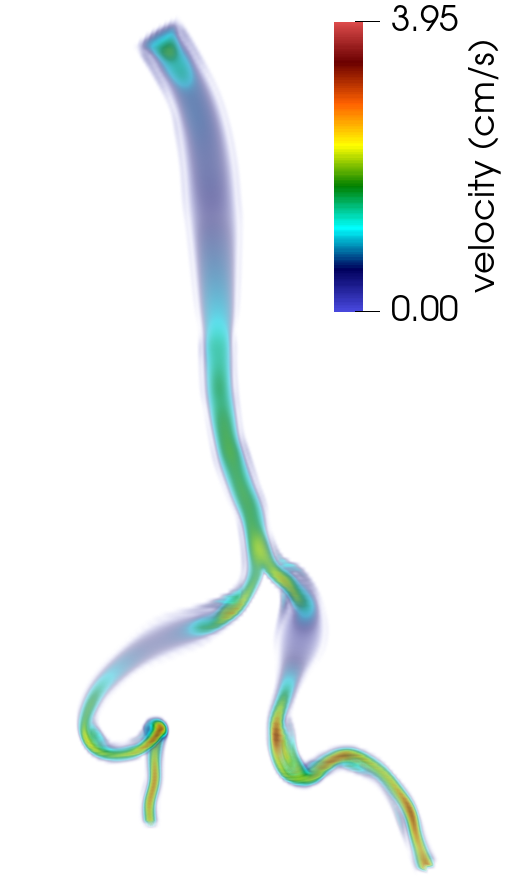}};
    \node[above=0.8cm of vt15] {$t = 0.9$ s};
    \node[above=0.1cm of vt15] {RB \qquad \qquad \quad \qquad SimVascular};
    \node[above=0.1cm of vt25] {RB \qquad \qquad \quad \qquad SimVascular};
    \node[above=0.8cm of vt25] {$t = 1.25$ s};

    \node[inner sep=0pt,below=0.6cm of vt15] (pt15)
    {\includegraphics[width = 0.23\textwidth]{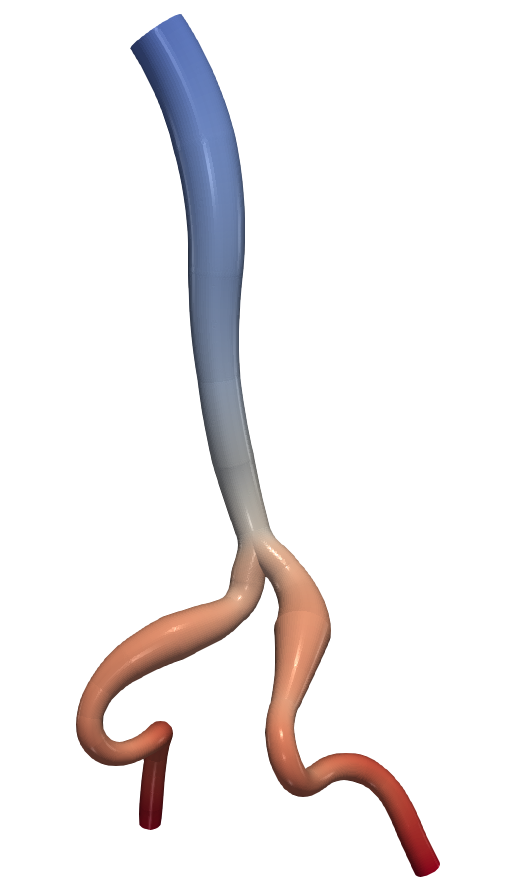}
     \includegraphics[width = 0.23\textwidth]{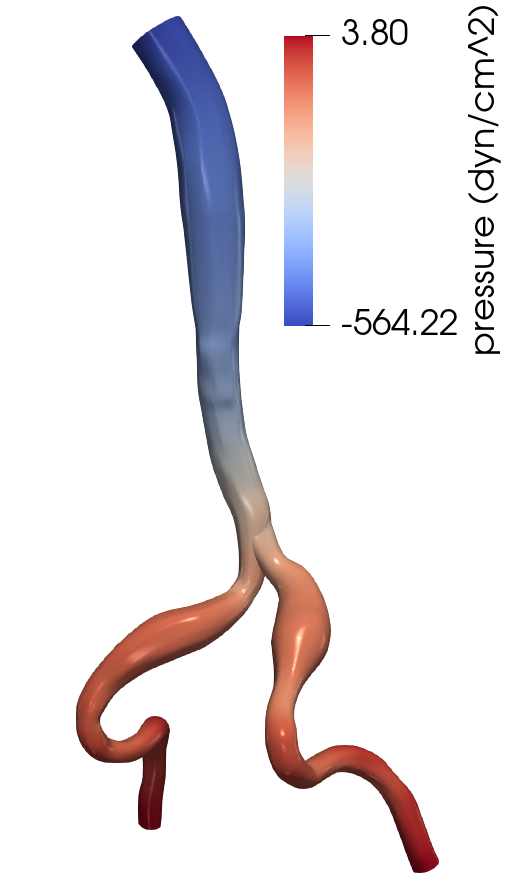}};
    \node[inner sep=0pt,right=0.6cm of pt15] (pt25)
    {\includegraphics[width = 0.23\textwidth]{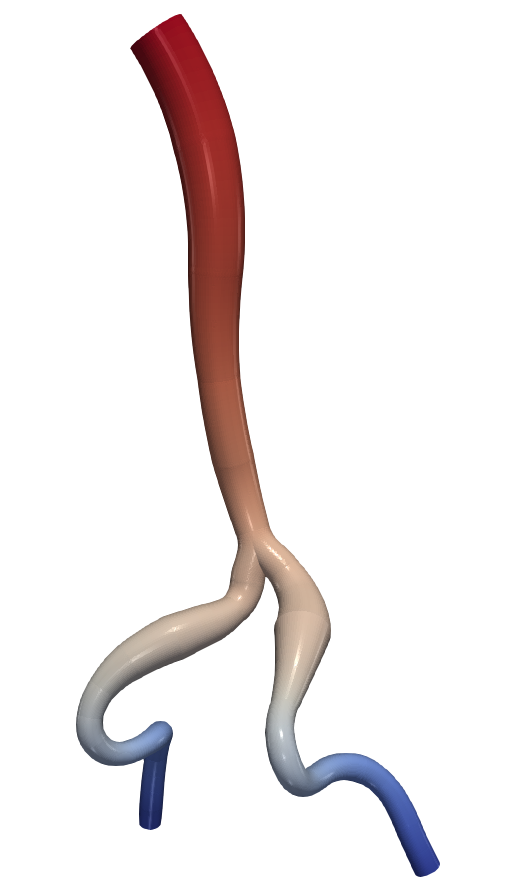}
     \includegraphics[width = 0.23\textwidth]{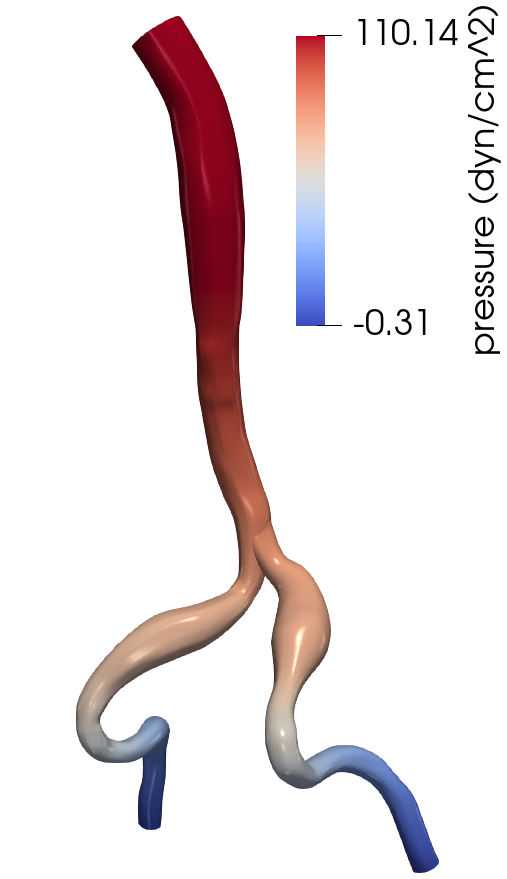}};
     \node[inner sep=0pt,below=0.6cm of pt15] (wss15)
     {\includegraphics[width = 0.23\textwidth]{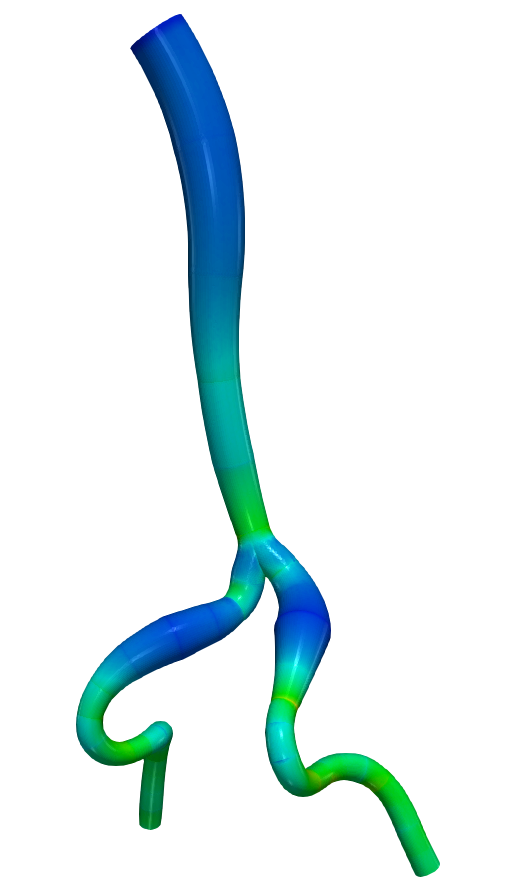}
      \includegraphics[width = 0.23\textwidth]{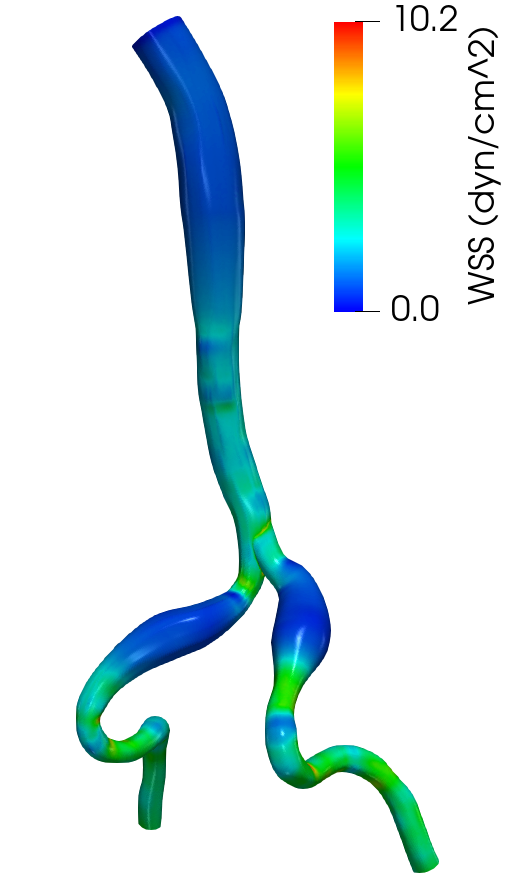}};
     \node[inner sep=0pt,right=0.6cm of wss15] (wss25)
     {\includegraphics[width = 0.23\textwidth]{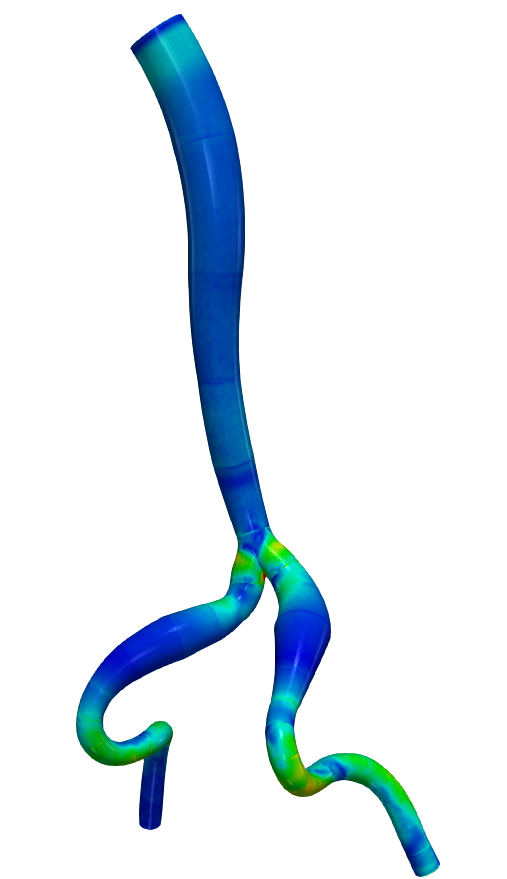}
      \includegraphics[width = 0.23\textwidth]{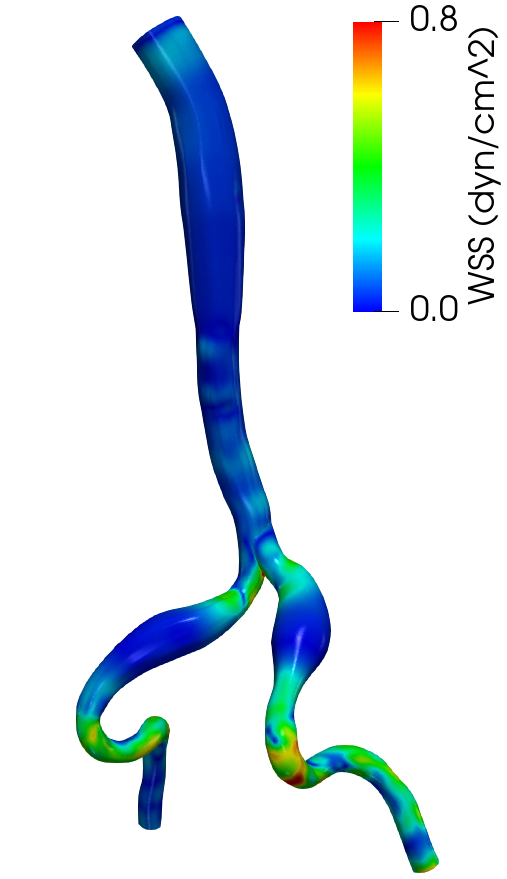}};
\end{tikzpicture}
\caption{The left and right columns---each composed of two sub-columns of plots---refer to time $t = 0.9$ s and $t = 1.25$ s, respectively. First row: velocity magnitude volume plot of the RB and reference solutions (sub-column left and sub-colum right, respectively). Second row: pressure plot of the RB and reference solutions (sub-column left and sub-colum right, respectively). Third row: magnitude of the WSS of the RB solution and reference solution (sub-column left and sub-colum right, respectively). The RB solution corresponds to the choice $\epsilon_u = 1\text{e--3}$ and $\epsilon_p = 1\text{e--5}$.}
\label{fig:qualitative_aorta_iliac}
\end{figure}

In Fig.~\ref{fig:qualitative_aorta_iliac}, we show the qualitative comparison of the velocity field magnitude, pressure and WSS distribution on the wall at two different timesteps. The RB solution is obtained with $\epsilon_u = 1\text{e--3}$ and $\epsilon_p = 1\text{e--5}$. We observe that, despite the differences in the employed geometries and in the underlying numerical discretization, the solutions share similar features. For instance, the pressure distribution is qualitatively almost identical, and the ranges for velocity and WSS magnitude achieved in every region are comparable. It is apparent, however, that most of the error (on the velocity magnitude in particular) is in the vicinity of the bifurcation. This is due to the fact that our choice of geometrical parameters for the corresponding building block does not allow for the reference bifurcation to be deformed into the target one with sufficient accuracy. For example, Fig.~\ref{fig:deformations} shows that we do not take into account the possibility of varying the radiae of the outlets, and this reflects in a large geometric error particularly on one of the branches, as depicted in Fig.~\ref{fig:geometric_dif} (bottom branch in the bottom right plot).

\begin{figure}
\tikzsetnextfilename{quantitative_aorta_iliac}
\centering
\setlength
\figureheight{0.4\textwidth}
\setlength
\figurewidth{\textwidth}
\input{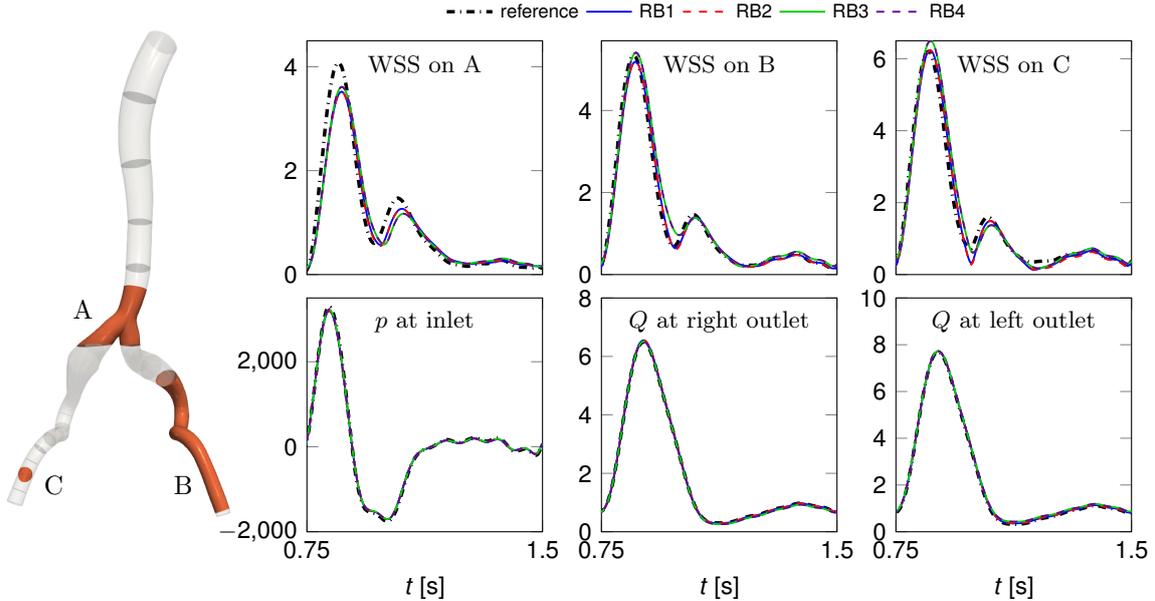}
\caption{Average WSS (in dyn/cm$^2$) on the three regions marked on the figure on the left (top row), and pressure $p$ (in dyn/cm$^2$) and flow rates $Q$ (in cc/s) at inlet and outlets, respectively (bottom row). The black dashed line refers to the reference solution computed by SimVascular, whereas the 4 colored lines are obtained with different RB settings. RB1: $\epsilon_u = 4\text{e--3}$, RB2: $\epsilon_u = 8\text{e--3}$ and $N_c = 80$, RB3: $\epsilon_u = 6.4\text{e--2}$ and $N_c = 40$, RB4: $\epsilon_u = 6.4\text{e--2}$ and $N_c = 20$. In all cases, $\epsilon_p = 1\text{e--5}$.}
\label{fig:quantitative_aorta_iliac}
\end{figure}

A more quantitative analysis of the performance of the ROM with respect to the reference solution is presented in Fig.~\ref{fig:quantitative_aorta_iliac}. Here, we show the average of the WSS magnitude over the three regions highlighted in the figure on the left, the pressure at the inlet and the outflow rate at the outlets, for the reference solution and for RB solutions corresponding to different choices of tolerances and truncations for the approximation of the nonlinear term. We chose to focus on a ``fine'' RB solution (RB1), where we do not apply the truncation of the convective term, and more ``coarse'' but efficient RB solutions with approximation of the convective term (RB2, RB3, RB4); for details regarding the employed POD tolerances and number of terms in the truncated sum, we refer the reader to the caption of Fig.~\ref{fig:quantitative_aorta_iliac}. The setup, solve and total runtimes in seconds for these simulations---which are run, differently from the reference simulation, on a single core---are the following: for RB1, $4240 + 36170  = 40410$ ($\sim$ 11 hours, speedup of 1.1), for RB2, $19382 + 4760 = 24142$ ($\sim$ 6.5 hours, speedup of 2), for RB2,  $7456 + 2398 = 9854$ ($\sim$ 2.5 hours, speedup of 5), for RB3, $1816 + 3774 = 5590$ ($\sim$ 1.5 hours, speedup of 8). In all cases, we achieve speedups larger than one with respect to the reference simulation obtained with SimVascular (although the gain is negligible in the case of RB1) but on a single core instead of 28. From the results presented in Fig.~\ref{fig:quantitative_aorta_iliac}, we note that, while the approximation of the pressure and flow rate is extremely precise for all RB settings, the performance on the WSS is more challenging. The curves for the average WSS on the two regions on the iliac arteries (B and C) are quite close to the reference one compared to the average WSS on the bifurcation (A). However, this is likely an effect of the geometric approximation rather than the accuracy of the ROM per se. As a matter of fact, we already noted in Fig.~\ref{fig:qualitative_aorta_iliac} that the largest errors are located in that area. We also remark that the POD tolerance plays a more dramatic role in the quality of the solution than the number of terms retained in the truncated nonlinear term $N_c$. Indeed, the simulations with smallest tolerances (RB1 and RB2) and largest ones (RB3 and RB4) lead to similar results, regardless of the value of $N_c$. Nevertheless, truncating the convective term is beneficial to the efficiency of the ROM.

\section{Conclusions}
\label{sec:conclusions}
In this paper, we presented an implementation of the reduced basis element method for the solution of the unsteady 3D Navier--Stokes equations in the context of cardiovascular simulations. We first considered the problem of coupling finite element solutions defined on subdomains obtained from parametrized geometrical deformations of reference building blocks. This was necessary, as the offline phase of our reduced order method requires the generation of snapshots from coupled finite element flow solutions obtained on a variety of geometries. In order to improve the efficiency of the coupled finite element solver, we devised an ad-hoc preconditioner which takes advantage of the saddle-point structure of the discretized linear system. In the following parts of the paper, we formulated the reduced order model by projecting the matrices and variables (velocity and pressure) onto the reduced basis spaces. This procedure is beneficial because it allows us to considerably reduce the number of degrees of freedom (hence, the size of the linear system to be solved at each iteration of the Newton--Raphson algorithm). In the numerical simulations, we demonstrated the capabilities of the method on the same geometry used for the offline phase and on a physiological geometry consisting of an aorta with the two iliac arteries. In the first case, we registered considerable speedups (from 12 to 33 over the total runtime and from 29 to 50 over the sole solve phase) with respect to the full order solution. Considerable gain in performance was also achieved in the second case for some choices of POD tolerances, although the fact that we considered a reference solution obtained with a different solver (i.e. SimVascular) made the comparison in terms of runtime more complex. In both applications we also analyzed the performance in terms of wall-shear stress reconstruction, which is possible in our reduced order method---as opposed, for example, to geometrical multiscale methods---because the 3D nature of the flow problem is preserved.

We believe that the results presented in this work are promising and that our study suggests many possibilities for the future developments of this reduced order method. As we discuss more in depth in the relative sections of this paper, these include: i) a parallel implementation of the devised saddle-point preconditioner which could exploit the special structure of the linear system (see Remark~\ref{remark:schur_parallel}), ii) an offline strategy based on the solution of physiological flow problems (rather than flow problems defined on artificially deformed geometries), iii) efficient ways to reduce the complexity of the setup phase and the assembly of the reduced convective term, and iv) automatic algorithms for the generation of accurate decomposed geometries out of medical images or reference meshes. Regarding this last point, we also think that---due to the geometrical difficulties to map reference bifurcations into target ones (see discussion in Section~\ref{subsec:numres_aorta})---an interesting follow-up of the current work is the study of an hybrid method in which some parts of the target geometry (e.g. the bifurcations) are modeled by means of the finite element method and others by means of the reduced basis method. This approach would allow us to easily treat cases featuring more complex geometries---for instance, cerebral aneurysms---that are not approximated as trivial deformations of tubes and which are currently too challenging to tackle with a strategy purely based on the reduced basis method.
\section*{Acknowledgments}
The research of Luca Pegolotti and Simone Deparis was supported by the Swiss National Foundation (SNF), project No. 188031. The work of Alison Marsden and Martin Pfaller was supported by NIH grant 1R01LM013120. The authors gratefully aknowledge the Scientific IT and Application Support at EPFL for providing the computational resources necessary to the numerical simulations presented in this paper.

\bibliographystyle{abbrv}
\bibliography{references}

\end{document}